\documentclass[12pt,twoside,english]{article}
\usepackage{ifthen}                % if then else Abfragen
\newboolean{boldeqnitalic}
\setboolean{boldeqnitalic}{true} 

\usepackage{caption}
\usepackage[cp1252]{inputenc}
\usepackage[intlimits] {amsmath}   % muss for defjs# stehen 
\usepackage{defjs2}                % eigene Definitionen    
\usepackage{epsfig}
\usepackage{graphicx}
\usepackage{ifthen}
\usepackage{fancyhdr}
\usepackage{rotating}
\usepackage{multirow}
\usepackage{booktabs}              % Dicke Tabellenlinien
\usepackage{amssymb}
\usepackage{amsbsy}
\usepackage{bm}                    % fette Schrift in Formeln
\usepackage{babel}
\usepackage{theorem}
\usepackage{upgreek}  % gerade (roman) griechische Buchstaben z.B. \upalpha
\usepackage{pifont}   % zusaetzliche Schriften, z.B. eingekreiste Zahlen
\usepackage{mathrsfs}
\usepackage{datetime}
\usepackage{tikz}
\usetikzlibrary{matrix}
\usepackage[all]{xy}
\usepackage{rotating}
\usepackage{subfig}
\usepackage[labelfont=bf, justification=justified]{caption}

\usepackage{algpseudocode}
\usepackage{algorithm}

\definecolor{dunkelgrau}{rgb}{0.8,0.8,0.8}
\definecolor{hellgrau}{rgb}{0.90,0.90,0.90} %... slightly darker 
\newcommand\defi{\mathrel{\overset{\makebox[0pt]{\mbox{\normalfont\scriptsize\sffamily def}}}{=}}}
\fancypagestyle{plain}{%
  \fancyhead{}%
  \fancyfoot[c]{\sffamily\thepage}%
}
\makeatletter                           % Leere Fuellseiten
\def\cleardoublepage{\clearpage\if@twoside \ifodd\c@page\else
  \hbox{}
  \vspace*{\fill}
  \thispagestyle{empty}
  \newpage
  \if@twocolumn\hbox{}\newpage\fi\fi\fi}
\makeatother
\begin{document}
\unitlength1.0cm
\frenchspacing

\thispagestyle{empty}
 
\begin{center}    
	{\bf \large From image data towards microstructure information --}\\[2mm]
	{\bf \large accuracy analysis at the digital core of materials}  \\[2mm]
%	{\bf \large for reliability of the digital twin in materials genomics}
\end{center}      

\vspace{4mm}
\ce{Bernhard Eidel, Andreas Fischer, Ajinkya Gote}

\vspace{4mm}

\ce{\small Heisenberg-Group, Institute of Mechanics, Department Mechanical Engineering} 
\ce{\small University Siegen, 57068 Siegen, Paul-Bonatz-Str. 9-11, Germany} 
\ce{\small $^{\ast}$e-mail: bernhard.eidel@uni-siegen.de, phone: +49 271 740 2224, fax: +49 271 740 2436} 
\vspace{2mm}

%\bigskip
%
%\begin{center}
%{\bf \large Highlights}
%
%\bigskip
%
%{\footnotesize
%\begin{minipage}{15.5cm} 
%
%\vspace*{-2mm}
%
%\begin{itemize}
% \item Concept for error analysis of image-based microstructure representation.  \\[-6mm]
% \item Error-decomposition into finite resolution modeling error and discretization error. \\[-6mm]
% \item Embedding into the unified error framework of two-scale FEM for homogenization. \\[-6mm]
% \item Concept enables a quantitative accuracy-efficiency balance score. \\[-6mm] 
% \item Best practice for resolution coarsening along with adaptive mesh coarsening.  
%\end{itemize}
%
%\end{minipage}
%}
%\end{center}

\bigskip

%\vfill
%\end{document}

\begin{center}
{\bf \large Abstract}

\bigskip

{\footnotesize
\begin{minipage}{14.5cm} 
\noindent
A cornerstone of computational solid mechanics in the context of digital transformation are databases for microstructures obtained from advanced tomography techniques.
Uniform discretizations of pixelized images in 2D 
%and reconstructed voxel structures in 3D 
are the raw-data point of departure for simulation analyses. 
This paper proposes the concept of a unified error analysis for image-based microstructure representations in uniform resolution along with adaptively coarsened discretizations. The analysis distinguishes between a modeling error due to finite, possibly coarsened image resolution and a discretization error, investigates their quantitative relation, spatial distributions and their impacts on the simulation results both on the microscale and the macroscale in the context of computational homogenization.    
The assessment of accuracy and efficiency is carried out for an exemplary two-phase material. Beyond the example considered here the concept is a rational tool in the transformation of raw image data into microstructure information adapted to particular simulation needs and endows the digital twin of real microstructures with validated characteristics for reliable, predictive simulations.    
\end{minipage}
}
\end{center}

{\bf Keywords:}
Microstructures; Image-based modeling; Error analysis; Homogenization; Digital twin \hfill 
%\\
%{\hfill Draft 01 --\,\today \, at \currenttime} 

%\hspace*{6.3cm} 
%\\[-2mm]
%vers.\,\today \, at \currenttime\\
%\begin{center}\emph{03.12.2020, accepted for publication in \emph{ZAMM}\\
%The finally published version can show minor differences to the present preprint}
%\end{center}
%{\color{red}{\bf Agenda:} micro-, macro-: one word or two words? $\mathbb{A}^0$ or what -- look to the first version.} 

%----------------------------------------------------------------------------------------------------------- 

\section{Introduction}
\label{sec:intro}

For image-based microstructure representations of heterogeneous solids, in 2D by pixels (px) and by voxels in 3D, the question arises in the numerical analysis, how image resolution influences the simulation results. A finite image resolution introduces a modeling error which comes on top of the standard discretization error. The computation of the modeling and discretization errors, their magnitudes, spatial distributions, their impact on microscale and macroscale results and corresponding simulation costs in computational homogenization is the topic of this paper. 

For pixel-/voxel-type microstructure representations the image resolution frequently defines the uniform finite element discretization likewise or it serves as the point of departure for e.g. quadtree-/octree-type
 adaptive coarsening. Efficient FFT-based homogenization methods directly use for the solution of the microscale problem uniform discretizations, compare Moulinec-Suquet \cite{Moulinec.1998} and others. The alternative of reconstructing smooth phase boundaries as a prerequisite for conforming finite element discretizations is time-consuming and cumbersome. Moreover, a decoupling of discretization from material properties as realized in the Finite Cell Method (FCM) \cite{Parvizian.2007}, \cite{Duster.2008}, \cite{Schillinger.2012}, \cite{Duster.2012}, \cite{Schillinger.2015} is hardly available by default in commercial and research solvers. The above aspects underpin why pixel- and voxel-based discretizations are of utmost relevance. 

\begin{Figure}[htbp]
	\begin{minipage}{16.5cm}  
	    \subfloat[uniform, 1024 px]
		{\includegraphics[height=5.2cm,angle=0]{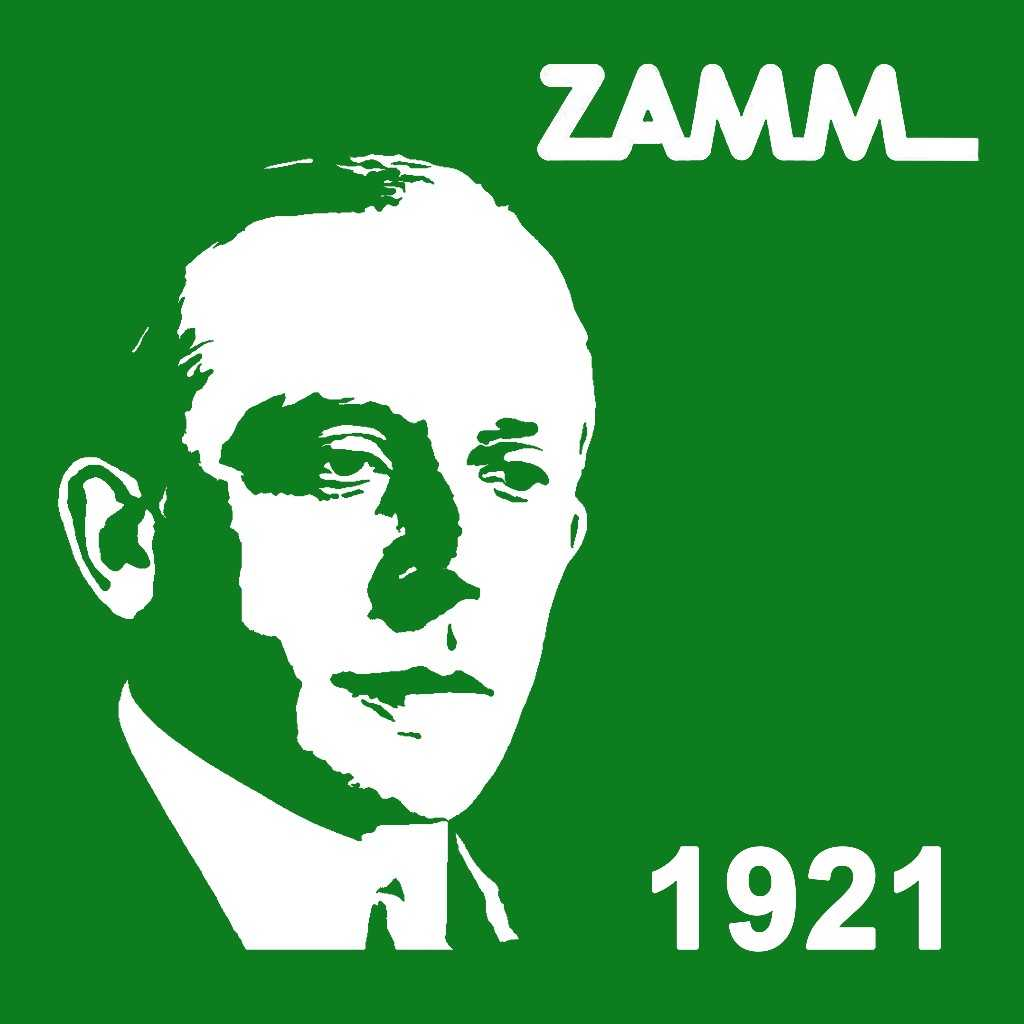}}
		\hspace*{0.005\linewidth}
		\subfloat[adaptive]
		{\includegraphics[height=5.2cm, angle=0]{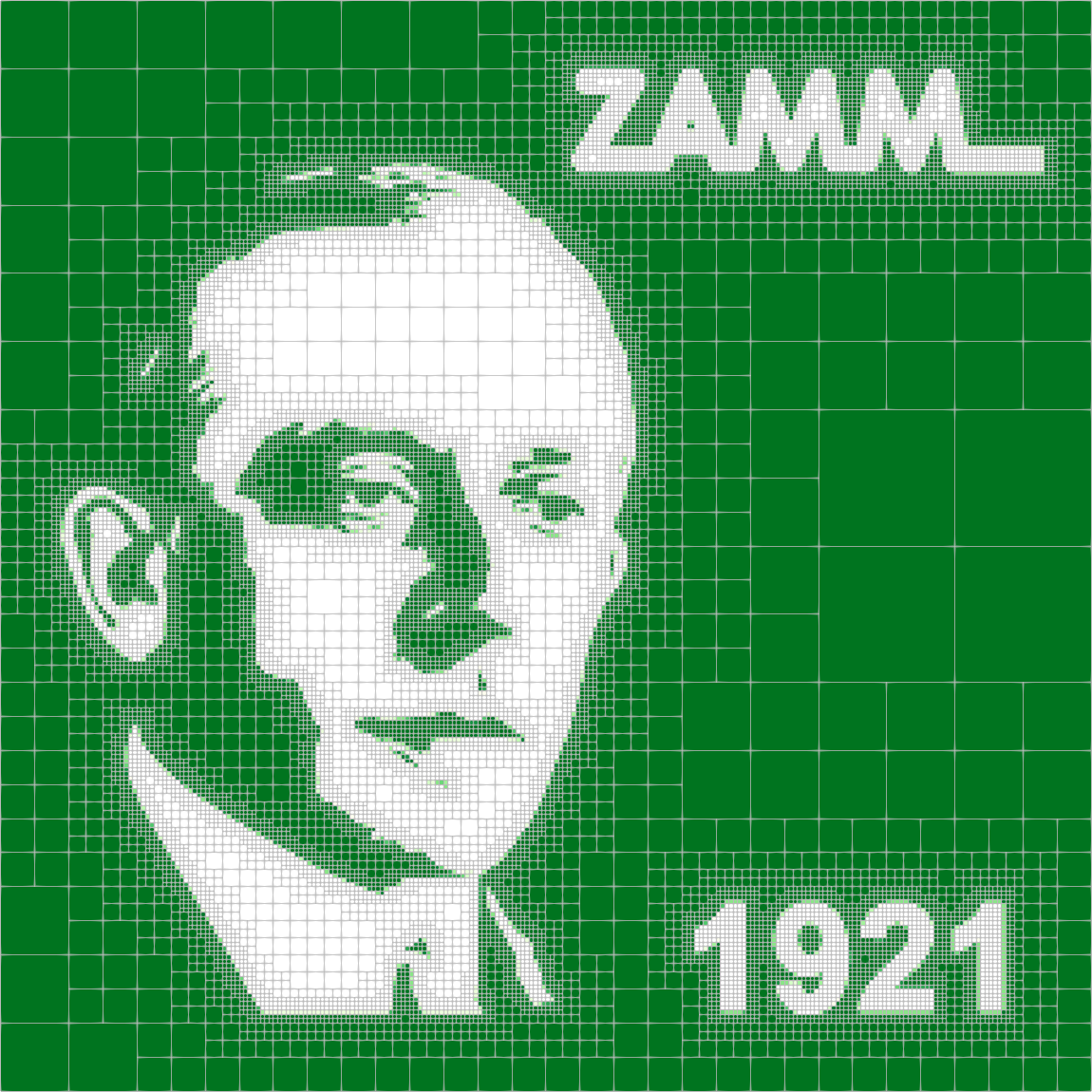}}
		\hspace*{0.005\linewidth}
		\subfloat[uniform, 64 px]
		{\includegraphics[height=5.2cm, angle=0]{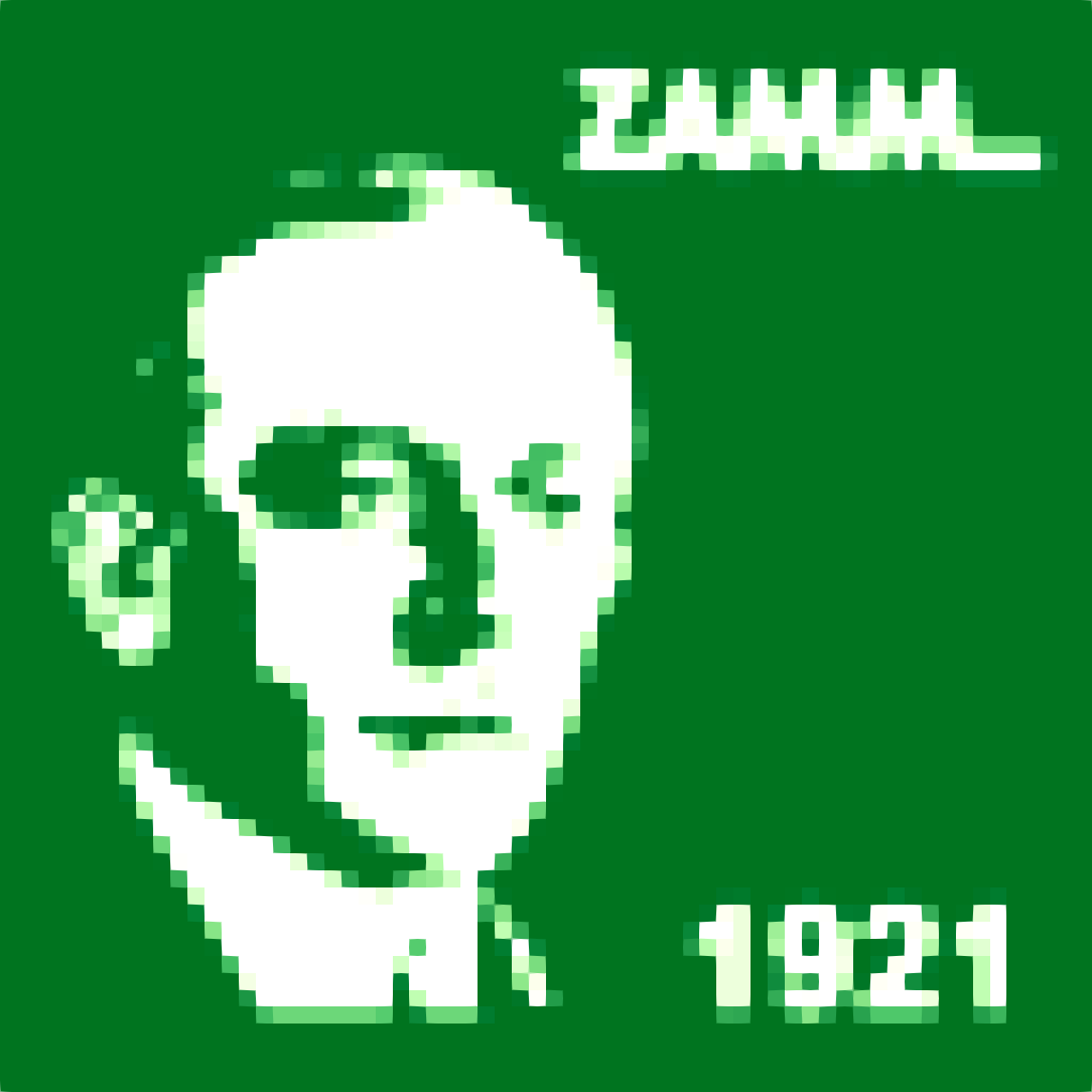}}
	\end{minipage}
	\caption{{\bf Richard-von-Mises-ZAMM-1921-microstructure}: (a) uniform resolution of 1024 px per edge, (b) adaptive mesh of (a), (c) uniform resolution of 64 px per edge.  
    \label{fig:Types-of-coarse-graining}}
\end{Figure}

Figure~\ref{fig:Types-of-coarse-graining} displays for the 2D case of a two-phase microstructure\footnote{Richard von Mises (1883--1953), an Austrian/US-American mathematician, founded the Zeitschrift f\"ur Angewandte Mathematik und Mechanik (ZAMM) in 1921.} two different versions of coarsening the original high-resolution image in (a); in (c) a uniform pixel-coarsening is carried out, in (b) an adaptive, microstructure-informed mesh coarsening. 

The question, which (uniform) image resolution shall be chosen, possibly as the starting point for a consecutive adaptive mesh coarsening, can be answered based on (i) an error analysis which distinguishes between a modeling error and a discretization error and (ii) corresponding computational efforts. The outcome enables the selection of a discretization with a well-balanced accuracy-efficiency score.

Work which considers the impact of image resolution on simulation results is rare and throughout restricted to effective properties. Hutzenlaub et al. \cite{Hutzenlaub.2013} investigate for particular cathode catalyst layers how pixel coarsening influences parameters such as diffusivity and conductivity. The work of Nguyen et al. \cite{Nguyen.2015} considers for foamed concrete the impact of image resolution to global stress-strain curves, local damage initiation and evolution. For the effect of voxel-resolution on flow in porous media see Sha et al. \cite{Shah.2016}, for applications in digital rock physics see Berg et al. \cite{Berg.2017}.
\\
For the finite cell method (FCM) the first mathematical error analysis was provided by Dauge et al. \cite{Dauge.2015}. Error estimation for FCM has very recently been introduced by Di Stolfo et al. \cite{DiStolfo.2019}, \cite{DiStolfo.2019b}. In the FCM context the quadrature error is very similar to the modeling error in the present work.

%{\color{red} @BE: further references (...) FFT-Zeman}

We assume that in the present concept of coarsening resolution and discretization the selected microstructure specimen is representative with respect to phase fractions, morphology etc., such that resolution coarsening preserves its representative characteristics approximately. Doing so, the present work follows another route than approaches that replace the representative volume element RVE by a simplified microstructure of reduced morphological and therefore reduced computational complexity. Serious attempts have been made to identify suchlike surrogate microstructures which shall qualify as candidates by their degree of statistical similarity \cite{Povirk.1995}, \cite{Ohser.2006}, \cite{Kumar.2006}, \cite{Schroder.2011}, \cite{Balzani.2014}, \cite{Scheunemann.2015}. %\cite{LeBris.2016}.

%In view of the computational complexity of solving statistical representative volume element reflecting the real microstructure, alternative approaches have tried to replace the RVE by a surrogate micro problem \cite{Povirk.1995}, \cite{Ohser.2006},  \cite{Balzani.2014}, \cite{Scheunemann.2015}. The similarity of the artificial surrogate microstructure, sometimes referred to as SSRVE, was measured by statistical measures such as Minkowski functional and alike.
%The present work chooses a different and novel route. It maintains the unit cell qualified as an RVE and searches for coarse-grained representations while preserving the original morphology of the microstructure.
 
%\begin{Figure}[htbp]
%	\subfloat[1024 px: 2 phases] 
%	{\includegraphics[width=7.2cm, angle=0]{1024p_eye}} 
%	\hspace*{0.04\linewidth}
%	\subfloat[256 px: 17 phases]
%	{\includegraphics[width=7.2cm, angle=0]{256p_eye}}
%	%
%	\caption{{\bf Pixel coarsening:} Eye detailed with (a) 1024 px resolution along with two phases/ sharp interfaces, (b) 256 px resolution with 17 interphases/blurred interfaces. 
	%{\color{red} von-Mises' Ohr waere wohl unproblematischer? Das schaut einen nicht so an.}
 %   \label{fig:Eye-details-vMises}}
%\end{Figure} 

Image segmentation refers to the identification of disjoint phases with sharp boundaries, where the number and properties of phases are assumed to be known, for a full-fledged example see \cite{Andra.2013}, \cite{Andra.2013b}, for an overview \cite{Ohser.2006}. The outcome of a standard segmentation looks like the image in Fig.~\ref{fig:Types-of-coarse-graining} (a). The present work considers in uniform pixel coarsening this standard segmentation and additionally an alternative approach which is visualized in Fig.~\ref{fig:Types-of-coarse-graining} (c); novel interphases are created at interfaces with a color code i.e. stiffness reflecting the phase fractions of the contributing pixels being merged in the coarser pixel. As a consequence, this approach, which is the standard of e.g. images in bitmap format, preserves the overall phase fraction of the microstructure but reduces the stiffness contrast at interfaces. How does that affect stress in terms of extrema and jumps at material interfaces, and corresponding errors? These questions will be answered in Sec. \ref{sec:NumericalExamples}.
  
Since homogenization in terms of a two-scale finite element method aims at describing processes on the microscale as well as on the macroscale, the accuracy-efficiency balance-score shall be considered on both scales. On the macroscale the micro error arrives as a propagated error and adds to the macro discretization error. We anticipate that this can result in a microstructure discretization which is rejected at the microscale for its large discretization errors, but still favorable for its excellent accuracy on the macroscale, for an example see \cite{Fischer.2020}.
 
The paper is organized as follows; to set the stage Sec.~\ref{sec:FEHMM} provides a short outline of the finite element heterogeneous multiscale method (FE-HMM) for linear elasticity and introduces quantities used in the augmented FE-HMM framework of errors and estimates described in Sec.~\ref{sec:ErrorsAndEstimates}. Section \ref{sec:Mesh-coarsening} introduces in detail the above-mentioned different concepts of uniform resolution-coarsening along with  consecutive adaptive mesh coarsening. Section \ref {sec:error_estimation_average} describes reconstruction-type error estimation for interfaces that have undergone resolution coarsening. Section \ref{sec:NumericalExamples} applies the error analysis to the von-Mises-ZAMM microstructure in Fig.~\ref{fig:Types-of-coarse-graining} (a) along with the coarsening features illustrated in (b) and (c).
  
%-------------------------------------------------------------
\section{The finite element heterogeneous multiscale method }
\label{sec:FEHMM}
%------------------------------------------------------------- 

\subsection{Model problem of linear elasticity} 
\label{subsec:ModelProblemLinearElasticity}

We consider a body $\mathcal{B}$, a bounded subset of $\mathbb{R}^{n_{dim}}$, $n_{dim}=2,3$, with boundary 
$\partial \mathcal{B} = \partial \mathcal{B}_D \cup \partial \mathcal{B}_N$ where the Dirichlet boundary $\partial\mathcal{B}_D$ and the Neumann boundary $\partial\mathcal{B}_N$ are disjoint sets.
The closure of the body $\mathcal{B}$ is denoted by $\overline{\mathcal{B}}$.
The body shall be in static equlibrium. It shall exhibit an inhomogeneous composition referred to as microstructure and shall be subject to body forces $\bm f$ and surface tractions $\bar{\bm t}$.

\subsubsection{The microproblem}
\label{subsec:StrongFormLinearElasticity}

The displacement $\bm u^{\epsilon} = (u_1^{\epsilon}, \ldots, u_{n_{dim}}^{\epsilon})$ of the body is given by the solution of 
\begin{equation}
\label{eq:StrongFormMicro-2}
- \, \sigma^{\epsilon}_{ij,j} = f_i \quad \mbox{in} \, \, \mathcal{B} \, , \quad
u_i^{\epsilon} = \bar{u}_i \quad \mbox{on} \, \, \partial \mathcal{B}_{D} \, , \quad
\sigma^{\epsilon}_{ij} \, n_j = \bar{t}_i  \quad \mbox{on} \, \, \partial \mathcal{B}_{N} \, .
\end{equation}
Linear elasticity is assumed to hold $\sigma^{\epsilon}_{ij}=\mathbb{A}^{\epsilon}_{ijlm}\, {\color{black}\varepsilon_{lm}}$ with the fourth order elasticity tensor $\mathbb{A}^{\epsilon}_{ijlm}$ and the infinitesimal strain tensor $\varepsilon_{ij}$; it holds $\varepsilon_{ij}(\bm u^{\epsilon}) = 1/2 \left(u_{i,j}^{\epsilon} + u_{j,i}^{\epsilon}\right)$ or in compact notation $\bm \varepsilon (\bm u^{\epsilon}) = \bm L \, \bm u^{\epsilon}$ with the linear differential operator $\bm L$.
Superscript $\epsilon$ refers to the heterogeneity of the elastic material, in particular to the characteristic length scale of the microstructure. 
In \eqref{eq:StrongFormMicro-2}$_{3}$, $\bm n=(n_1, \ldots, n_{n_{dim}})^T$ is the unit outward normal to $\partial \mathcal{B}$.  

The variational form is obtained by multiplying the strong form \eqref{eq:StrongFormMicro-2} by a test function $\bm v \in \mathcal{V}$ and by the application of Green's formula.

Find $\bm u^{\epsilon}$ such that 
\begin{equation}
\label{eq:weak-form-for-u-epsilon}
B_{\epsilon} (\bm u^{\epsilon}, \bm v) \defi 
\int_{\mathcal{B}} {\color{black}\bm \sigma^{\epsilon}(\bm u^{\epsilon})} : \bm \varepsilon(\bm v) \, dV 
= \int_{\mathcal{B}} \bm f \cdot \bm v \, dV \, + \,  \int_{\partial \mathcal{B}_N} \bar{\bm t} \cdot \bm v \,dA 
\defi \bm F(\bm v)\, ,
\end{equation}
which must hold for all $\bm v \in \mathcal{V}$, where $\mathcal{V}$ is the space of virtual displacements fulfilling the homogeneous Dirichlet boundary conditions (BCs) $\mathcal{V}=\{\bm v; \bm v \in H^1(\mathcal{B})^{n_{dim}}, \bm v|_{\partial \mathcal{B}_D} = \bm 0 \}$.
%\begin{equation}
%\label{eq:HilbertSpaceV}
%  \mathcal{V}=\{\bm v; \bm v \in H^1(\mathcal{B})^{n_{dim}}, \bm v|_{\partial \mathcal{B}_D} = \bm 0 \} \, .
%\end{equation}

%---------------------------------------------------------------------------------------
\subsubsection{The macroproblem}
\label{subsec:Variational-FE-HMM-macro}

The strong form of the macroscopic/homogenized boundary value problem (BVP) is given by, cf. \cite{Eidel.2018} 
\begin{equation}
\label{eq:Homogenized-Strong-Form}
- \, \sigma^{0}_{ij,j} = \langle f_i \rangle \ \,  \mbox{in} \, \, \mathcal{B}  
\, , \quad 
u_i^{0} = \langle \bar{u}_i\rangle_{\Gamma} \ \, \mbox{on} \, \, \partial \mathcal{B}_{D}  
\, , \quad
\sigma^{0}_{ij} \, n_j = \langle \bar{t}_i\rangle_{\Gamma}  \ \, \mbox{on} \, \, \partial \mathcal{B}_{N}  
\end{equation}
where $u_{i}^0$ is the macroscopic displacement and $\mathbb{A}^{0}$ is the homogenized elasticity tensor. In \eqref{eq:Homogenized-Strong-Form}$_{1,3}$ $\sigma^{0}_{ij}$ is the macroscopic stress obtained by the volume average over the microdomain. 

The values for the Dirichlet as well as Neumann BC in \eqref{eq:Homogenized-Strong-Form}$_{2,3}$ are obtained by surface averages of corresponding BC in \eqref{eq:StrongFormMicro-2}$_{2,3}$, for details we refer to \cite{Eidel.2018}. Similarly, $\langle f_i \rangle$ is the volume average of body forces in \eqref{eq:StrongFormMicro-2}$_{1}$.

The solution of the homogenized problem is obtained from the variational form 
\begin{equation}
B_0 (\bm u^0, \bm v) =  \int_{\mathcal{B}} {\color{black} \bm \sigma^0(\bm u^0)}: \bm \varepsilon(\bm v) \, dV   
=  \int_{\mathcal{B}} \bm f \cdot \bm v \, dV \, + \,  \int_{\partial \mathcal{B}_N} \bar{\bm t} \cdot \bm v \,  dA
\qquad \forall \, \bm v \in \mathcal V \, ,
\label{eq:VariationalFormHomogenizedProblem}                         
\end{equation}
which follows from multiplying the strong form \eqref{eq:Homogenized-Strong-Form} by test functions $\bm v$ along with the application of Green's formula. For notational convenience we skip in \eqref{eq:VariationalFormHomogenizedProblem} and in the rest of the paper the averaging symbols $\langle \bullet \rangle$, $\langle \bullet \rangle_{\Gamma}$ for $\bm f$, $\bar{\bm u}$ and $\bar{\bm t}$ but keep in mind that these quantities follow from volume and surface averages, respectively. 

%We consider the piecewise linear continuous FEM in macro- and microspace, respectively. 
We define a macro finite element space as 
\begin{equation}
\mathcal{S}^p_{\partial \mathcal{B}_D}(\mathcal{B}, {\mathcal T}_H) = \left\{ \bm u^H \in H^1(\mathcal{B})^{n_{dim}}; \bm u^H|_{\partial \mathcal{B}_D} = \bar{\bm u}; \bm u^H|_{K} \in {\mathcal{P}}^{p}(K)^{n_{dim}}, \, \forall \, K \in {\cal T}_{H} \right\} \, ,
\label{eq:MacroFESpace}
\end{equation}
where ${\mathcal P}^{p}$ is the space of polynomials on the element $K$, ${\mathcal T}_H$ the (quasi-uniform) triangulation of $\mathcal{B} \, \subset \, \mathbb{R}^{n_{dim}}$. For the characteristic macro element size $H$ it holds $H \gg \epsilon$ for efficiency. The space $\mathcal{S}^{p}_{\partial \mathcal{B}_D}$ is a subspace of $\mathcal{V}$.
%defined in \eqref{eq:HilbertSpaceV}.

For the solution of \eqref{eq:StrongFormMicro-2} in the macrodomain the two-scale finite element heterogeneous multiscale method FE-HMM is used. FE-HMM was introduced as an instance of the very general HMM by E and Engquist \cite{E.2003} and analyzed for elliptic PDEs in \cite{E.2005}, for linear elasticity in \cite{Abdulle.2006}. 

The solution of the macro problem in FE-HMM follows from the variational form. \\
Find $\bm u^H \in \mathcal{S}_{\mathcal{B}_D}(\mathcal{B}, \mathcal{T}_H)$ such that
\begin{equation}
\label{eq:VariationalFormulationHMM}
B_H (\bm u^H, \bm v^H) = \int_{\mathcal{B}} \bm f \cdot \bm v^H  \, dV \, + \,  \int_{\partial \mathcal{B}_N} \bar{\bm t} \cdot \bm v^H \, dA
\qquad \forall \bm v^H \in \mathcal S_{\partial \mathcal{B}_D} (\mathcal{B}, \mathcal{T}_H) \, ,
\end{equation}
which reads as a standard, single-scale finite element method.

\subsection{The modified macro bilinear form of FE-HMM}
\label{subsec:Modified-Macro-Bilinear-Form}
\begin{Figure}[htbp]
	\begin{minipage}{16.0cm}  
		\includegraphics[width=14.5cm, angle=0]{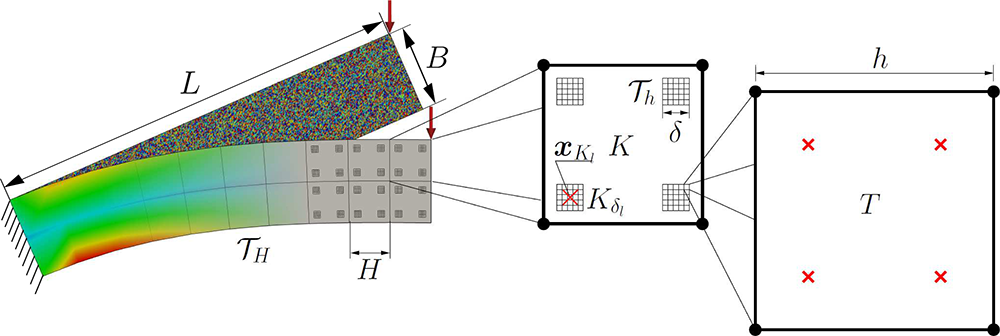}  
		\\[2mm]
		{\hspace*{3.5cm} (a) \hspace*{4.2cm} (b) \hspace*{3.0cm} (c)}
	\end{minipage}
	\caption{{\textbf FE-HMM as a two-scale finite element method}: (a) Macroscopic BVP with macrotriangulation $\mathcal{T}_H$, (b) one macro finite element $K$ of size $H$ with microdomains/RVEs $K_{\delta_l}$ of triangulation $\mathcal{T}_h$, centered at the macro quadrature points $\bm x_{K_{l}}$, (c) micro finite element $T$ of size $h$, here with Gauss quadrature points for $p$=$q$=1 on both scales.
		\label{fig:MicMac-Problem-Meshing-etc}}
\end{Figure}

Since the homogenized constitutive tensor $\mathbb{A}^{0}(\bm x)$ is typically not known for heterogeneous matter, the bilinear form $ B_H (\bm u^H, \bm v^H)$ cannot be calculated according to \eqref{eq:ModifiedBilinearForm-1} using standard numerical quadrature with $\bm x_{K_{\delta_l}}$ and $\omega_{K_{\delta_l}}$ the quadrature points and quadrature weights, respectively
\begin{eqnarray}
B_H (\bm u^H, \bm v^H) &=& \sum_{K\in \mathcal T_H} \sum_{l=1}^{N_{qp}} \omega_{K_{\delta_l}} \,
\left[ \mathbb{A}^{0} \bm \varepsilon^0 (\bm u^H(\bm x_{K_{\delta_l}})) : \bm \varepsilon(\bm v^H(\bm x_{K_{\delta_l}})) \right]
\label{eq:ModifiedBilinearForm-1} \\
&\approx& \sum_{K\in \mathcal T_H} \sum_{l=1}^{N_{qp}} \omega_{K_{\delta_l}} \, \left[ \dfrac{1}{{\color{black}|K_{\delta_l}|}} 
{\color{black} \int_{K_{\delta_l}} {\color{black} \mathbb{A}^{\epsilon} \bm \varepsilon^{\epsilon}}(\bm u^h_{K_{\delta_l}}) : \bm \varepsilon(\bm v^h_{K_{\delta_l}}) \, dV} \right] \, .
\label{eq:ModifiedBilinearForm-2}                           
\end{eqnarray}

Instead, FE-HMM approximates the virtual work expression at point $\bm x_{K_l}$ in the semidiscrete form \eqref{eq:ModifiedBilinearForm-1} by another bilinear form according to \eqref{eq:ModifiedBilinearForm-2} where the microheterogeneous elasticity tensor $\mathbb{A}^{\epsilon}$ is employed. 

According to this approximation, the solution $\bm u_{K_l}^h$ is obtained on microsampling domains $K_{\delta_l}=\bm x_{K_{\delta_l}} + \delta \, [-1/2, +1/2]^{n_{dim}}$, $\delta \geq \epsilon$, which are each centered at the quadrature points $\bm x_{K_{\delta_l}}$ of $K$, $l=1, \ldots, N_{qp}$. For a visualization see Fig.~\ref{fig:MicMac-Problem-Meshing-etc}. These microsampling domains with volume $|K_{\delta_l}|$ provide the additive contribution to the stiffness matrix of the macro finite element. In order to avoid too heavy notation we will replace $K_{\delta_l}$ by $K_{l}$.  

The approximation of \eqref{eq:ModifiedBilinearForm-1} by \eqref{eq:ModifiedBilinearForm-2} is the core of FE-HMM, a modified quadrature rule that fulfills Hill's postulate \cite{Hill.1963}, \cite{Hill.1972}, if energetically consistent boundary conditions are applied to the unit cell. 

\subsection{Variational formulation of the microproblem}
\label{subsec:Variational-FE-HMM-micro}

It can be shown that the FE-HMM microproblem resembles the discrete version of the cell problem of asymptotic expansion, if it is formulated for each microdomain $K_l$ in $K$ with $l=1, \ldots, N_{qp}$, $K \in \mathcal{T}_H$ like this:
\\[2mm]
Find $\bm u^h_{K_l}$ such that the conditions for macro-micro coupling and for the micro bilinear form \eqref{eq:micro-problem-vers1a} and \eqref{eq:micro-problem-vers1b} are fulfilled:
%\begin{equation}
%\label{eq:micro-problem-vers1}
%\renewcommand{\arraystretch}{1.6}
%\left.\begin{array}{rcl}
%\left(\bm u^h_{K_l} -  \bm u^{H}_{lin, K_l} \right) &\in& \mathcal{S}^q (K_l, \mathcal{T}_h)  \\[2mm] 
%B_{K_l}(\bm u^h_{K_l}, \bm w^h_{K_l})  &:=& 
%\displaystyle{\int_{K_l}} {\color{black}\bm \sigma^{\epsilon}}(\bm u^h_{K_l}) : \bm \varepsilon(\bm w^h_{K_l}) \, dV = 0  \\ 
%& &  \forall  \, \bm w^h_{K_l} \in \mathcal{S}^q (K_l, \mathcal{T}_h)   \, ,  
%\end{array} \quad \right\} \; 
%\end{equation} 
\begin{eqnarray}
\label{eq:micro-problem-vers1a}
\left(\bm u^h_{K_l} -  \bm u^{H}_{lin, K_l} \right) &\in& \mathcal{S}^q (K_l, \mathcal{T}_h)  \\[2mm] 
\label{eq:micro-problem-vers1b}
B_{K_l}(\bm u^h_{K_l}, \bm w^h_{K_l})  &:=& 
\displaystyle{\int_{K_l}} {\color{black}\bm \sigma^{\epsilon}}(\bm u^h_{K_l}) : \bm \varepsilon(\bm w^h_{K_l}) \, dV=0  \quad \forall  \, \bm w^h_{K_l} \in \mathcal{S}^q (K_l, \mathcal{T}_h)   \, ,  
\end{eqnarray} 
where the micro finite element space $\mathcal{S}^q_{}(K_l, \mathcal{T}_h)$ is defined by
\begin{equation}
\mathcal{S}^q (K_l, \mathcal{T}_h) = \{ \bm w^h \in \mathcal{W}(K_l); \bm w^h \in (\mathcal{P}^q(T))^{n_{dim}}, \, T\in \,\mathcal{T}_h \}\, .
\label{eq:Periodic-micro-FE-space}
\end{equation}
In \eqref{eq:Periodic-micro-FE-space} $\mathcal{T}_h$ is a quasi-uniform discretization of the sampling domain $K_l$ with mesh size $h \ll \epsilon$ resolving the finescale and $\mathcal{P}^q$ is the space of polynomials on the element $T$. The particular choice of the Sobolev space $\mathcal{W}(K_l)$ sets the boundary conditions for the micro problems, cf. \cite{Abdulle.2009}, Sec. 3.2. Among the coupling conditions that fulfill Hill's postulate we consider in this paper periodic boundary conditions (PBCs). The consideration of kinematically uniform displacement conditions (KUBC) and constant traction conditions (TBC), which are also energetically consistent, can be found e.g. in \cite{Fischer.2019}.

The linearization of $\bm u^H$ in \eqref{eq:micro-problem-vers1a} is carried out at the quadrature point $\bm x_{K_l}$  
\begin{equation}
\bm u^{H}_{lin, K_l} = \bm u^H (\bm x_{K_l}) + (\bm x - \bm x_{K_l}) \cdot \nabla \bm u^H(\bm x_{K_l}) \, .
\label{eq:linearization_uH}
\end{equation}
It ensures a homogeneous deformation on the microdomain and resembles therein the unit cell problem of asymptotic homogenization and thus is in the frame of strain-driven first order computational homogenization, for a discussion of these links see \cite{Eidel.2018}.  

%--------------------------------------------------------------------------------------- 
For the solution of \eqref{eq:micro-problem-vers1a}, \eqref{eq:micro-problem-vers1b} a basis $\{N_I^H\}_{I=1}^{M_{mac}}$ for the macro finite element space $\mathcal{S}^p_0(\mathcal{B}, \mathcal{T}_H)$ is employed in order to represent the macrosolution $\bm u^H$ of \eqref{eq:VariationalFormulationHMM}. Similarly, a basis $\{N_i^h\}_{i=1}^{M_{mic}}$ of the micro finite element space $\mathcal{S}^q_0(K_l, \mathcal{T}_h)$, \eqref{eq:Periodic-micro-FE-space}, is introduced in order to represent the solution $\bm u^h$ of a microproblem. $M_{mac}$ denotes the number of nodes of the macrodomain, and $M_{mic}$ denotes the number of nodes of each microdomain. 
Hence, the macro- and the microsolution follow the representation
\begin{equation}
\bm u^H = \sum_{I=1}^{M_{mac}} N_I^H \, \bm d_I^H\, , \qquad \bm u^h = \sum_{i=1}^{M_{mic}} N_i^h \, \bm d_i^h\, ,
\label{eq:micro-displacement-vector-in-the-fe-basis-shortened}
\end{equation}
where $\bm d_I^H$ is the displacement vector of macronode $I$, and $\bm d_i^h$ is the displacement vector for micronode $i$.

%---------------------------------------------------------------------------------------

\subsection{Macrostiffness calculation} 
\label{subsec:bottom-up-macrostiffness-calculation}

The macro bilinear form $B^e_H(\bm u^H, \bm v^H)$ is the virtual internal work for a macro finite element. The corresponding bilinear form in terms of the shape functions $B^e_H(\bm N_I^H, \bm N_J^H)$ results in the macro element stiffness matrix contribution $\bm k^{e,mac}_{IJ}$ for macronodes $I, J$, a $n_{dim} \times n_{dim}$ matrix. It holds  
\begin{equation}
\bm k^{e,mac}_{IJ} = B_H^{e} (\bm N_I^H, \bm N_J^H) = \sum_{l=1}^{N_{qp}} \dfrac{\omega_{K_l}}{|K_l|} 
\int_{K_l} (\bm L \bm u^{h(I)}_{K_l})^T \mathbb{A}^{\epsilon} (\bm x) \, \bm L \bm u^{h(J)}_{K_l} \, dV \, .
\label{eq:ModifiedBilinearForm-2-for-varphiH}   
\end{equation}
In \eqref{eq:ModifiedBilinearForm-2-for-varphiH} $\bm u^{h(I)}_{K_l}$ is the counterpart of $\bm u_{K_l}^h$ in \eqref{eq:micro-problem-vers1b}. It is the dimensionless solution of the microproblem on $K_l$, which is driven by the shape function $N_I^H$ at macronode $I$. In the following, we add $x_i, i=1, \ldots, n_{dim}$ to account for the vector-valued field problem of dimension $n_{dim}$. Consequently, $\bm u^{h(I,x_i)}_{K_l}$ is the microsolution driven by a macroelement unit-displacement state $\bm u^{H(I,x_i)}_{lin, K_l}$ at node $I$ in $x_i$-direction.

For stiffness calculation, problem \eqref{eq:micro-problem-vers1b} is reformulated in that $\bm u^{h(I,x_i)}_{K_l}$ replaces $\bm u^{h}_{K_l}$.

For the coupling of ${\bm u}^{H(I,x_i)}_{lin, K_l}$ with $\bm u^{h(I,x_i)}_{K_l}$ the two fields are expanded into the same basis $\{N_i^h\}_{i=1}^{M_{mic}}$ of $\mathcal{S}^q(K_l, \mathcal{T}_h)$,

\begin{equation}
{\bm u}^{H(I,x_i)}_{lin, K_l} = \sum_{m=1}^{M_{mic}} \, N^h_{m, K_l} {\bm d}^{H(I,x_i)}_{m} \, , 
\qquad 
\bm u^{h(I,x_i)}_{K_l} = \sum_{m=1}^{M_{mic}} N^h_{m, K_l} \, \bm d^{h(I,x_i)}_{m} \, .
\label{eq:varphi-h-by-psi-shortened}  
\end{equation}

%---------------------------------------------------------------------------------------
The solution of the microproblems for the minimizers $\bm d^{h(I,x_i)}$ is presented in Sec.~\ref{subsec:Solution-of-microproblems}. The macroelement stiffness matrix according to \eqref{eq:ModifiedBilinearForm-2-for-varphiH} yields after some algebra
%\begin{eqnarray}          
%\bm k^{e,mac}_{IJ} 
%&=& B^e_H\left[\bm N_I^H, \bm N_J^H\right]  \nonumber\\
%&=& \label{eq:k-mac-element-5}
%\sum_{l=1}^{N_{qp}} \dfrac{\omega_{K_l}}{|K_l|} 
%\, \left( \bm d^{h(I)} \right)^T \, \bm K^{mic}_{K_l} \,  \bm d^{h(J)}  \, ,   
%\end{eqnarray}  
\begin{equation}          
\bm k^{e,mac}_{IJ} 
%\, = \, B^e_H\left[\bm N_I^H, \bm N_J^H\right] 
\, = \, \sum_{l=1}^{N_{qp}} \dfrac{\omega_{K_l}}{|K_l|} 
\, \left( \bm d^{h(I)} \right)^T \, \bm K^{mic}_{K_l} \,  \bm d^{h(J)}  \, ,   
\end{equation}  
where $\bm d^{h(I)}=\left(\,\bm d^{h(I,x_1)} | \bm d^{h(I,x_2)} | \bm d^{h(I,x_3)} \, \right)$ for $n_{dim}=3$. The assembly of $\bm k^{e,mac}_{IJ}$ results in $\bm k^{e,mac}$ and corresponds to an assembly of $\bm d^{h(I)}$ for $I=1,\ldots, N_{node}$ in columns which gives the transformation matrix $\bm T_{K_l}$
\begin{eqnarray}
\bm k^{e,mac}_{K} &=& \sum_{l=1}^{N_{qp}} \dfrac{\omega_{K_l}}{|K_l|} 
\, \, \bm T^{T}_{K_l}\, \bm K^{mic}_{K_l} \, \bm T_{K_l} \label{eq:k-mac-element-6}   \\
\mbox{with} \quad \bm T_{K_l} &=& \bigg[ \Big[ \big[ \bm d^{h(I,x_i)} \big]_{i=1,\ldots,n_{dim}} \Big]_{I=1, \ldots, N_{node}} \bigg] \label{eq:k-mac-element-7} \,.        
\end{eqnarray}
 
In the present context of stiffness computation, a macro element shape function represents a unit displacement state for macro node  $I, I=1, \ldots, N_{node}$ in each direction of space $x_i\, |\, i=1, \ldots, n_{dim}$. 
They drive the microproblem in terms of the corresponding nodal values ${\bm d}^{H(I,x_i)}_{m}, m=1, \ldots, M_{mic}$ in each microdomain to evaluate the macroelement stiffness $\bm k^{e,mac}_{IJ}$. Each unit displacement state in $x_i$-direction induces in ${\bm d}^{H(I,x_i)}$ nonzero components only in $x_i$, for $n_{dim} = 3$ e.g. 
${\bm d}^{H(I,x_i)}|_{i=2} = \left[ 0, {d}^{H(I,x_2)}_{1,x_2}, 0, \hdots, 0, {d}^{H(I,x_2)}_{M_{mic},x_2}, 0\right]^T$. 

%---------------------------------------------------------------------------------------
\newcommand\myeq{\mathrel{\overset{\makebox[0pt]{\mbox{\normalfont\sffamily !}}}{=}}}

\subsection{Solution of the microproblems} 
\label{subsec:Solution-of-microproblems}

The total micro stiffness matrix for an RVE is obtained by standard processes of computing micro stiffness matrices through numerical quadrature and a consecutive assembly.  

The solution of the microproblem is obtained by the method of Lagrange multipliers, which renders the total energy for a macro unit displacement state 
\begin{eqnarray}
\mathcal{L} (\bm d^{h(I,x_i)}, \bm \lambda^{(I,x_i)}) &=& \dfrac{1}{2} \left(\bm d^{h(I,x_i)}\right)^T \bm{K}_{K_l}^{mic} \, \bm d^{h(I,x_i)}
+ \bm \lambda^{(I,x_i)\,T} \, \bm G \, \left(\bm d^{h(I,x_i)} - \overline{\bm d}^{H(I,x_i)}\right) 
\label{eq:Lagrange-functional} \\
& & \mbox{for }I=1, \ldots, N_{node}, \, \, \mbox{and} \, \, x_i\,|_{i=1, \ldots, n_{dim}} \, ,\nonumber
\end{eqnarray}
 
where $\bm G$ contains the coupling constraints of the RVE. The vector of Lagrange multipliers $\bm \lambda^{(I,x_i)} \in \mathbb R^{(1+L)\cdot n_{dim}}$, where $L$ depends on the type of microcoupling, represents forces on the microdomain which enforce the micro coupling condition. 
%reads for $n_{dim}=3$ as
%\begin{equation}
%\bm \lambda^{(I,x_i)}=\{\lambda^{(I,x_i)}_{0,x_1}, \lambda^{(I,x_i)}_{0,x_2}, %\lambda^{(I,x_i)}_{0,x_3}, \lambda^{(I,x_i)}_{1,x_1}, \lambda^{(I,x_i)}_{1,x_2}, %\lambda^{(I,x_i)}_{1,x_3}, 
%\ldots ,\lambda^{(I,x_i)}_{L,x_1}, \lambda^{(I,x_i)}_{L,x_2}, \lambda^{(I,x_i)}_{L,x_3}\}^T \, . 
%\end{equation}
The variations of $\mathcal{L}$ with respect to $\bm d^{h(I,x_i)}$ and to $\bm \lambda^{(I,x_i)}$ result in the stationarity conditions
\begin{equation}
\left[ \begin{array}{cc}
\bm K^{mic}_{K_l} & \bm G^T \\
\bm G  & \bm 0
\end{array} \right] 
\left[ \begin{array}{c}
\bm d^{h(I,x_i)} \\
\bm \lambda^{(I,x_i)}\\
\end{array}\right]
=
\left[ \begin{array}{c}
\bm 0   \\
\bm G \, {\bm d}^{H(I,x_i)}    \\
\end{array}\right]  \, \, \mbox{for} \, \, I=1, \ldots, N_{node}, \, \, i=1, \ldots, n_{dim} \, ,  
\label{Solve4alpha-lambda}
\end{equation}
hence, a system of linear equations with $N_{node} \cdot n_{dim}$ right hand sides. The solution vectors are augmented to full matrices, hence, $\bm d^{h(I,x_i)} \rightarrow \bm T$, $\bm \lambda^{(I,x_i)} \rightarrow \bm \Lambda$, ${\bm d}^{H(I,x_i)} \rightarrow \bm d^H$.   

The solution of \eqref{Solve4alpha-lambda} serves the purpose to compute the transformation matrix $\bm T_{K_l}$ according to \eqref{eq:k-mac-element-7}.  After the consecutive solution of the global macroproblem for $\bm u^H$, the microproblems have to be solved. Then, \eqref{Solve4alpha-lambda} is driven by the true macroscopic displacement vector $\bm d^H$, which results in the true microdisplacements $\bm d^{h}$.
 
Nonlinear extensions of FE-HMM in solid mechanics have been proposed in \cite{Eidel.2018b} and in \cite{Nejad.2019}, for applications of FE-HMM in the context of nonstandard homogenization for identifying material parameters in the Relaxed Micromorphic Model see \cite{dAgostino.2019} and \cite{Neff.2019}.

%SSRVE references:  \cite{Ohser.2006}, \cite{Balzani.2014}, \cite{Scheunemann.2015}, \cite{LeBris.2016}, \cite{Fischer.2019b}.

%\vfill
%\newpage
   
%---------------------------------------------------------------------------------------
\section{Errors and estimates}
\label{sec:ErrorsAndEstimates}
%---------------------------------------------------------------------------------------
\label{sec:Apriori-and-Aposteriori-ErrorEstimates}
  
\subsection{Errors on the microscale} 
%\label{subsec:A-priori-estimates}
 
For the image-based microstructure representation at finite resolution we propose to decompose the total micro error $e^{\epsilon}_{\text{mic}}$ on the microdomain\footnote{For notational convenience we prefer to write $\mathcal{B}_{\epsilon}$ instead of $K_{l}$ or $K_{\delta_{l}}$ for the discretized microdomain.} $\mathcal{B}_{\epsilon}$ into a discretization error $e^{\epsilon \, h}_{\text{mic}}$ and a modeling error $e^{\epsilon \, \square}_{\text{mic}}$ due to a coarsened resolution; it holds 
\begin{eqnarray}
 e^{\epsilon}_{\text{mic}} &\leq& e^{\epsilon \, h}_{\text{mic}} + e^{\epsilon \, \square}_{\text{mic}} \, ,  
 \label{MicroErrorTriangleInequality} \\
e^{\epsilon}_{\text{mic}} &:=& || \bm u^h(h_{\square}, h) - \bm u^h(h_{\square}\rightarrow 0, h \rightarrow 0)||_{(\mathcal{B}_{\epsilon})} \, ,
\label{eq:total-micro-error-at-microscale} \\ 
e^{\epsilon \, h}_{\text{mic}} &:=& || \bm u^h(h_{\square}, h) - \bm u^h(h_{\square}, h \rightarrow 0)||_{(\mathcal{B}_{\epsilon})} \, , 
\label{eq:micro-discretization-error-at-microscale} \\
e^{\epsilon \, \square}_{\text{mic}} &:=& || \bm u^h(h_{\square}, h\rightarrow 0) - \bm u^h(h_{\square}\rightarrow 0, h \rightarrow 0)||_{(\mathcal{B}_{\epsilon})} \, , 
\label{eq:micro-modeling-error-at-microscale}   
\end{eqnarray}
where $\bm u^h(h_{\square}, h)$ is the FE-HMM microsolution at discretization $h$ and resolution $h_{\square}$ with the pixel or voxel size $h_{\square}$, $\bm u^h(h_{\square}\rightarrow 0, h \rightarrow 0) =: \bm u^{\epsilon}$ the exact solution and $\bm u^h(h_{\square}, h \rightarrow 0)$ the exact solution for finite $h_{\square}$.
The condition $h_{\square} \rightarrow 0$ in the reference resolution for the modeling error is a task of image acquisition in contrast to obtaining a reference solution for the discretization error $h \rightarrow 0$. The reference resolution can be finite, hence $h_{\square} \rightarrow h^{\ast}_{\square}$, which is then to be considered in \eqref{eq:total-micro-error-at-microscale} and \eqref{eq:micro-modeling-error-at-microscale}. It generally must hold $h \leq h_{\square}$. 

For the micro error measured on the microscale $e^{\epsilon}_{\text{mic}}$ standard a priori estimates of the finite element method hold for sufficient regularity. With the extension by the modeling error $e^{\epsilon \, \square}_{\text{mic}}$, the estimates read in the $L^2$- and in the energy norm $A$ on a microdomain $\mathcal{B}_{\epsilon}$  
\begin{alignat}{4}  %  number of & signs is 2n-1, with n: the number of columns.
%e^{\epsilon}_{\text{mic}} &\leq  && \quad e^{\epsilon \, h}_{\text{mic}} &&+   & \  &e^{\epsilon \, \square}_{\text{mic}}  \\[2mm] 
%\textrm{with} & & \nonumber \\
|| \bm u^{\epsilon} - \bm u^h ||_{L^2(\mathcal{B}_{\epsilon})} &\leq \  &&C  \left(\dfrac{h}{\epsilon}\right)^{q+1} &&+ & \, &e^{\epsilon \,\, \square}_{\text{mic}}    \, ,
\label{eq:Total-Error-estimate-L2-mic} \\
%|| \bm u^{\epsilon} - \bm u^h ||_{H^1(\mathcal{B}_{\epsilon})}   &\leq \  &&C  %\left(\dfrac{h}{\epsilon}\right)^{q} &&+ & \,\, &e^{\epsilon \, \square}_{\text{mic}}   \, ,
%\label{eq:Total-Error-estimate-H1-mic}  \\
|| \bm u^{\epsilon} - \bm u^h ||_{A(\mathcal{B}_{\epsilon})}   &\leq \  &&C  \left(\dfrac{h}{\epsilon}\right)^{q} &&+ & \,\, &e^{\epsilon \, \square}_{\text{mic}}   \, .
\label{eq:Total-Error-estimate-Energy-mic}  
\end{alignat}

For the coefficients of the homogenized elasticity tensor it holds  
\begin{equation}
\label{eq:ConvergenceElasticities}
|| \mathbb{A}^0_{ijkl} - \mathbb{A}^{0,h,\square}_{ijkl}|| \leq C \left( \dfrac{h}{\epsilon} \right)^{2q} + e^{\epsilon \, \square}_{\text{mic}}
\end{equation}
with the exact homogenized elasticity tensor $\mathbb{A}^0_{ijkl}$ and its approximation 
$\mathbb{A}^{0,h}_{ijkl}$ obtained at micro element size $h$. 

While the estimate for the homgenized elasticities is given according to \eqref{eq:ConvergenceElasticities}, the corresponding relative error $|| \mathbb{A}^0_{ijkl} - \mathbb{A}^{0,h,\square}_{ijkl}||/|| \mathbb{A}^0_{ijkl}||$ would be not fully descriptive, since it does not include the stiffness contrast of material phases. For that reason the relative error for a microstructure consisting of two phases, 1 and 2, at discretization $h$ and resolution $h_{\square}$ is computed according to 
\begin{eqnarray}
\label{eq:total-micro-error-Cijkl}
e^{\epsilon}_{\text{mic}}(\mathbb{A}^{0,h,\square}_{ijkl}) &:=& 
\dfrac{||\mathbb{A}^{0}_{ijkl} - \mathbb{A}^{0,h,\square}_{ijkl}||}{\text{min} \{||\text{dev}\, \mathbb{A}^{0,r}_{ijkl}||, r=1,2 \}} 
\\
\text{with} \quad \text{dev} \,  
\mathbb{A}^{0,r}_{ijkl} &:=& \mathbb{A}^0_{ijkl} - \mathbb{A}^{0,r}_{ijkl} \, , \quad r=1,2 \, .
\label{eq:phase-contrast}
\end{eqnarray}
Through the deviation of the exact homogenized elasticities from those of the single phases $\mathbb{A}^{0,r}_{ijkl}$ in \eqref{eq:phase-contrast} the bounds set by the individual phases are included.

Similar to \eqref{eq:micro-discretization-error-at-microscale} and \eqref{eq:micro-modeling-error-at-microscale} the total error \eqref{eq:total-micro-error-Cijkl} can be decomposed into modeling and discretization parts, the latter by the computation of the reference solution $\mathbb{A}^{0}_{ijkl}$ at fixed $h_{\square}$ and $h\rightarrow 0$.  

\subsection{Errors on the macroscale} 

The FE-HMM provides unified a priori estimates; for elliptic partial differential equations (PDEs) derived in \cite{E.2005}, \cite{Ohlberger.2005}, for linear elasticity in \cite{Abdulle.2006}, comprehensively described in \cite{Abdulle.2009}. 
  
The total FE-HMM error can be decomposed into three parts measured at the macroscale  
\begin{equation}
   || \bm u^0 - \bm u^H || \, \leq \, \underbrace{|| \bm u^0 - \bm u^{0,H} ||}_{\displaystyle e_{\text{mac}}} 
                             \, + \, \underbrace{|| \bm u^{0,H} - \widetilde{\bm u}^H ||}_{\displaystyle e_{\text{mod}}}
                             \, + \, \underbrace{|| \widetilde{\bm u}^H - \bm u^H ||}_{\displaystyle e_{\text{mic}}} \, ,
   \label{eq:Error-decomposition-mac-mod-mic}
\end{equation}
where $e_{\text{mac}}$, $e_{\text{mod}}$, $e_{\text{mic}}$ are the macro error, the modeling error, and the propagated micro error.

In \eqref{eq:Error-decomposition-mac-mod-mic}, $\bm u^0$ is the solution of the homogenized problem \eqref{eq:Homogenized-Strong-Form}, $\bm u^H$ the FE-HMM solution, $\bm u^{0,H}$ is the standard (single-scale) FEM solution of problem \eqref{eq:VariationalFormHomogenizedProblem} that is obtained through exact $\mathbb{A}^{0}$; and $\widetilde{\bm u}^H$ is the FE-HMM solution obtained through exact microfunctions (in $W(K_l)$).

For sufficient regularity the a priori estimates hold
\begin{eqnarray}
   || \bm u^0 - \bm u^H ||_{L^2(\mathcal{B})} &\leq& C\left( H^{p+1} + \left(\dfrac{h}{\epsilon}\right)^{2q} \right) + e_{\text{mod}}    \, ,
    \label{eq:Total-Error-estimate-L2} \\
%   || \bm u^0 - \bm u^H ||_{H^1(\mathcal{B})} &\leq& C\left( H^p + \left(\dfrac{h}{\epsilon}\right)^{2q} \right) + e_{\text{mod}}   \, ,
%   \label{eq:Total-Error-estimate-H1} \\
   || \bm u^0 - \bm u^H ||_{A(\mathcal{B})} &\leq& C\left( H^p + \left(\dfrac{h}{\epsilon}\right)^{2q} \right) + e_{\text{mod}}   \, .
   \label{eq:Total-Error-estimate-Energy}  
\end{eqnarray}
 
For $e_{\text{mod}}$ in \eqref{eq:Total-Error-estimate-L2}, \eqref{eq:Total-Error-estimate-Energy} it holds 
\begin{equation}
\label{eq:ModelingErrorDecomposition}
e_{\text{mod}} = e_{\text{mod\,BC}} + e^{\square}_{\text{mic}}     
\end{equation}
where $e_{\text{mod\,BC}}$ is a modeling error due to boundary conditions (BC) on the RVE and $e^{\square}_{\text{mic}}$ is the newly introduced finite-resolution modeling error. For $e_{\text{mod\,BC}}$ it holds
\begin{equation}
\label{eq:ModelingError}
 e_{\text{mod\,BC}} = \left\{ \begin{array}{ll}
         0 & \mbox{for periodic coupling with} \,\, \delta/\epsilon \in \mathbb{N} \\
         {C\,\dfrac{\epsilon}{\delta}}  & \mbox{for Dirichlet coupling with} \,\, \delta > \epsilon \end{array} \right. \,    
\end{equation} 
given that the hypotheses hold, that the elasticity tensor $\mathbb{A}^{\epsilon}$ is periodic on the RVE and, that the micro solution is sufficiently smooth, \cite{Jecker.2016}. 

The modeling error $e_{\text{mod\,BC}}$ for Dirichlet coupling in \eqref{eq:Total-Error-estimate-L2}--\eqref{eq:Total-Error-estimate-Energy} is due to boundary layers \cite{E.2005} (Thm. 1.2), \cite{Abdulle.2009} and remains as a residual even for $H \rightarrow 0$ and $h \rightarrow 0$. 

\bigskip

{\bf Remarks}

For its composition covering both the macro error as well as the micro error, the estimates \eqref{eq:Total-Error-estimate-L2}--\eqref{eq:Total-Error-estimate-Energy} enable strategies to achieve the optimal convergence order for minimal computational costs in uniform micro-macro discretizations, for an assessment see \cite{Eidel.2016}, \cite{Jecker.2016}, \cite{Eidel.2018}. 

The above unified error estimates derived for FE-HMM apply for FE$^2$ for the equality of the methods \cite{Eidel.2018}. FE-HMM and FE$^2$ were developed independently, the latter in \cite{Golanski.1997}, \cite{Moulinec.1998}, \cite{Smit.1998}, \cite{Miehe.1999}, \cite{Feyel.2000}, \cite{Kouznetsova.2001}.

The computation of the micro errors is carried out by projection of the approximate solution to the Gauss points of a reference mesh ($h \rightarrow 0$) with a consecutive numerical quadrature for approximating the integral element error as detailed in \cite{Fischer.2019}, \cite{Fischer.2020}.

The estimation of the discretization error is described in  Sec.~\ref{sec:error_estimation_average} with a focus on peculiarities at interfaces, which have undergone resolution coarsening as detailed in Sec. \ref{sec:Mesh-coarsening}. 

For both error computation and error estimation on the microscale the macro displacement field $\bm u^{0,H}$ is used for consistency.   

The modeling micro error is not directly accessible to error estimation and therefore the total micro error neither. 

\section{Coarsening of image resolution and mesh discretization}
\label{sec:Mesh-coarsening}

The present section proposes two types of coarsening, the first is a uniform coarsening of microstructure resolution which implies a coarsening of finite element discretization. Consequently, a modeling error and a discretization error are introduced. The second type is a non-uniform mesh coarsening, which keeps the image resolution constant and therefore the modeling error likewise, but increases the discretization error. 
These two types of coarsening will be combined in that uniform coarsening of microstructure resolution is followed by adaptive mesh coarsening.

\subsection{Uniform pixel coarsening}
\label{subsec:uniform_coarsening}

In resolution coarsening all pixels of the initial, finely resolved microstructure image undergo uniform coarsening, no matter if they are located on phase boundaries or inside of phases. For phase boundaries rules must be established, how a patch of $2 \times 2$ fine pixels having different stiffness (represented by different color codes) merge into one coarser pixel. Since the finite element discretization follows this uniform coarsening, pixels represent finite elements and pass their properties onto them.

In the following two different variants of pixel coarsening are presented.

Variant A follows a rule of mixtures; the newly created pixel exhibits properties of the volume average of the pixels merging in that coarser pixel, hence a coarsening that creates new interphases at interfaces. If the number of phases is known and their interfaces sharp, the newly created interphases are artefacts, which smoothen discrete interfaces, reduce their stiffness mismatch and consequently reduce corresponding stress jumps. The volume average over the entire microdomain however is preserved. 

Variant B aims to overcome the drawbacks of a coarsening that introduces new interphases. For coarsened pixels at interfaces the quantity and properties of the initial phases are preserved following the rule, the majority wins. In case of an equal count the newly created coarser pixel is endowed with those phase properties that shift the overall phase ratio closer to the original one.  

\begin{Figure}[htbp]
	\centering
	\subfloat[]
	{\includegraphics[height=3.6cm,angle=0]{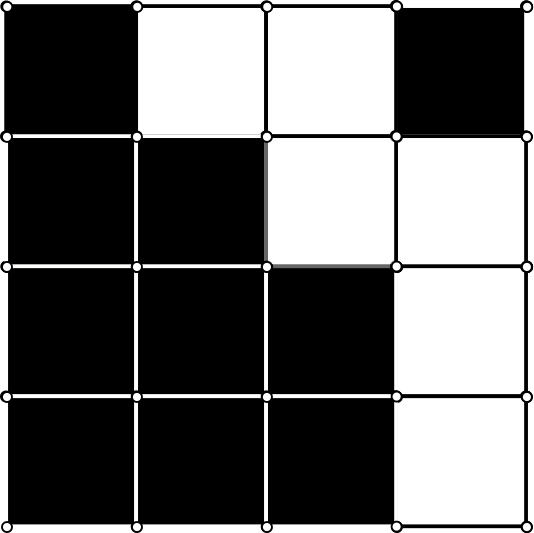}} 
	\hspace*{0.04\linewidth}
	\subfloat[]
	{\includegraphics[height=3.6cm,angle=0]{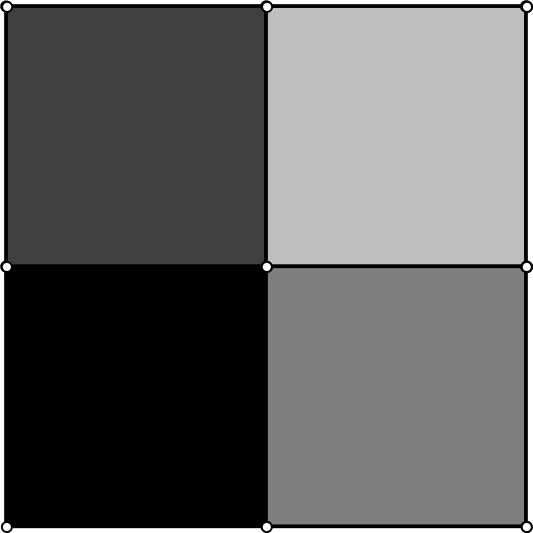}}
	\hspace*{0.04\linewidth}
	\subfloat[]
	{\includegraphics[height=3.6cm,angle=0]{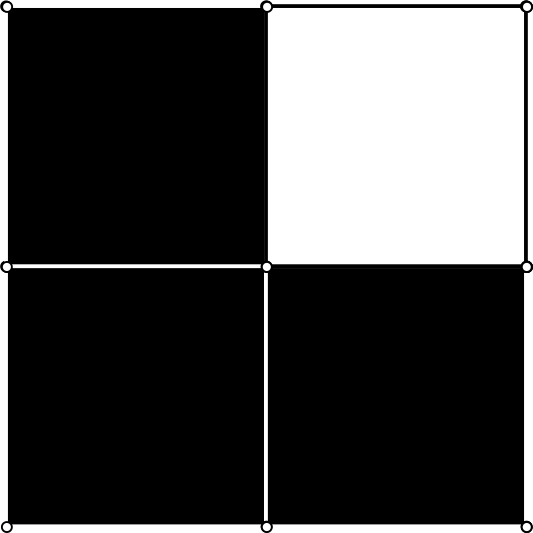}}
	\hspace*{0.02\linewidth}
	
	\caption{\textbf{Uniform mesh coarsening}: (a) Initial, fine discretization, (b) coarsened mesh with intermediate phases and (c) coarsened mesh preserving the initial phases in quantity and properties.}
	\label{fig:uniform_mesh_coarsening}
\end{Figure}

Figure \ref{fig:uniform_mesh_coarsening} illustrates the two variants of coarsening for a uniform mesh with elements belonging to two different phases, Fig. \ref{fig:uniform_mesh_coarsening} (a). Coarsening along with intermediate phases results in the mesh of Fig.~\ref{fig:uniform_mesh_coarsening} (b). Only for the coarsened element in the lower left no averaging is required. The other three elements obtain their properties from averaging, which is represented in the image by the corresponding color code.

Applying the second approach for uniform mesh coarsening leads to the mesh in Fig. \ref{fig:uniform_mesh_coarsening} (c). The majority-wins rule is ambiguous only for the coarsened element in the lower right, the assignment of the black phase leads to a phase ratio closer to the initial one. 

Notice that coarsening variant A generates results of the type how raw image data e.g. in bitmap format are typically represented. Variant B follows the process of image segmentation applied to microstructures with intermediate phases in terms of intermediate color codes. As a result, using a microstructure representation obtained from coarsening variant A can be interpreted as skipping image segmentation at all.  
 
\subsection{Quadtree-type, adaptive mesh coarsening}
\label{subsec:Quadtree_Coarsening}

The second type is a non-uniform, adaptive mesh coarsening, which keeps the image resolution constant and therefore the modeling error likewise, but increases the discretization error. It starts out from a --possibly already coarsened-- uniform mesh.

In this procedure realized by a quadtree-type mesh coarsening algorithm the interface resolution is maintained for accuracy, in the interior of phases mesh coarsening is carried out for efficiency. This type of microstructure-informed adaptivity is frequently used for multiphase materials \cite{Mishnaevsky.2005}, 
\cite{Legrain.2011}, \cite{Lian.2013}, \cite{Ren.2015}, \cite{Miska.2019} for application in homogenization and fracture, in the context of the Scaled Boundary FEM \cite{Saputra.2017}, \cite{Gravenkamp.2017} and many more.
But only very recently a rigorous accuracy analysis based on error estimation was presented in \cite{Fischer.2020}. 

\begin{Figure}[htbp]
	\centering
	\subfloat[ ]
	{\includegraphics[width=0.26\linewidth]{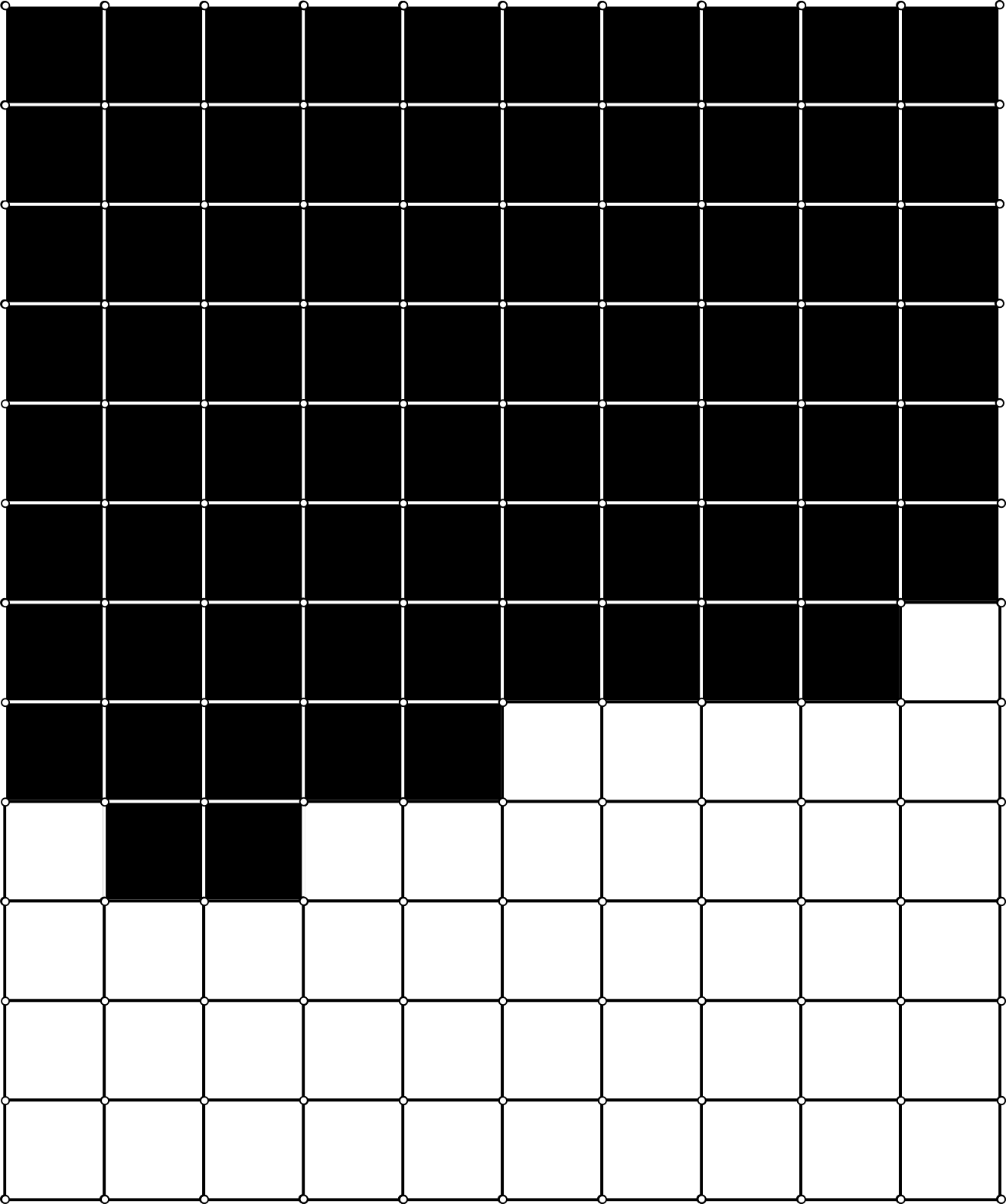}} \hspace*{0.04\linewidth}
	\subfloat[ ]
	{\includegraphics[width=0.26\linewidth]{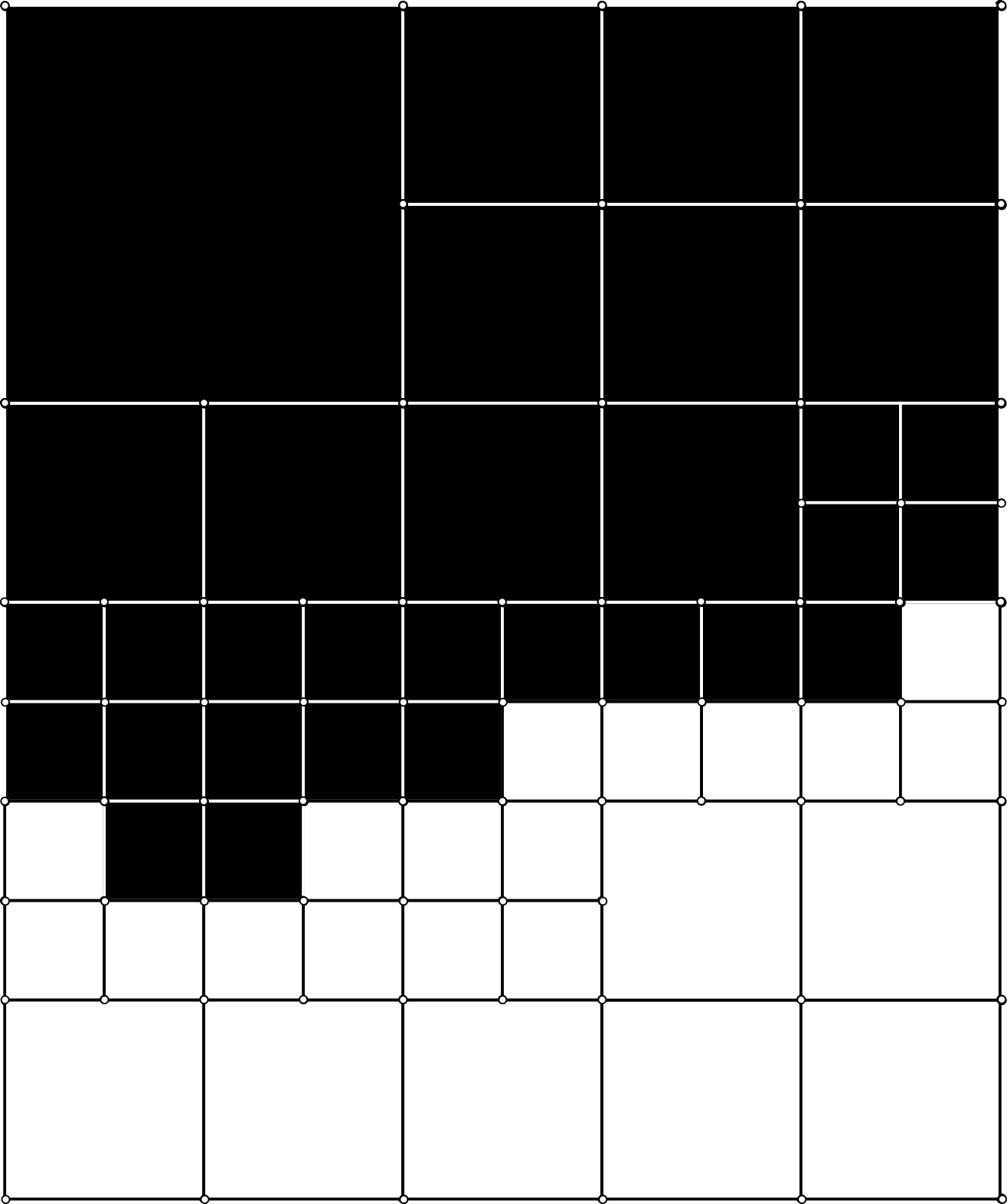}} 
	\caption{\textbf{Adaptive mesh coarsening}: (a) Original uniform mesh, (b) quadtree-type adaptively coarsened mesh.}
	\label{fig:Quadtree}
\end{Figure}

As an example, Fig.~\ref{fig:Quadtree} (a) displays an interface in a two-phase microstructure in its initial uniform pixel resolution coinciding with the micro mesh of element size $h$. The adaptive, quadtree-type coarsening results in the discretization of Fig.~\ref{fig:Quadtree} (b). It preserves the fine resolution at the phase boundary. Inside of the two phases the elements have been coarsened to new elements of side length $2h$ and $4h$. 

Uniform pixel coarsening in the present paper is followed by a non-standard adaptive mesh refinement. Here, adaptivity is not directed by the distribution of an a posteriori error estimate. Instead, it is microstructure-informed in that high resolution at interfaces is maintained, whereas mesh-coarsening is carried out in the interior of phases. Suchlike preprocessing typically provides a favorable balance of accuracy and efficiency. The obtained, adaptively refined mesh can then be assessed by a posteriori error estimation for various loading conditions \cite{Fischer.2020}.

\section{Error estimation for coarsened microstructures}
\label{sec:error_estimation_average} 

For reconstruction-type error estimation according to Zienkiewicz-Zhu \cite{Zienkiewicz.1987}, \cite{Zienkiewicz.1992}, \cite{Zienkiewicz.1992b} the accuracy critically depends on improved nodal stresses $\boldsymbol{\sigma}^\star$ and strains $\boldsymbol{\varepsilon}^\star$   
\begin{eqnarray}
(\bar{e}^{\epsilon \, h}_{\text{mic}})^2 = || \bar{\bm e} ||^2_{A(\Omega_{\epsilon})} 
%= || \bm u^\star - \bm u^h ||_{A(\Omega)} 
&=&  \int_{\Omega_{\epsilon}} \left( \boldsymbol{\sigma}^\star - \boldsymbol{\sigma}^h \right) \colon \left( \boldsymbol{\varepsilon}^\star - \boldsymbol{\varepsilon}^h \right) \, \text{d}V   \\
&\approx& \ \sum_{T \in \mathcal{T}_{h}}^{} \left( \sum_{i=1}^{ngp} \omega_i \left( \boldsymbol{\sigma}^\star - \boldsymbol{\sigma}^h \right)(\bm x_i^h) \colon \left( \boldsymbol{\varepsilon}^\star - \boldsymbol{\varepsilon}^h \right)(\bm x_i^h) \ \text{det} \bm J \right) \, .
\label{eq:error_estimator}
\end{eqnarray}
In contrast to the true error \eqref{eq:micro-discretization-error-at-microscale} based on reference stresses obtained on a discretization $h\rightarrow 0$, improved stresses in the estimate are obtained on the same mesh. Improved stresses for a node in a phase interior are calculated by simple averaging the values extrapolated from quadrature points of adjacent elements to that node or by more sophisticated techniques exploiting superconvergence \cite{Zienkiewicz.1992}, \cite{Zienkiewicz.1992b}. Stress averaging for a node at the discrete interface of two different phases ignores the stiffness mismatch and corresponding stress jumps. A phase distinction assigning two different sets of stress to a node (referred to as duplex stress) as displayed in Fig. \ref{fig:stress_duplex_vs_quadruplex} (a) does not only overcome these drawbacks but also results in more accurate error estimates compared to a standard stress averaging that ignores interfaces \cite{Fischer.2020}. It is therefore the method of choice for stress computation at interfaces for phase-preserving coarsening, variant B in this work. 
\begin{Figure}[htbp]
	\centering
	\subfloat[ ]
	{\includegraphics[width=0.40\linewidth]{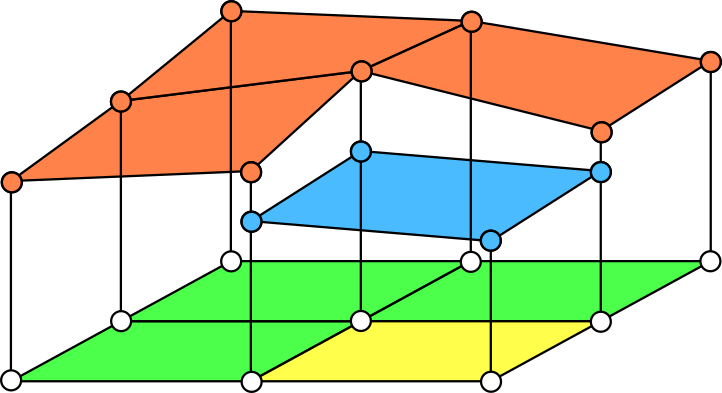}} \hspace*{0.01\linewidth}
	\subfloat[ ]
	{\includegraphics[width=0.40\linewidth]{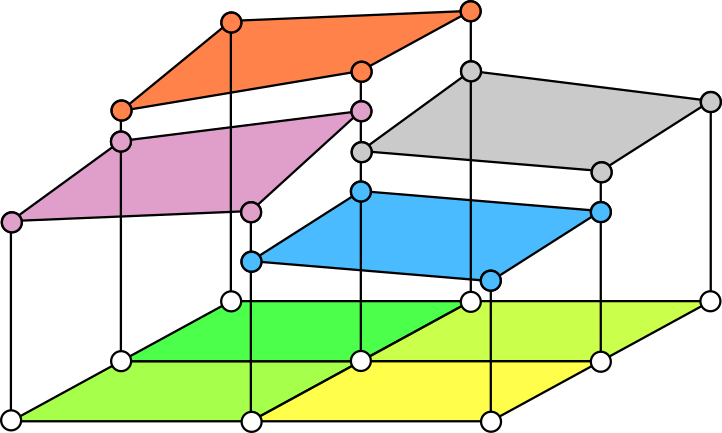}}
	\caption{\textbf{Stress distributions at phase boundaries}: (a) Duplex stress at the discrete interface between two phases, a stiff one in green and a more compliant one in yellow, (b) quadruplex stress at a boundary with intermediate phases.}
	\label{fig:stress_duplex_vs_quadruplex}
\end{Figure}

If the scheme of phase-distinction in nodal stress computation is consistently applied to interfaces following from coarsening variant A thus creating new interphases, a new situation arises. Figure \ref{fig:stress_duplex_vs_quadruplex} (b) displays the case of a $2 \times 2$ element patch each with different phases. Consequently, the central patch node is endowed with a quadruplex stress set. 

A new, related issue shall be illustrated in Fig. \ref{fig:stress_transfer}
for nodes of an isolated element having throughout different phases in direct neighborhood. Notice that stresses are computed in the first place in the quadrature points, for $q=1$ in the $2 \times 2$ gauss points marked with squares. Next, these stress values are extrapolated to the element nodes by the shape functions, a transfer marked by red arrows. For error estimation according to \eqref{eq:error_estimator} the improved nodal stresses $\boldsymbol{\sigma}^\star$ have to be transferred back to the quadrature points as marked by green arrows. When no averaging is carried out at all nodes due to phase-distinction, the nominally improved stresses transferred back into the quadrature points coincide with the original stress values therein, which results in a vanishing error in the stresses. This is clearly an artefact that spoils error estimation. 

\begin{Figure}[htbp]
	\centering
	{\includegraphics[width=0.40\linewidth]{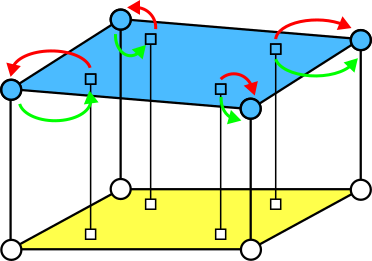}}
	\caption{\textbf{Stress-transfer for an element surrounded by elements of different phases}: Stresses from quadrature points (squares) are transferred (red arrows) to nodes (circles), for error estimation they are transferred back into the quadrature points and compared with the original values (green arrows). Identical values result in a zero value for the estimated error.}
	\label{fig:stress_transfer}
\end{Figure}
  
Notice that for a microstructure starting from two different phases in the initial resolution, $n$ coarsening steps could lead to $2^{2n} + 1$ 
%{\color{blue} Muesste es hier $2^{2n}+1$ sein? Fuer n=2 landen wir sonst bei 9 Phasen, aber tatsaechlich besteht ein 2-mal vergroebertes Element aus 16 'urspruenglichen' Elementen, also sind 17 Phasen moeglich?!} 
newly created phases at maximum, where an interface node is typically part of four different adjacent phases. Hence, the setup in Figs. \ref{fig:stress_duplex_vs_quadruplex} (b) and \ref{fig:stress_transfer} is rather the rule than an exception for coarsening variant A. 
 
As a consequence, the underestimation of the true error is expected locally at interfaces where it is typically largest and as an overall error of the total microdomain as well. This behavior shall be underpinned by an example. We compare at interface nodes the quadruplex stress computation with the standard concept where nodal stresses are obtained by averaging the elementwise nodal stresses of all adjacent elements independent from their phase.

\begin{Figure}[htbp]
	\centering
	\subfloat[$e^{h}_{\text{mic}}$]
	{\includegraphics[width=0.32\linewidth]{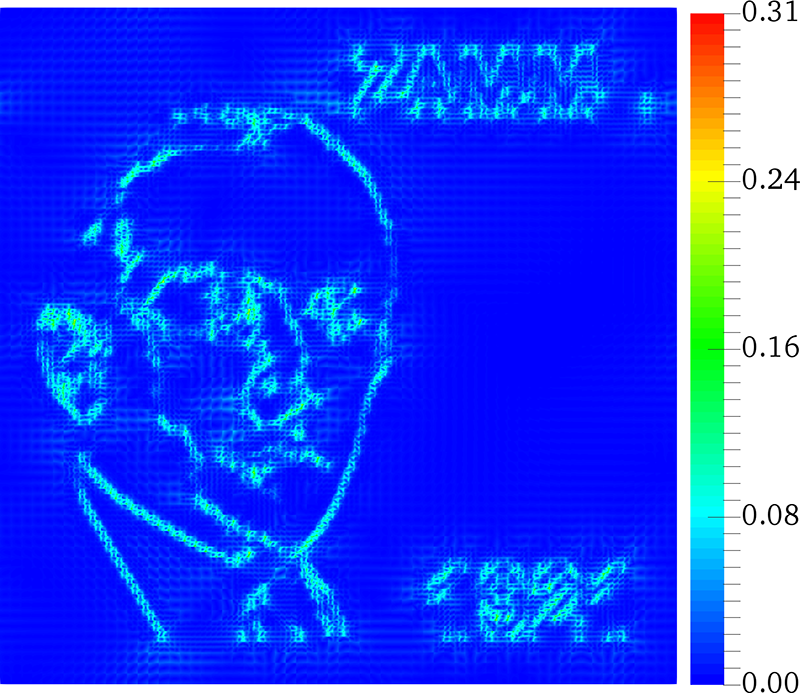}} \hspace*{0.01\linewidth}
	\subfloat[$\bar{e}^{h}_{\text{mic}}$(quad)]
	{\includegraphics[width=0.32\linewidth]{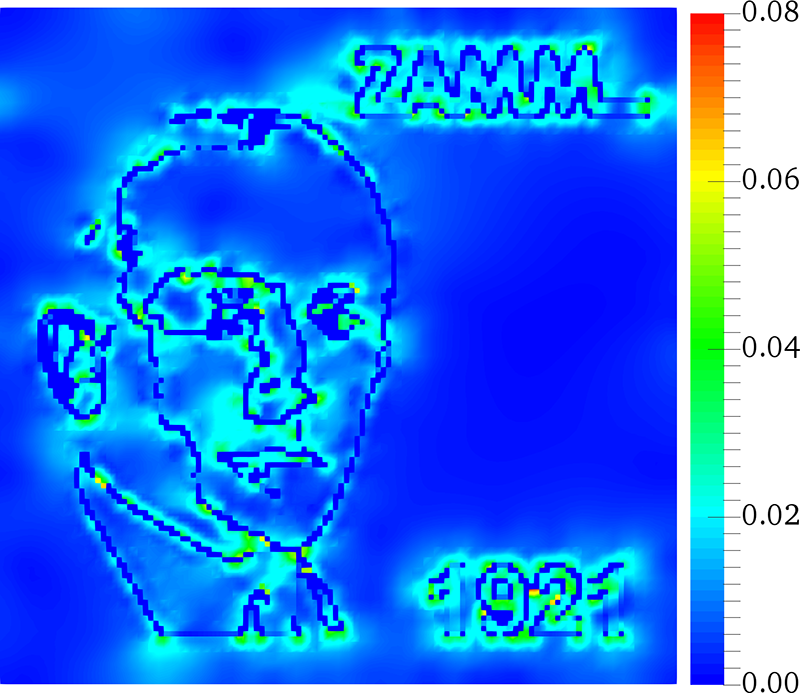}} \hspace*{0.01\linewidth}
	\subfloat[$\bar{e}^{h}_{\text{mic}}$(aver)]
	{\includegraphics[width=0.32\linewidth]{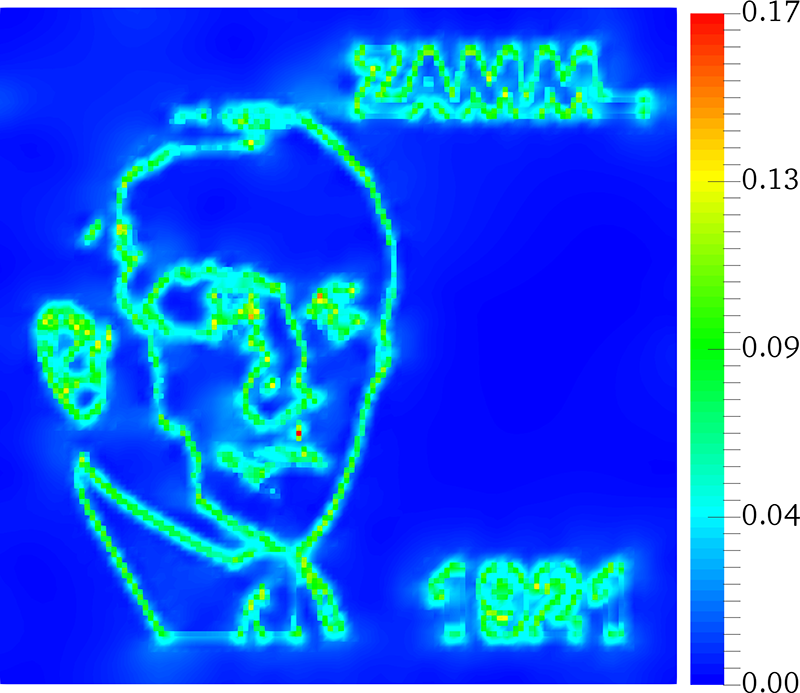}}
	\caption{\textbf{Relative discretization errors for mesh with intermediate phases}: (a) Calculated error $e^{h}_{\text{mic}}$, (b) estimated error based on quadruplex stresses $\bar{e}^{h}_{\text{mic}}$(quad), (c) estimated error based on averaged stress $\bar{e}^{h}_{\text{mic}}$(aver). All relative errors per element are obtained by division by the energy of the element $||\bm u||_{A}$.
	\label{fig:error_distribution_calc_vs_duplex_vs_average}}
\end{Figure}

\begin{table}[htbp]
	\centering
	\begin{tabular}{c c c c c}
		\hline
		$e^{\epsilon \, h}_{\text{mic}}$ & $\bar{e}^{\epsilon \, h}_{\text{mic}}$(quad) & $\Theta$(quad)  & $\bar{e}^{\epsilon \, h}_{\text{mic}}$(aver) & $\Theta$(aver) \\
		$15.8114$ & $8.2921$ & $0.5244$ & $22.1951$ & $1.3865$ \\
		\hline
	\end{tabular}
	\caption{\textbf{Comparison of discretization errors:} Calculated and estimated errors with their efficiency index $\Theta$. All error data in $10^{-4}$ (MPa).}
	\label{tab:error_calc_vs_est}
\end{table}

Figure~\ref{fig:error_distribution_calc_vs_duplex_vs_average} displays the relative discretization error for an example which will be discussed in detail in Sec.~\ref{sec:NumericalExamples}. The discretization used here follows from three consecutive steps of uniform coarsening according to variant A thus creating intermediate phases.

The accurate discretization error based on a fine reference solution is displayed in \ref{fig:error_distribution_calc_vs_duplex_vs_average} (a). Estimated errors for nodal quadruplex stresses are shown in Subfig.~(b), for nodal stress averaging without phase distinction in Subfig.~(c). Notice that the latter version qualitatively captures the true distribution of discretization errors showing maxima at interfaces. The error distribution based on quadruplex stresses however exhibits the artefact of vanishing errors at the phase boundary, in agreement with our elementary reasoning.

Table~\ref{tab:error_calc_vs_est} shows the numbers of the calculated and estimated discretization errors for the full microdomain along with the efficiency indices. While the estimated error based on quadruplex stresses $\bar{e}^{\epsilon \, h}_{\text{mic}}$(quad) underestimates the calculated discretization error $e^{\epsilon \, h}_{\text{mic}}$ due to the vanishing error on the phase boundary, the estimation based on averaged stresses overestimates the error.
Qualitatively, at least the estimation based on averaged nodal stresses renders acceptable distributions. Quantitatively however, none of both methods achieves accurate results as indicated by efficiency indices considerably deviating from unity, again consistent with the above reasoning.  

Notice that here, the 'crime' of resolution coarsening that creates new, additional interphases with its inherent artefacts, must be alleviated by another 'crime', the required stress averaging in order to avoid further artefacts in error estimation. 

\section{Numerical example}   
\label{sec:NumericalExamples} 
%----------------------------------------------------------------------

\subsection{Two-phase microstructure} 

This section provides an error analysis of an image-based microstructure representation which undergoes coarsening in resolution and in discretization. The modeling error introduced by the coarsened resolution and the standard finite element discretization error according to  \eqref{MicroErrorTriangleInequality}--\eqref{eq:micro-modeling-error-at-microscale} are investigated.  

We consider the two-phase microstructure of Fig. \ref{fig:Types-of-coarse-graining} (a) subject to periodic boundary conditions; the isotropic, linear elasticity of the inclusion phase (white) is characterized by $E_i=192.1$~MPa, $\nu_i=0.2$, the matrix phase by $E_m=100.0$~MPa, $\nu_m=0.2$. The square unit cell of edge length $\epsilon=1$~mm exhibits an initial uniform pixel resolution of 1024$^2$, which is equal to the finite element discretization, hence $h_{\square}=h=1/1024$~mm.
The plane-strain macro problem is a cantilever beam of length $L=5000$~mm, height $B=1000$~mm and width $D=100$~mm, see Fig.~\ref{fig:macroproblem}. It is loaded at its free end by a line-load of $q_0=0.02$~N/mm. Postprocessing is carried out on a microdomain at $x=y=2.1132$~mm. 
%{\color{blue} Ich habe bei $q_0=0.02$ und $x=y=2.1132$ jeweils die ganze Formel in den Mathematik-Modus gepackt, um gleichmaessige Abstaende in den Gleichungen zu haben. Vorher war es $q_0$=0.02 bzw. $x$=$y$=2.1132. War das gewollt?}

\begin{Figure}[htbp]
	\centering
	\includegraphics[width=12.0cm,angle=0]{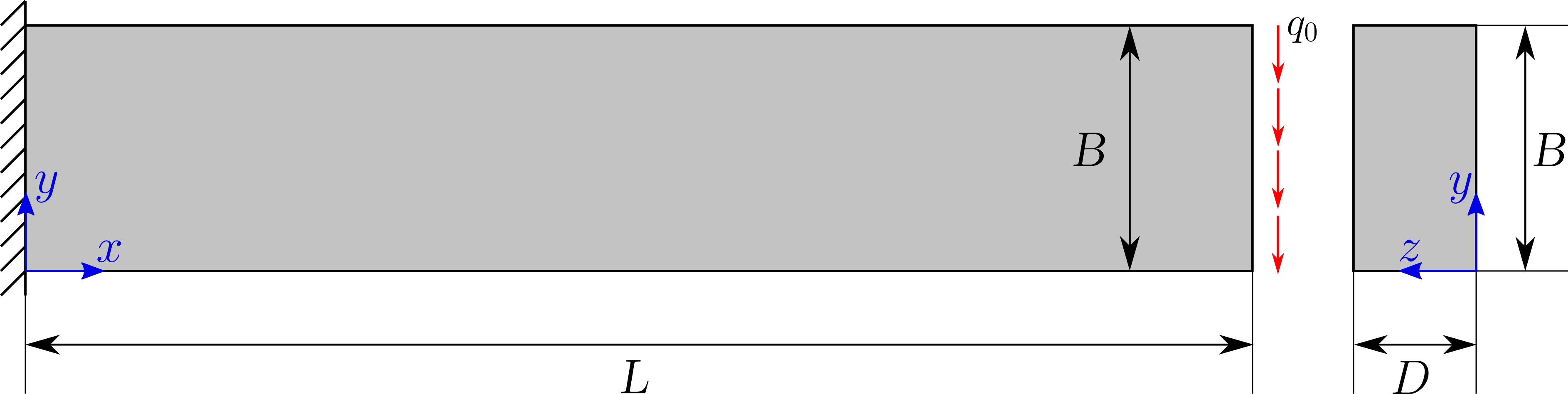}
	\caption{\textbf{Macroproblem:} cantilever subject to a line load $q_0$ at $x=L$.}
	\label{fig:macroproblem}	
\end{Figure}
\subsection{Uniform pixel coarsening}

The outcome of uniform pixel coarsening is shown in Fig.~\ref{fig:Uniform-Coarsening-Total}, for the case of introducing new interphases in (a)--(f), for the phase-preserving case in (g)--(l). 
Since in the first case the color code of the coarsened pixel is the volume average of the contributing four pixels, it preserves the volume average of the Young's moduli. In the limiting case of one single pixel the homogenized Young's modulus equals this volume average. For the present phase fraction ratio of matrix to inclusion of 75.81\% to 24.19\% the obtained Young's modulus is 122.28 MPa in that limiting case.

For the phase-preserving coarsening variant B the 'majority-wins' strategy violates the initial phase ratio. Here, in the limiting case of one single pixel the homogenized elasticity coincides with that of the matrix phase.

The reference solution is obtained on a grid with $h_{\square}=h=1/2048$~mm. The phase distribution coincides for $h_{\square}=1/2048$~mm and $h_{\square}=1/1024$~mm which implies that there is no modeling error for the case of 1024$^2$ pixels.

\begin{Figure}[htbp]
	\subfloat[1024$^2$px]
	{\includegraphics[height=5.1cm,angle=0]{vonMisesZAMM1921}}
	\hspace*{0.01\linewidth}
	\subfloat[64$^2$px]
	{\includegraphics[height=5.1cm, angle=0]{vm_ZAMM_Coarse_Mesh4}}
	\hspace*{0.01\linewidth}
	\subfloat[32$^2$px]
	{\includegraphics[height=5.1cm, angle=0]{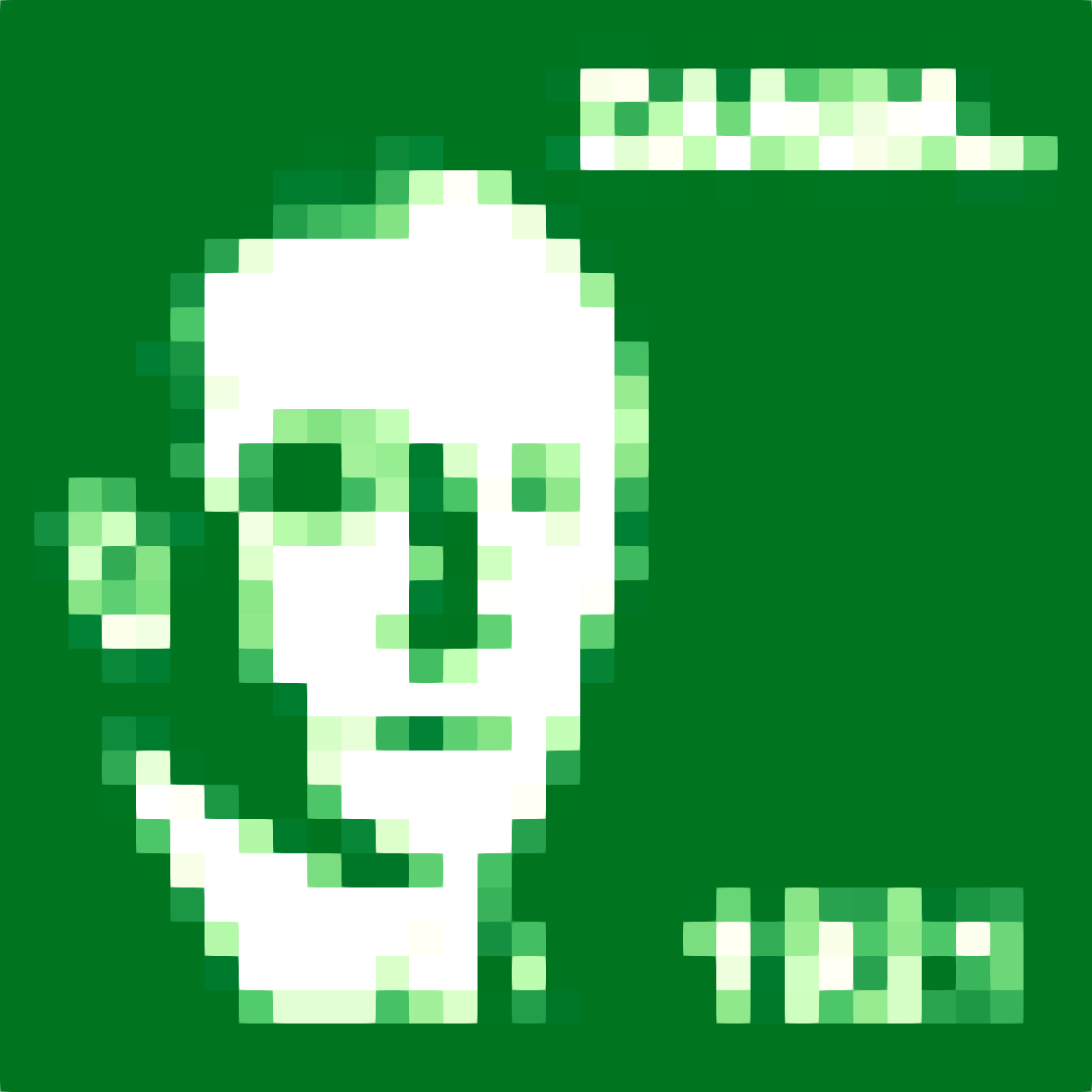}}
	\\[2mm]  
	\subfloat[16$^2$px]
	{\includegraphics[height=5.1cm,angle=0]{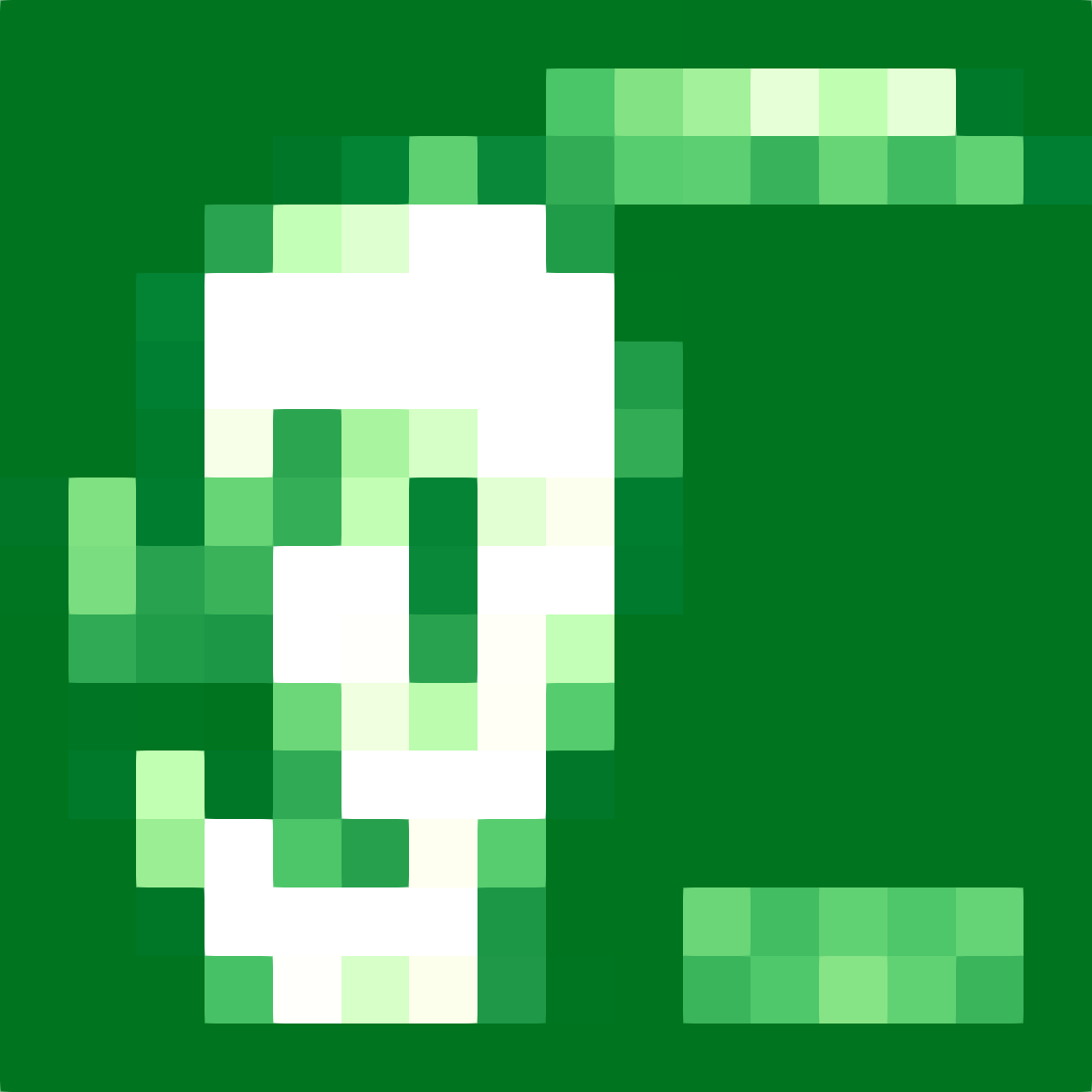}}
	\hspace*{0.01\linewidth}
	\subfloat[8$^2$px]
	{\includegraphics[height=5.1cm, angle=0]{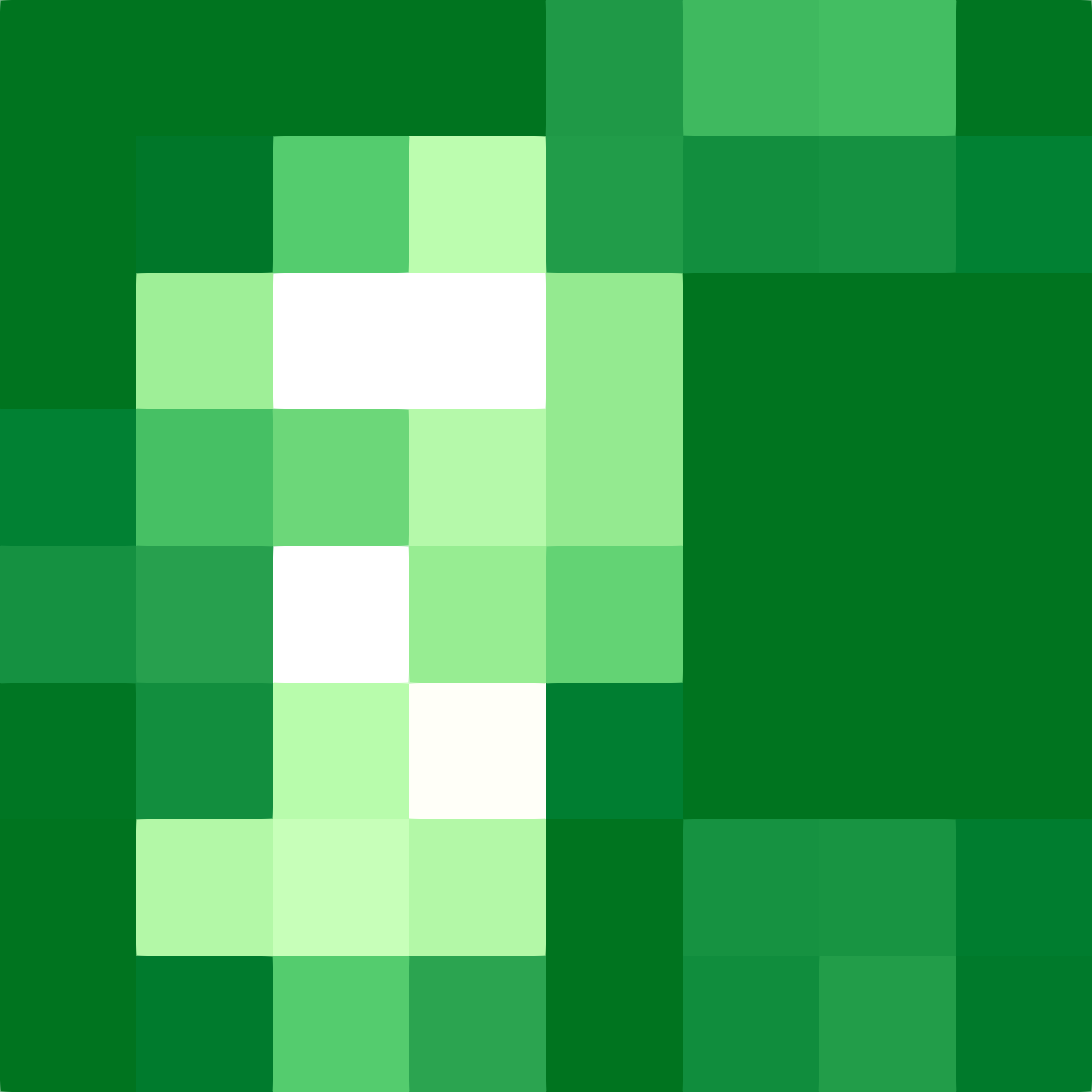}}
	\hspace*{0.01\linewidth}
	\subfloat[1$^2$px]
	{\includegraphics[height=5.1cm, angle=0]{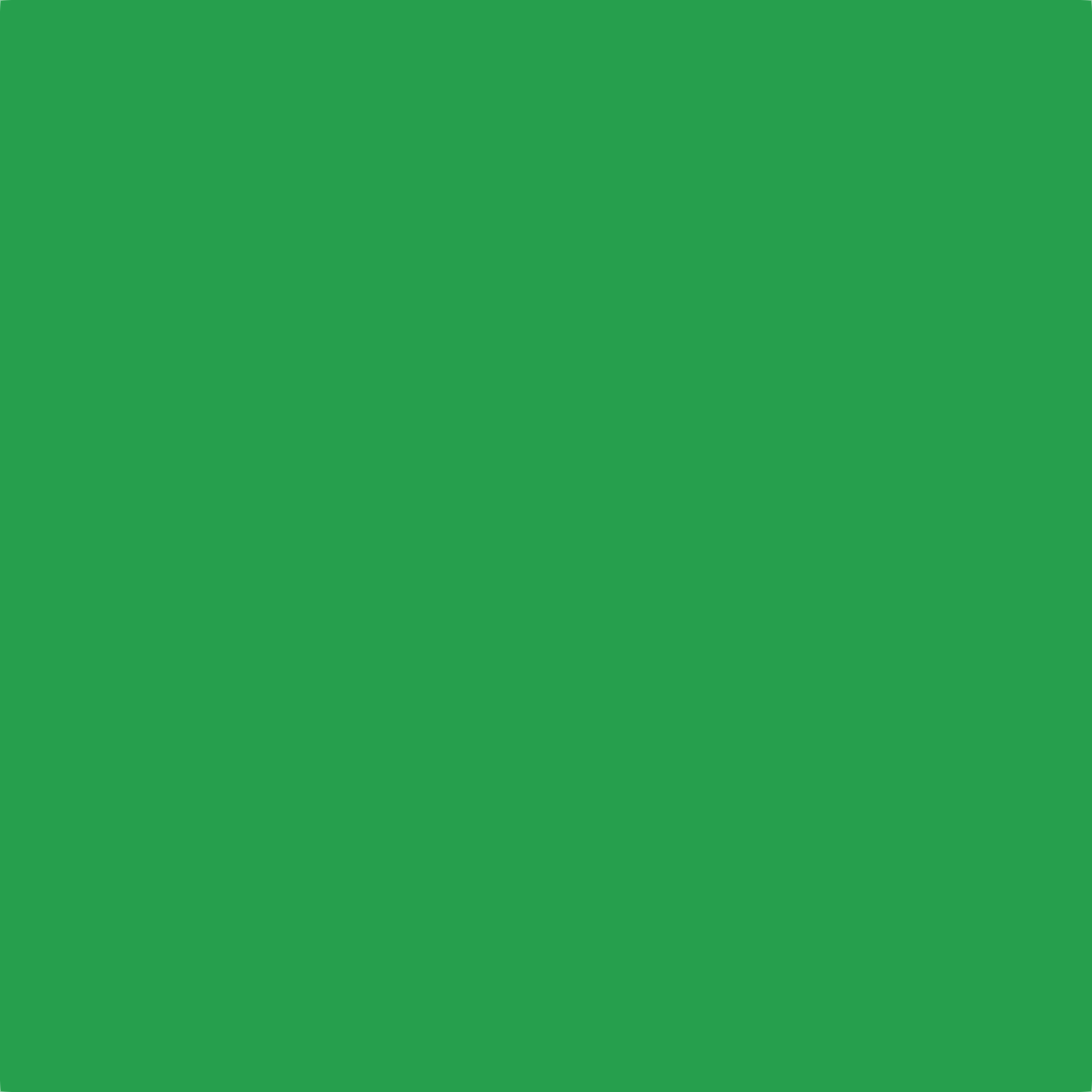}}
	\\[2mm]  
	\subfloat[256$^2$px]
	{\includegraphics[height=5.1cm,angle=0]{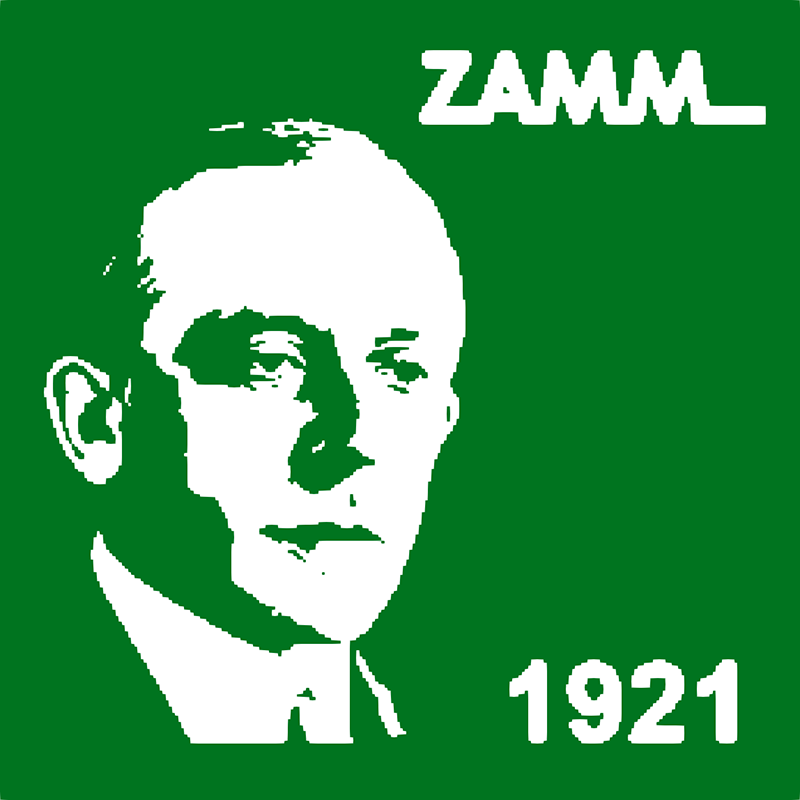}}
	\hspace*{0.01\linewidth}
	\subfloat[64$^2$px]
	{\includegraphics[height=5.1cm, angle=0]{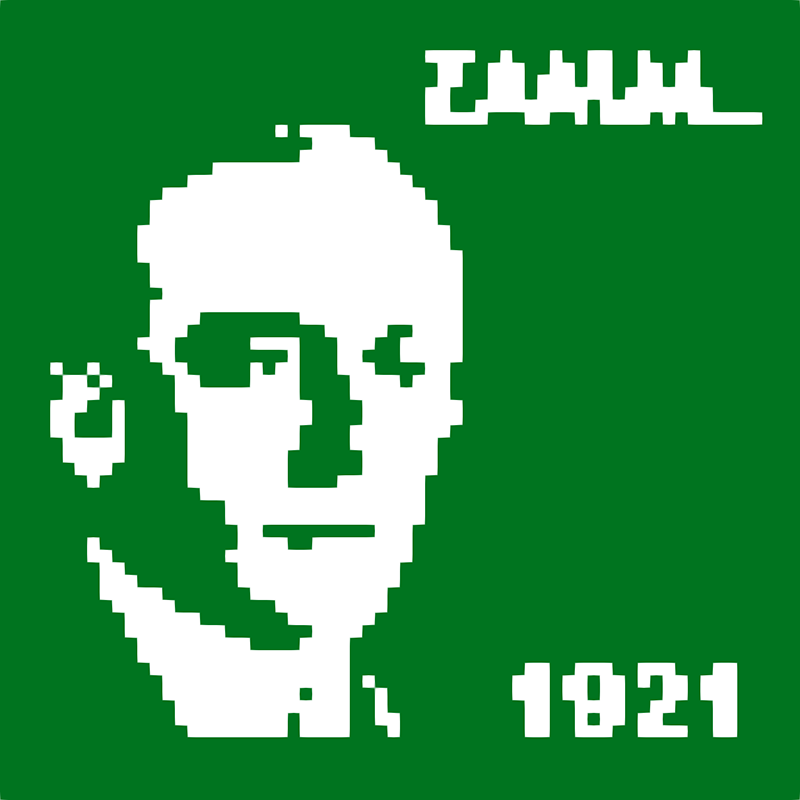}}
	\hspace*{0.01\linewidth}
	\subfloat[32$^2$px]
	{\includegraphics[height=5.1cm, angle=0]{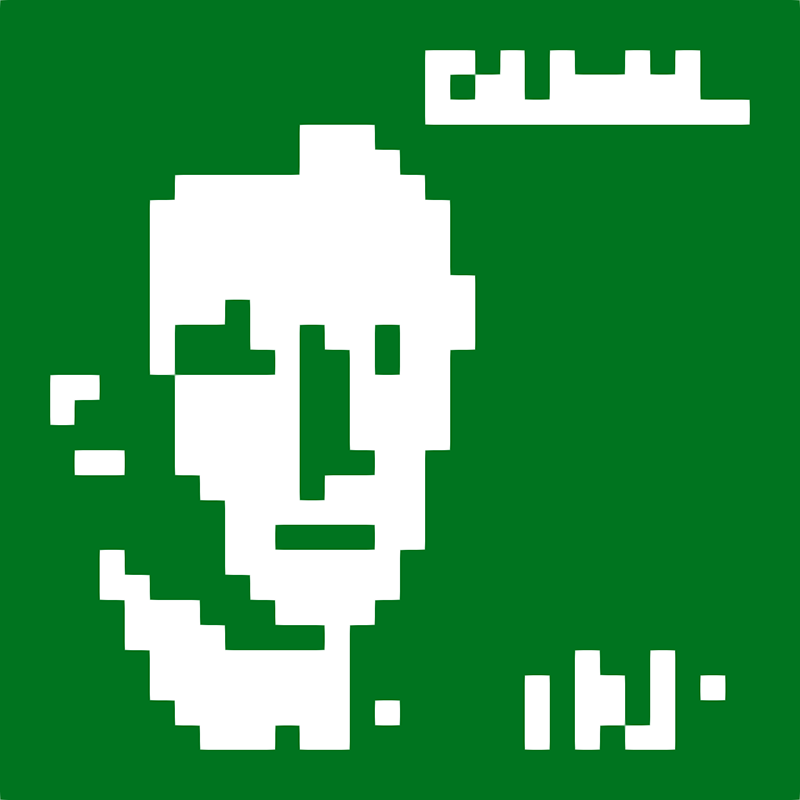}}
	\\[2mm]  
	\subfloat[16$^2$px]
	{\includegraphics[height=5.1cm,angle=0]{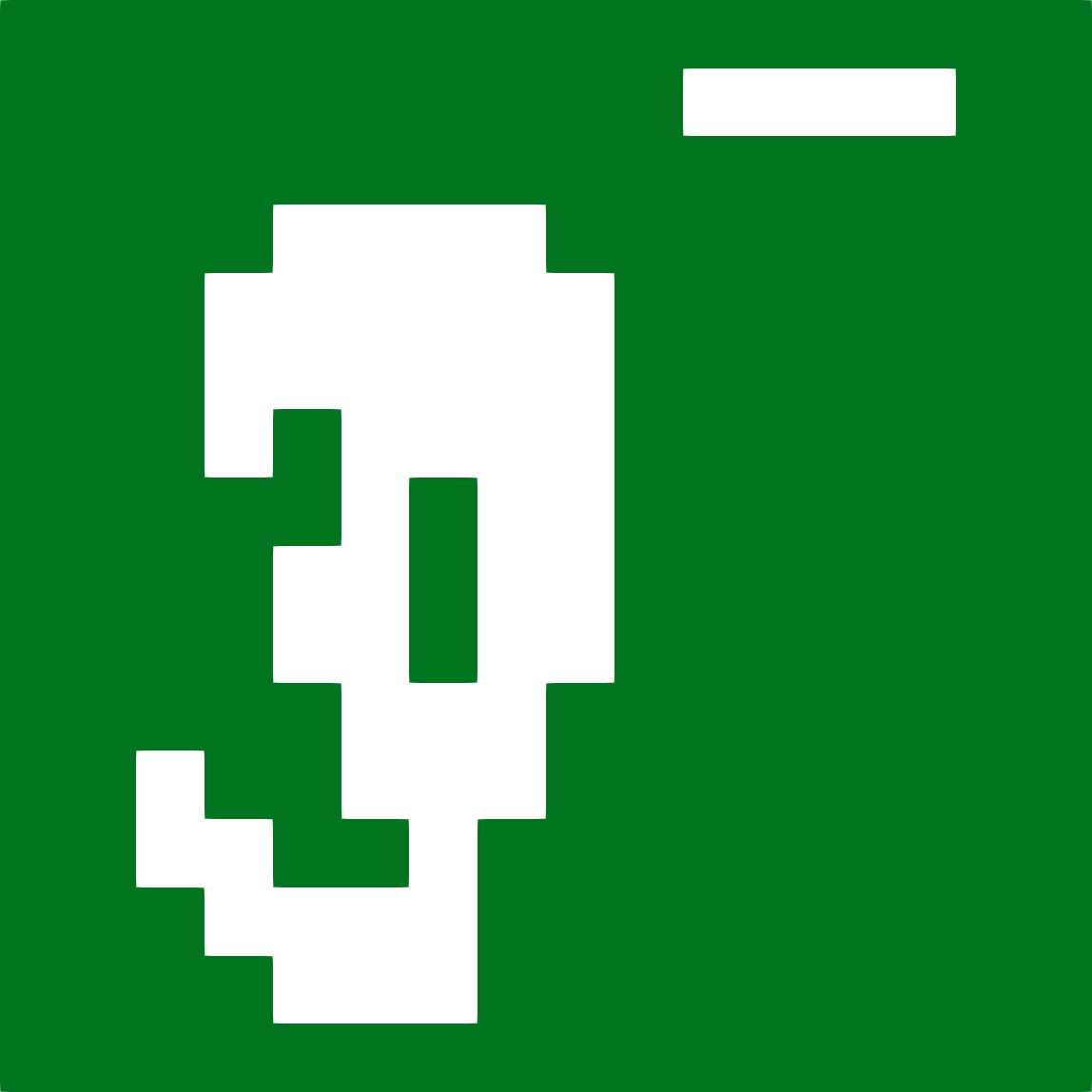}}
	\hspace*{0.01\linewidth}
	\subfloat[8$^2$px]
	{\includegraphics[height=5.1cm, angle=0]{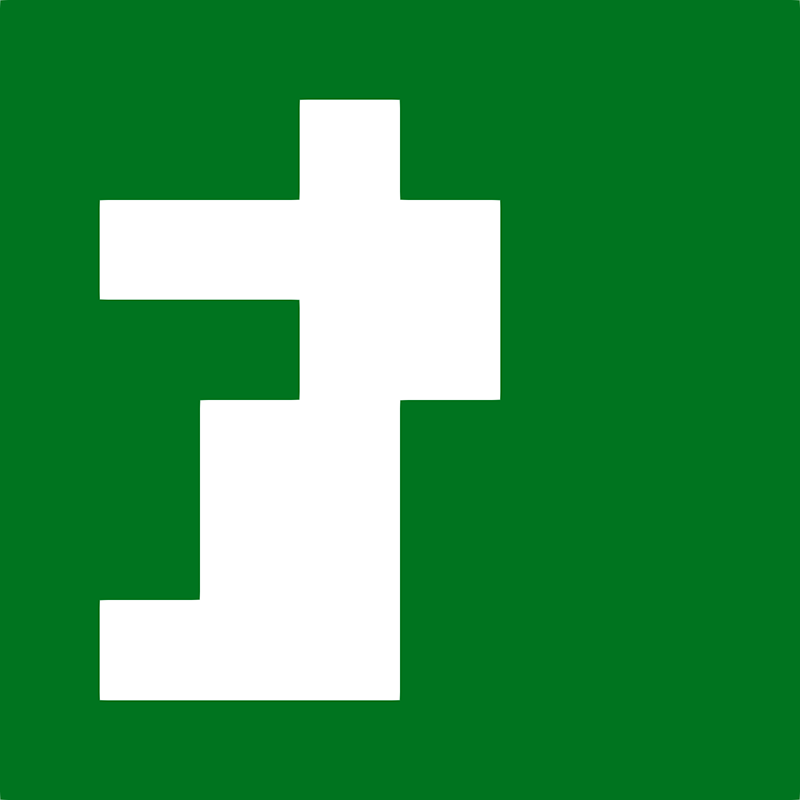}}
	\hspace*{0.01\linewidth}
	\subfloat[1$^2$px]
	{\includegraphics[height=5.1cm, angle=0]{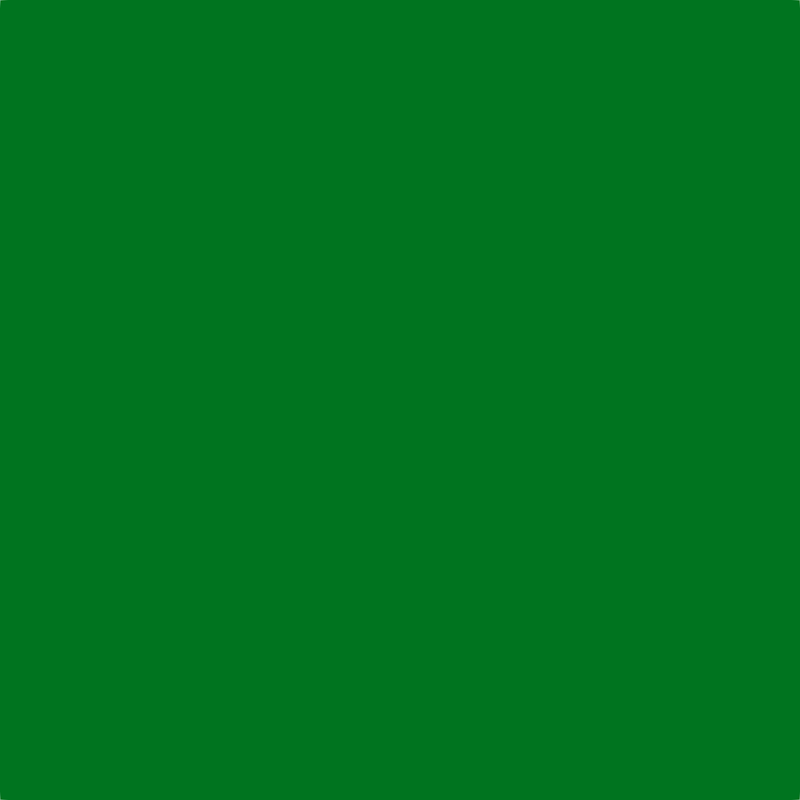}}
	%		%
	\caption{{\bf Pixel coarsening in two variants:} (a)--(f) variant A with additional interphases and (g)--(l) variant B preserving the two-phase composition.   
		\label{fig:Uniform-Coarsening-Total}}
\end{Figure} 

\begin{Table}[htbp]
\begin{minipage}{16.5cm}  
\footnotesize
\centering
\renewcommand{\arraystretch}{1.2} 
\begin{tabular}{c r c c c c c c c}
\hline
&  & \multicolumn{7}{c}{uniform mesh coarsening} \\
\multicolumn{2}{r}{step no.} & 0  & 1  & 2  & 3  & 4 & 5 &    \\
\hline
& px & 1024$^2$ & 512$^2$ & 256$^2$ & 128$^2$ & 64$^2$ & 32$^2$ &   \\  
& ndof & $2\,097\,152$ & $524\, 288$ & $131\, 072$ & $32\, 768$ & $8\, 192$ & $2\,048$ &   \\
& factor & $1.0000$ & $0.2500$ & $0.0625$ & $0.0156$ & $0.0039$ & $0.0010$ &   \\
%			& time (s) & $ $ & $ $ & $ $ & $ $ & $ $ & $ $ & \\
%			& speedup  & $1.0$ & $ $ & $ $ & $ $ & $ $ & $ $ & \\
\hline
&  & \multicolumn{7}{c}{phase-preserving} \\
\hline
& $e^{\epsilon}_{\text{mic}}$ & $ 7.8277 $ & $23.1379$ & $32.8531$ & $45.8962$ & $62.5401$ & $83.1333$ &   \\
%			& factor & $ $ & $1.xxxx$  & $1.xxxx$  & $1.xxxx$  & $1.xxxx$  & $1.xxxx$ & \\
& $\bar e^{\epsilon \, h}_{\text{mic}}$ & $7.5989 $ & $10.7376$ & $15.1268$ & $21.1626$ & $28.8471$ & $36.7235$ &   \\
%			& factor & $ $ & $1.xxx$  & $1.xxx$  & $1.xxx$  & $1.xxx$  & $1.xxx$ & \\
& $e^{\epsilon \, h}_{\text{mic}}$ & $ 7.8277 $ & $12.6256 $ & $18.1766 $ & $25.4883 $ & $34.2359 $ &$41.7658$ &    \\
%			& factor & $ $ & $        $ & $ $ & $ $ & $ $ & & \\
& $\Theta$ & $0.9708 $ & $0.8505 $  & $0.8322 $ & $0.8303 $ & $0.8426 $ &$0.8793$ &  \\   
& $e^{\epsilon \, \square}_{\text{mic}}$ & $ 0.00 $ & $21.4106 $ & $ 30.0302$ & $41.5990 $ & $56.6115 $ & $76.6315$ &   \\
%			& factor & $ $ & $ $ & $ $ & $ $ & $ $ &  & \\[2mm]
\hline
&  & \multicolumn{7}{c}{additional interphases} \\
\hline 
& $e^{\epsilon}_{\text{mic}}$ 		     & $ 7.8277 $ & $16.5271$ & $23.9076$ & $33.5541$ & $46.8402$ & $62.2162$ &   \\
& {\color{black}$\bar e^{\epsilon \, h}_{\text{mic}}$(quad)} & $7.5989 $ & $6.2388$ & $6.6837$ & $8.2921$ & $ 9.3745$  & $ 11.2591$ &   \\
& {\color{black}$\bar e^{\epsilon \, h}_{\text{mic}}$(ave)} & $10.7730$ & $11.9762$ & $15.8392$ & $21.9222$ & $30.0643$ & $37.1357$ &   \\
& $e^{\epsilon \, h}_{\text{mic}}$ & $7.8277 $ & $8.8569 $ & $11.5926 $ & $15.8114 $ & $21.0864 $ &$24.9076$ &   \\
& {\color{black}$\Theta$(quad)} & $0.9708 $ & $0.7044 $  & $0.5766 $ & $0.5244 $ & $0.4446 $ & $0.4520 $ &   \\   
& {\color{black}$\Theta$(ave)} & $1.3763$ & $1.3522$ & $1.3663$ & $1.3865$ & $1.4258$ & $1.4909$ &   \\  
& $e^{\epsilon \, \square}_{\text{mic}}$ & $ 0.00 $ & $14.4814 $ & $ 21.3185$ & $29.9109 $ & $41.8428 $ & $56.8764$ &  \\
\hline
\end{tabular} 
\end{minipage}
\caption{\textbf{Pixel coarsening, phase-preserving:} for various resolutions the number of degrees of freedom (ndof), ndof reduction factors, different errors along with their increase factor, effectivity index $\Theta$. Error data in the energy norm in $10^{-4}$\,(MPa). 
}
\label{tab:vonMises-Ndof-Error-speedup-UNIFORM} 
\end{Table} 

\begin{Figure}[htbp]
	\centering
    \subfloat[additional interphases]
	{\resizebox{0.38\columnwidth}{!}{\includegraphics{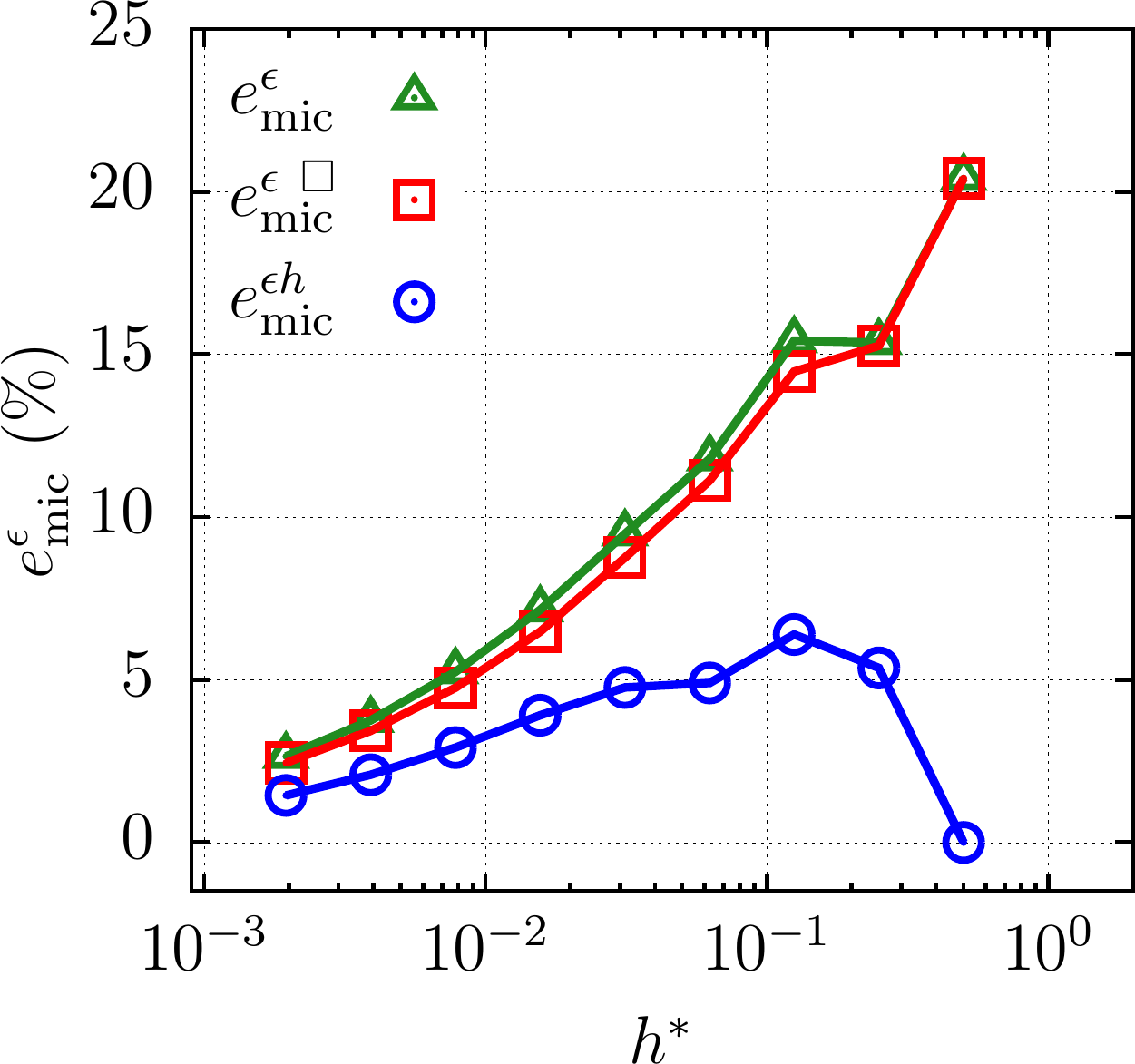}}}
		\hspace*{0.01\linewidth}
	\subfloat[phase-preserving]
	{\resizebox{0.38\columnwidth}{!}{\includegraphics{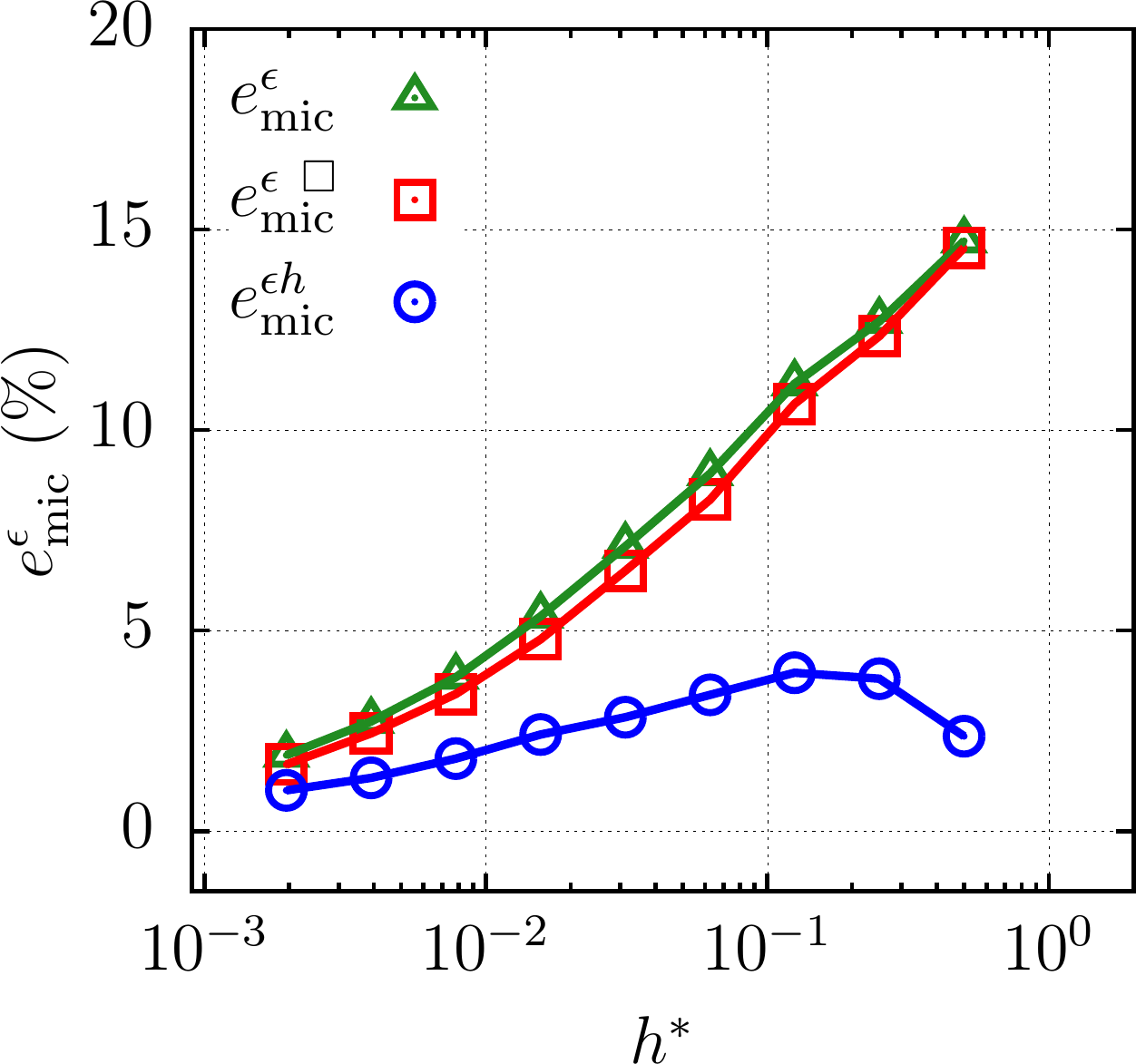}}}		
	\caption{\textbf{Pixel coarsening, phase-preserving:} relative percental micro errors, total $e^{\epsilon}_{mic}/||\bm u||_{A(\mathcal{B}_{\epsilon})}$, modeling $e^{\epsilon\, \square}_{mic}/||\bm u||_{A(\mathcal{B}_{\epsilon})}$, and discretization $e^{\epsilon\, h}_{mic}/||\bm u||_{A(\mathcal{B}_{\epsilon})}$ versus characteristic element length $h^{*}$, resolution from 512$^2$ px to 2$^2$ px.}
	\label{fig:err-PhasePreserving}	
\end{Figure}
Results are displayed in Tab.~\ref{tab:vonMises-Ndof-Error-speedup-UNIFORM} and in Fig.~\ref{fig:err-PhasePreserving}. We use $h^{*} := \sqrt{1/\text{num}_{elem}}$ with num$_{elem}$ the number of elements in a unit cell for the inclusion of results for nonuniform meshes. For square microdomains with uniform meshes it holds $h^{*}=h$. The micro error and its parts monotonously increase for both coarsening variants with the exception of the discretization error in the last two steps of extremely coarse resolutions, which are in either case not competitive. The modeling error is larger than the discretization error, a gap that continuously increases with coarsening. Inequality \eqref{MicroErrorTriangleInequality} is fulfilled. Error estimation for coarsening with new interphases along with phase-distinction in stress computation at interfaces leads to poor results, see $\Theta$(quad) in Tab.~\ref{tab:vonMises-Ndof-Error-speedup-UNIFORM}. Averaging stress alleviates this issue of error computation already addressed in Sec.~\ref{sec:error_estimation_average}.    

\begin{Table}[htbp]
	\begin{minipage}{16.5cm}  
		\footnotesize
		\centering
		\renewcommand{\arraystretch}{1.2} 
		\begin{tabular}{r l l l l l}
			\hline
			& vol.\% & $\mathbb{A}_{11}$ & $\mathbb{A}_{22}$ & $\mathbb{A}_{33}$ & $\mathbb{A}_{12}$ \\
			\hline
			inclusion phase	 & 24.2 & 213.44 & 213.44 & 80.04 & 53.36 \\
			homogenized	 &      & 127.99 & 129.54 & 47.82 & 32.13 \\  
			matrix phase   & 75.8 & 111.11 & 111.11 & 41.67 & 27.78 \\
			\hline	
		\end{tabular} 
	\end{minipage}
	\caption{\textbf{Elastic coefficients}: for the two phases and the two-phase composite in (MPa).}
	\label{tab:ElasticCoefficients} 
\end{Table} 

\begin{Figure}[htbp]
	\centering
	\subfloat[additional interphases]
	{\resizebox{0.38\columnwidth}{!}{\includegraphics{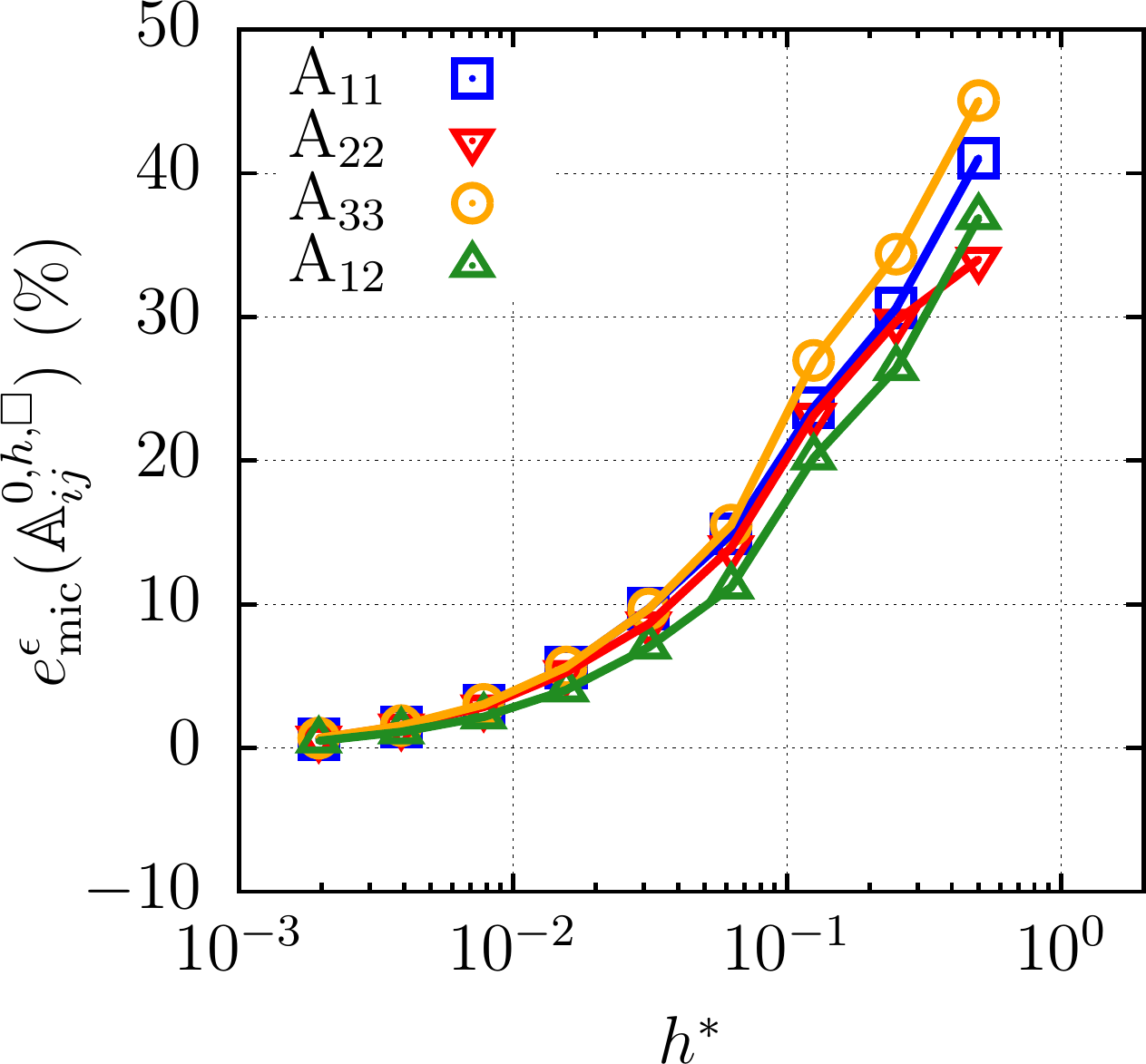}}}
	\hspace*{0.01\linewidth}
	\subfloat[phase-preserving]
	{\resizebox{0.38\columnwidth}{!}{\includegraphics{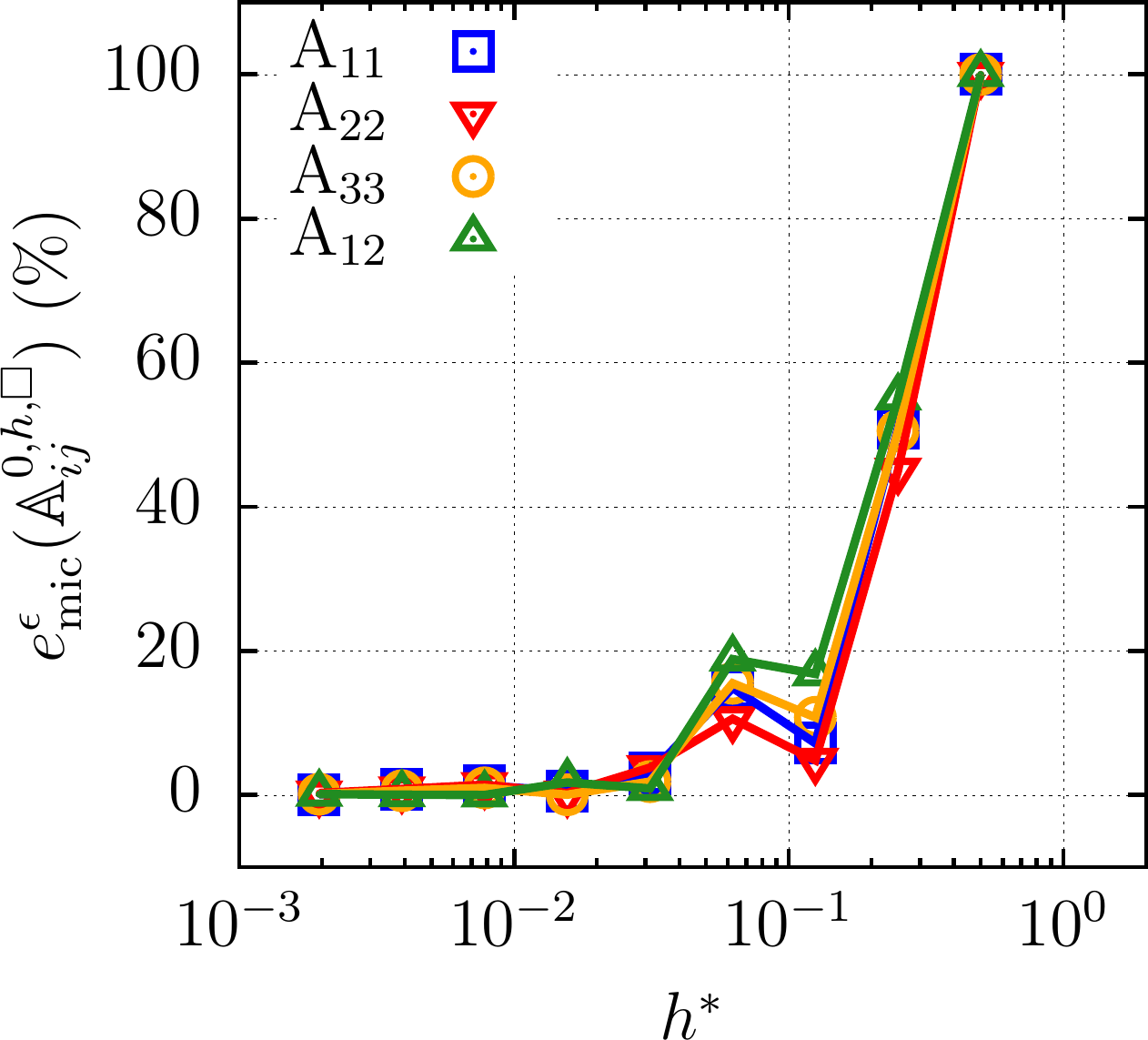}}}
	\caption{\textbf{Pixel coarsening:} relative total micro error $e^{\epsilon}_{\text{mic}}(\mathbb{A}^{0,h,\square}_{ij})$ versus characteristic element length $h^{*}$, resolution from 512$^2$ px to 2$^2$ px.}
	\label{fig:err-Aij-both-coarsening-variants}	
\end{Figure}

\begin{Figure}[htbp]
	\centering
	\subfloat[additional interphases, from 512$^2$ px to 2$^2$ px]
	{\resizebox{0.95\columnwidth}{!}{\includegraphics{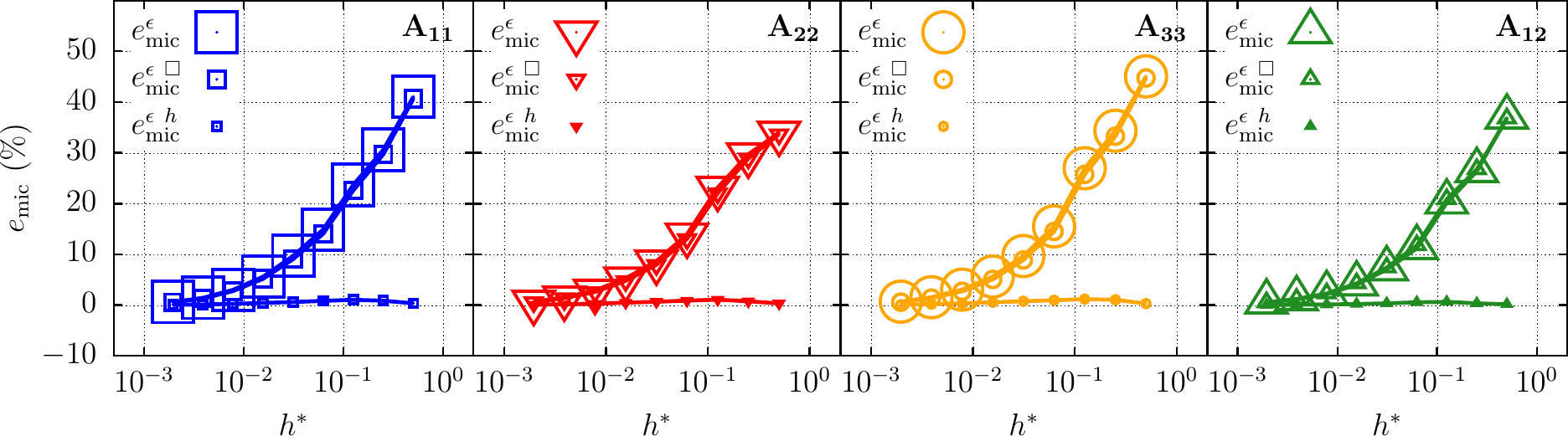}}}
    \\[2mm]
	\subfloat[phase-preserving, from 512$^2$ px to 32$^2$ px]
	{\resizebox{0.95\columnwidth}{!}{\includegraphics{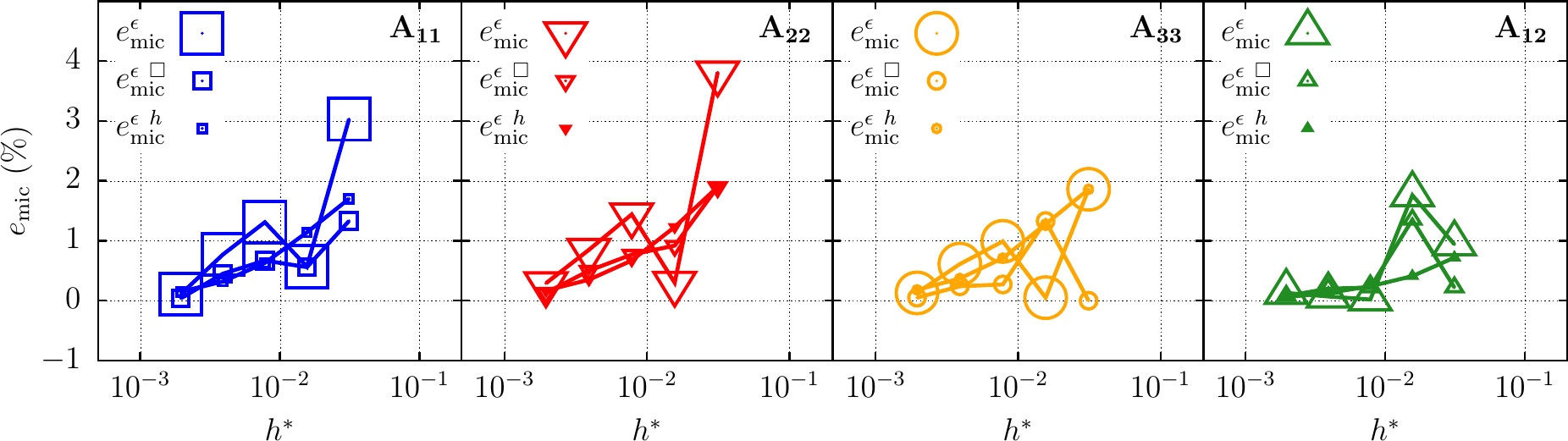}}}
	\caption{\textbf{Pixel coarsening:} relative micro errors for $\mathbb{A}_{ij}$, total $e^{\epsilon}_{mic}$, modeling $e^{\epsilon\, \square}_{mic}$, and discretization $e^{\epsilon\, h}_{mic}$ versus characteristic element length $h^{*}$.}
	\label{fig:error-Aij-decomposition-both-coarsening-variants}
\end{Figure}

A second measure of the total micro error is the error in the coefficients of the homogenized elasticity tensor displayed in Fig.~\ref{fig:err-Aij-both-coarsening-variants}. In Fig.~\ref{fig:error-Aij-decomposition-both-coarsening-variants} the decomposition into modeling and discretization error is provided. The errors are computed according to \eqref{eq:total-micro-error-at-microscale}--\eqref{eq:micro-modeling-error-at-microscale} in the adaption to the case of the elasticity coefficients \eqref{eq:ConvergenceElasticities}. Exact values of the homogenized elasticities for a 2048$^2$ px resolution are listed in Tab.~\ref{tab:ElasticCoefficients}. 

For the coarsening with additional interphases the micro error is continuously increasing in all components of the homogenized elasticity tensor, thereby overestimating the true stiffness, see Fig.~\ref{fig:err-Aij-both-coarsening-variants} (a). The modeling error is predominant in comparison to the discretization error, the latter is bounded to less than 2\% in the first 5  stages of coarsening Fig.~\ref{fig:error-Aij-decomposition-both-coarsening-variants} (a).

For the phase-preserving coarsening variant the micro error is throughout smaller in the first five coarsening steps (down to 32$^2$ px) than for the coarsening with new interphases Fig.~\ref{fig:err-Aij-both-coarsening-variants} (b). The error continuously increases, coarsening step no. 4 is an outlier. Modeling and discretization error are roughly of the same magnitude Fig.~\ref{fig:error-Aij-decomposition-both-coarsening-variants} (b).

\begin{Table}[htbp]
	\begin{minipage}{16.5cm}  
		\footnotesize
		\centering
		\renewcommand{\arraystretch}{1.2} 
		\begin{tabular}{r c c c c c c c c c c }
			\hline
			%			step no. &0 & 1  & 2  & 3  & 4 & 5 & 6 & 7 & 8 & 9  \\
			%			\hline
			&  & \multicolumn{9}{c}{uniform pixel coarsening} \\ 
			px & 2048$^2$   & 512$^2$ & 256$^2$ & 128$^2$ & 64$^2$ & 32$^2$ & 16$^2$ & 8$^2$ & 4$^2$ & 2$^2$    \\  
			%			ndof &  & $528\,390$ & $133\,126$ & $33\,798$ & $8\,710$ & $2\,310$ & $646$ & $198$ & $70$ & $30$  \\
			\hline
            &  & \multicolumn{9}{c}{additional interphases} \\
            \hline
            $u_{\text{max}}$ 
            & 86.57 & 86.50	& 86.40	& 86.22 	&85.92	& 85.45	& 84.86	& 83.90 	& 83.15 	& 82.11 \\
            max\{$\sigma_{\text{von Mises}}$\}  
            & 	1.930 & 1.652	& 1.553	& 1.473	& 1.377	& 1.285	& 1.112	& 1.034	& 1.011	& 0.987  \\
            \hline
			&  & \multicolumn{9}{c}{phase-preserving} \\
			\hline
			$u_{\text{max}}$ 
			& 86.57 & 86.56	& 86.48	& 86.40 	& 86.48	& 86.17	& 88.22	& 87.25 	& 92.68 	& 99.74 \\
			max\{$\sigma_{\text{von Mises}}$\}  
			& 	1.930	& 1.713		& 1.595		& 1.497		& 1.38		& 1.356		& 1.067		& 1.04		& 0.908		& 0.791  \\
			\hline
		\end{tabular} 
	\end{minipage}
	\caption{\textbf{Pixel coarsening:} maximum deflection $u_\text{max}$ in ($\text{mm}$) of cantilever beam and maximum value of von Mises stress in ($\text{N/mm$^2$}$) for macro GP at $x=y=2.1132$~mm.}
	\label{tab:deflection-vm-stress} 
\end{Table}

\begin{Figure}[htbp]
	\subfloat[512$^2$px]
	{\includegraphics[height=4.1cm,angle=0]{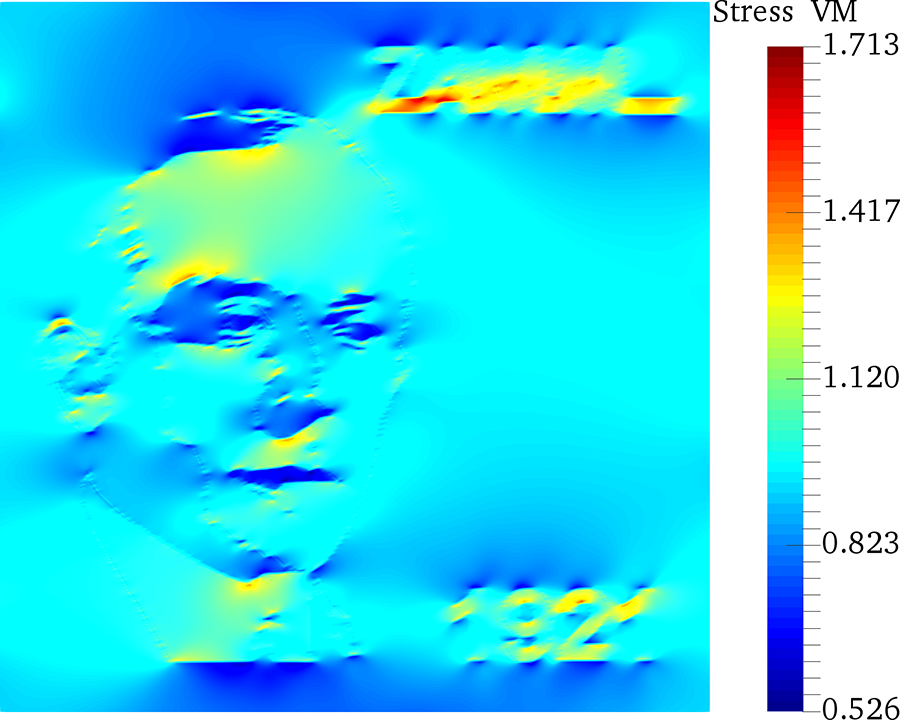}}
	\hspace*{0.01\linewidth}
	\subfloat[128$^2$px]
	{\includegraphics[height=4.1cm, angle=0]{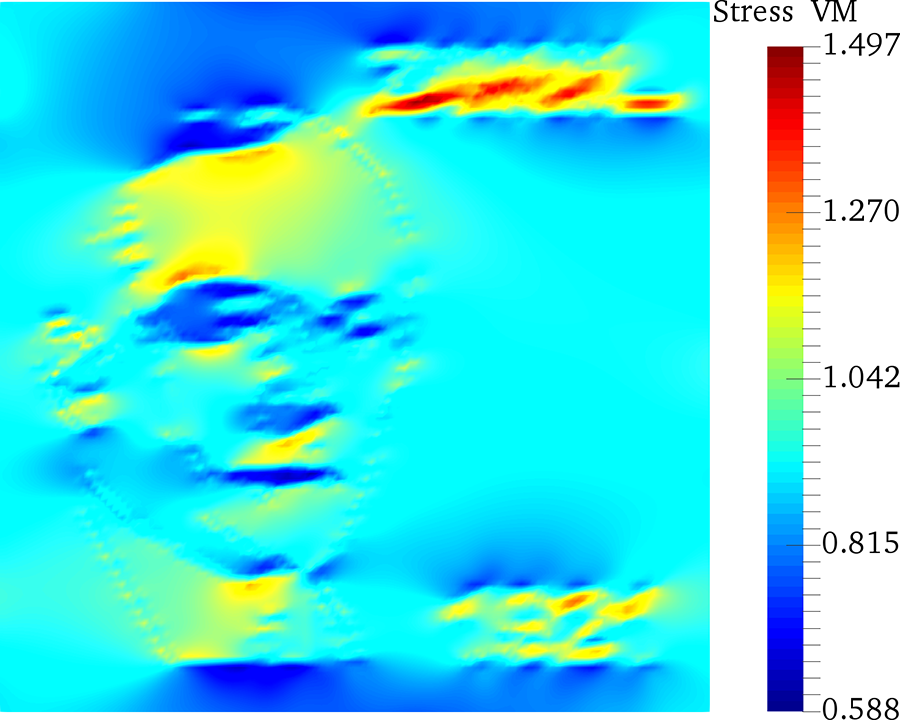}}
	\hspace*{0.01\linewidth}
	\subfloat[32$^2$px]
	{\includegraphics[height=4.1cm, angle=0]{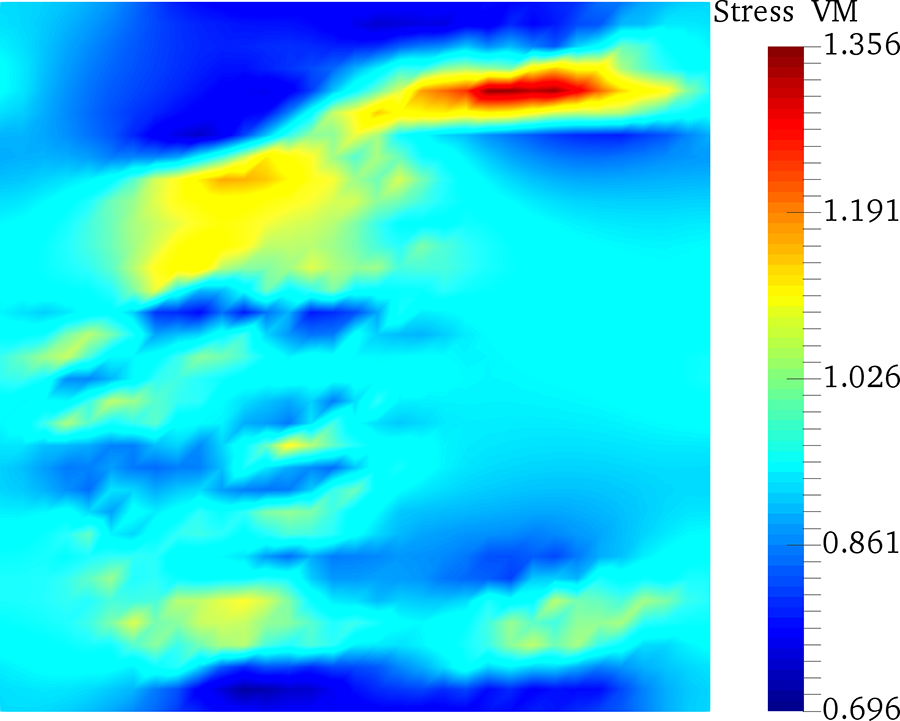}}
%	\\[2mm]
%	\subfloat[512$^2$px]
%	{\includegraphics[height=4.1cm,angle=0]{zamm_interphase_512p_stress_vm}}
%	\hspace*{0.01\linewidth}
%	\subfloat[128$^2$px]
%	{\includegraphics[height=4.1cm, angle=0]{zamm_interphase_128p_stress_vm}}
%	\hspace*{0.01\linewidth}
%	\subfloat[32$^2$px]
%	{\includegraphics[height=4.1cm, angle=0]{zamm_interphase_32p_stress_vm}}
	%
	\caption{\textbf{Von Mises stress} for uniform discretizations from phase-preserving coarsening.  
	% (d)--(f) for new interphases. 
	\label{fig:vm-stress-phase-preserving}}
\end{Figure} 
 
%The results fulfill the triangle equality \eqref{MicroErrorTriangleInequality}. 
Table \ref{tab:deflection-vm-stress} lists with the maximal deflection of the macrostructure a macro accuracy measure and with the maximal von Mises stress on the chosen microdomain a local micro accuracy measure. The latter becomes relevant if inelastic material behavior comes into play. 

Figure \ref{fig:dicretization-err-rel-energy-norm-2phase} shows the distribution of micro errors and their constituents for the two coarsening variants A and B. Errors are confined to phase boundaries. 
  
\begin{Figure}[htbp]
New interphases, variant A 
    \\[1mm]
	\subfloat[$e^{\epsilon \, \square}_{mic}/||\bm u||_{A(\mathcal{B}_{\epsilon})}$]
	{\includegraphics[height=4.6cm, angle=0]{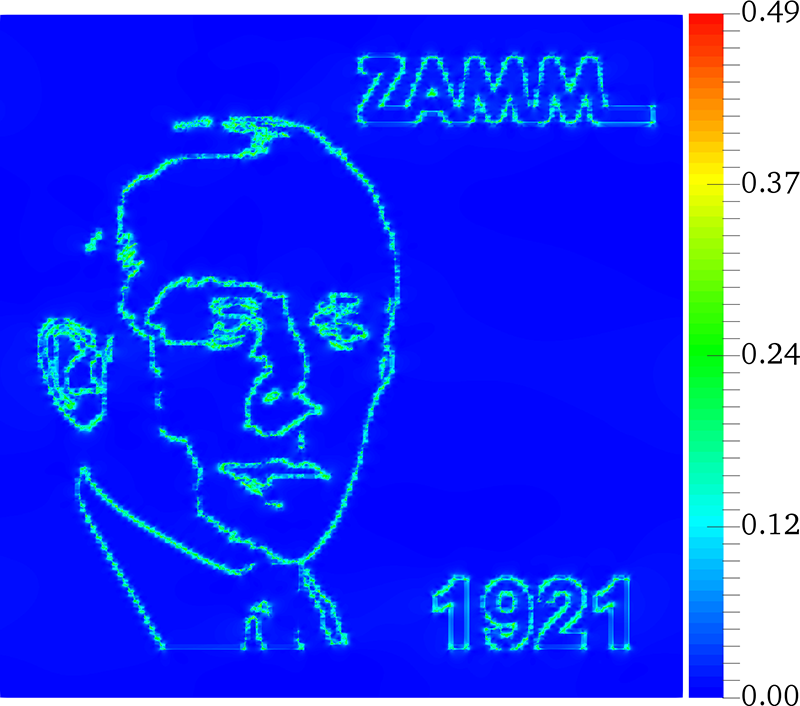}}
	\hspace*{0.01\linewidth}
	\subfloat[$e^{\epsilon \, h}_{mic}/||\bm u||_{A(\mathcal{B}_{\epsilon})}$]
	{\includegraphics[height=4.6cm, angle=0]{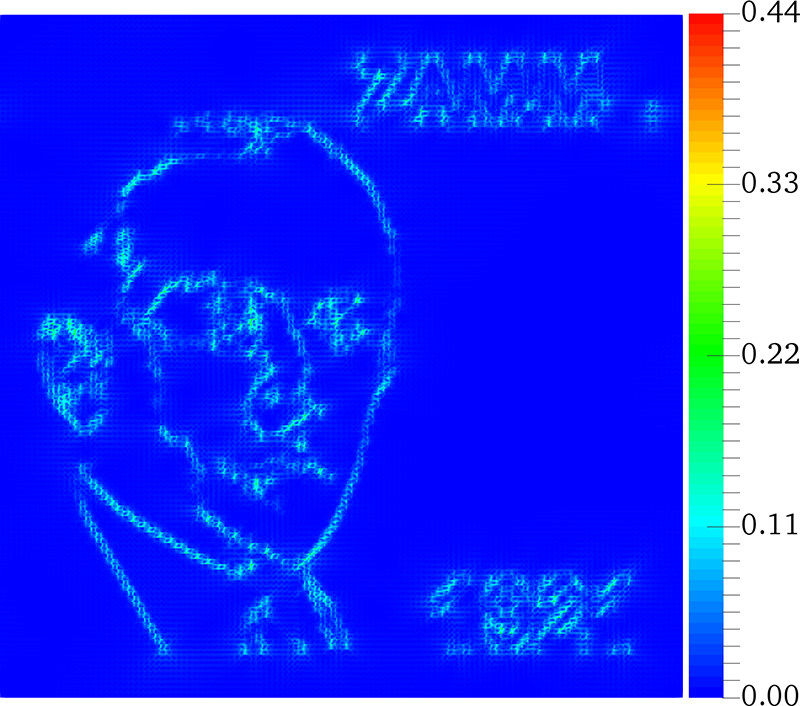}}
	\hspace*{0.01\linewidth}
	\subfloat[$e^{\epsilon}_{mic}/||\bm u||_{A(\mathcal{B}_{\epsilon})}$]
	{\includegraphics[height=4.6cm, angle=0]{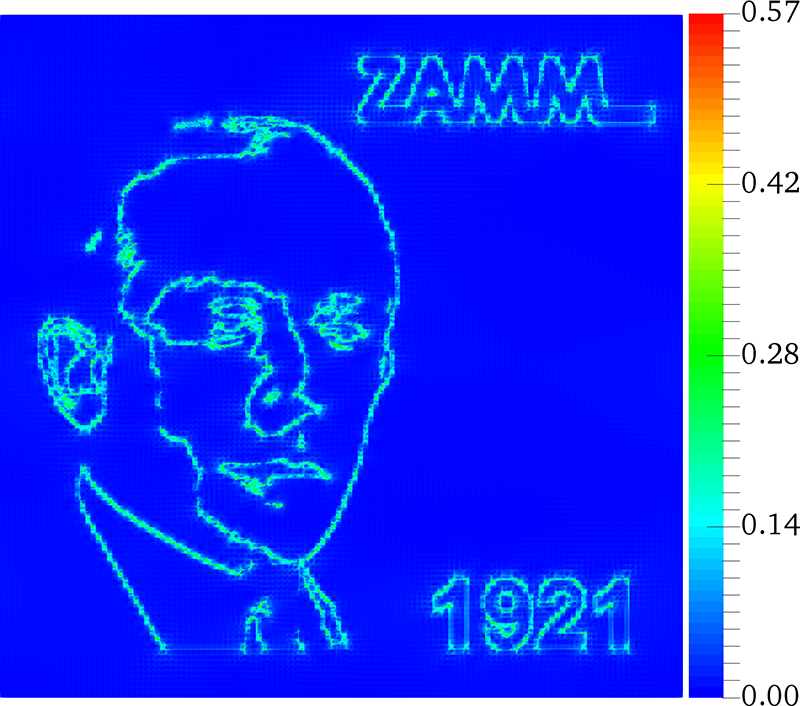}}
\\[4mm]
Phase-preserving, variant B 
\\[0mm]
	\subfloat[$e^{\epsilon \, \square}_{mic}/||\bm u||_{A(\mathcal{B}_{\epsilon})}$]
	{\includegraphics[height=4.6cm, angle=0]{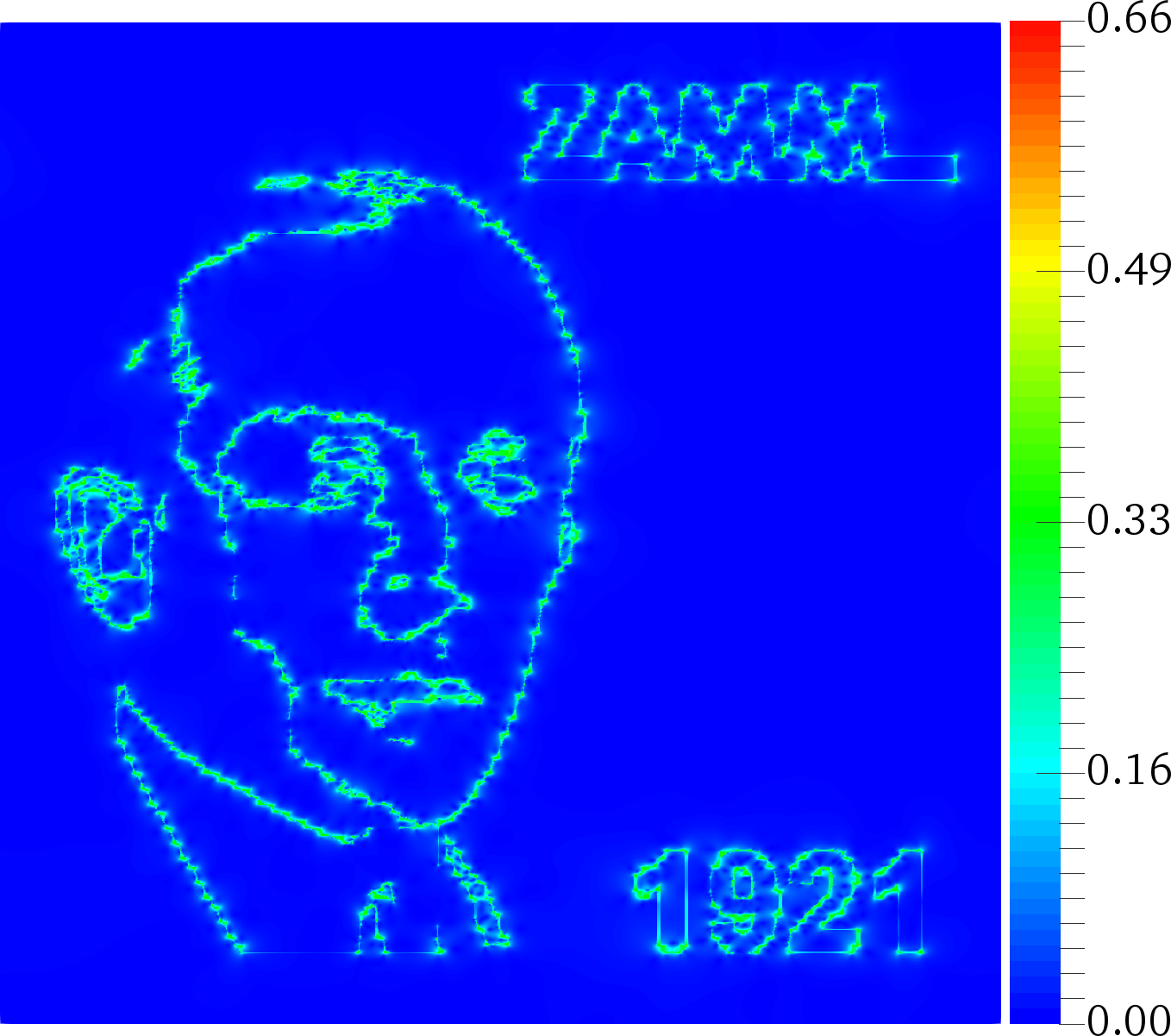}}
	\hspace*{0.01\linewidth}
	\subfloat[$e^{\epsilon \, h}_{mic}/||\bm u||_{A(\mathcal{B}_{\epsilon})}$]
	{\includegraphics[height=4.6cm, angle=0]{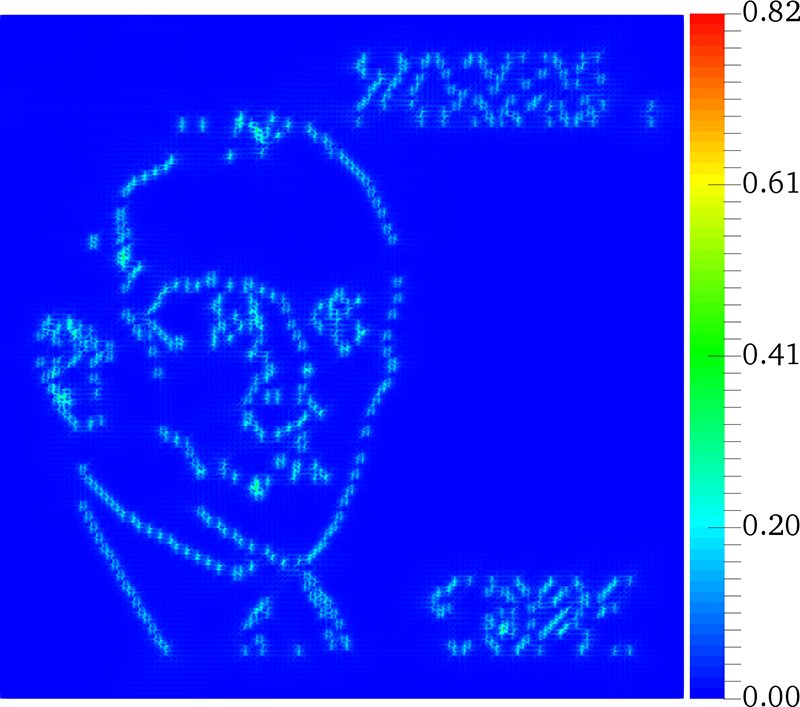}}
	\hspace*{0.01\linewidth}
	\subfloat[$e^{\epsilon}_{mic}/||\bm u||_{A(\mathcal{B}_{\epsilon})}$] 
	{\includegraphics[height=4.6cm, angle=0]{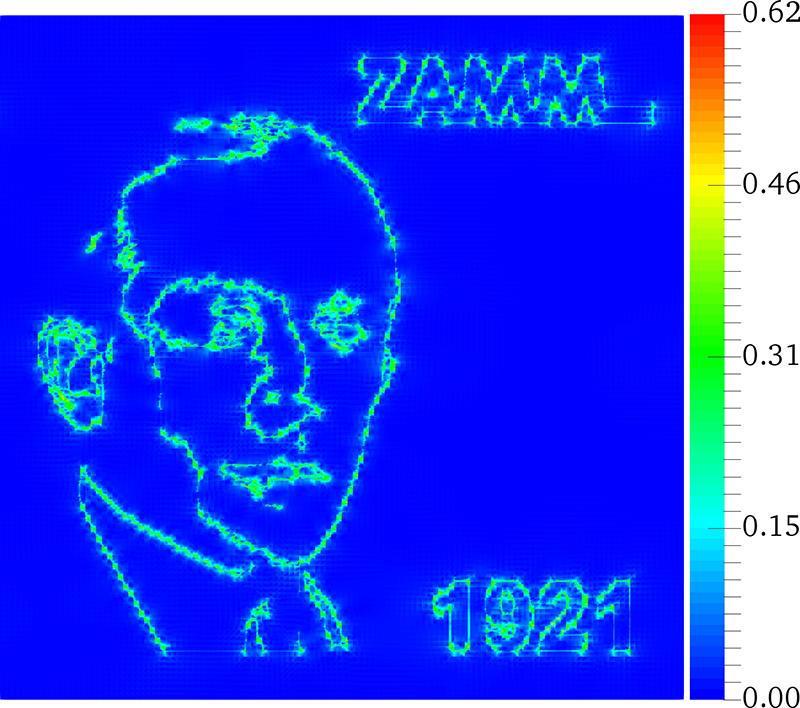}}  
	\caption{\textbf{Error distributions for 128$^2$ px resolution}: (from left to right) modeling error, discretization error, total error.  
   \label{fig:dicretization-err-rel-energy-norm-2phase}}
\end{Figure} 

\subsection{Adaptive mesh-coarsening}

The first step of coarsening was a uniform pixel coarsening changing the image resolution. The second step carries out a consecutive, adaptive, microstructure-guided mesh coarsening based on the uniformly coarsened pixel resolutions. In this second step the discretization at interfaces is maintained for accuracy, but in the phase interior mesh coarsening is carried for efficiency. Consequently, this step introduces an additional discretization error but preserves exactly the existing phase distributions including their boundaries. Three different versions of adaptively coarsened meshes are shown in Fig.~\ref{fig:Adaptive-meshcoarsening-vMises}, which are each based on different uniform pixel resolutions. As a consequence the maximal number of quadtree-type adaptive coarsening steps differs between these versions. 

Table~\ref{tab:AdapCoarsening-Ndof-Error-speedup} indicates the considerable computational savings of adaptive mesh coarsening at moderate accuracy losses. 
Error estimation for phase-preserving coarsening turns out to be much more accurate than coarsening with new interphases as indicated by effectivity indices much closer to one. Notice that the effectivity indices for the different stress computations in error estimation confirm the preliminary results of Sec.~\ref{sec:error_estimation_average}.  

\begin{Figure}[htbp]
%	\subfloat[1024px, adap7]
%	{\includegraphics[height=5.1cm,angle=0]{vonMises_1024p_ZAMM_QT_Mesh7}}
%	\hspace*{0.01\linewidth}
%	\subfloat[128px, adap4]
%	{\includegraphics[height=5.1cm, angle=0]{vonMises_128p_ZAMM_QT_Mesh4}}
%	\hspace*{0.01\linewidth}
%	\subfloat[64px, adap3]
%	{\includegraphics[height=5.1cm, angle=0]{vonMises_64p_ZAMM_QT_Mesh3}}\\
	\subfloat[1024$^2$ px, adap7]
	{\includegraphics[height=5.1cm,angle=0]{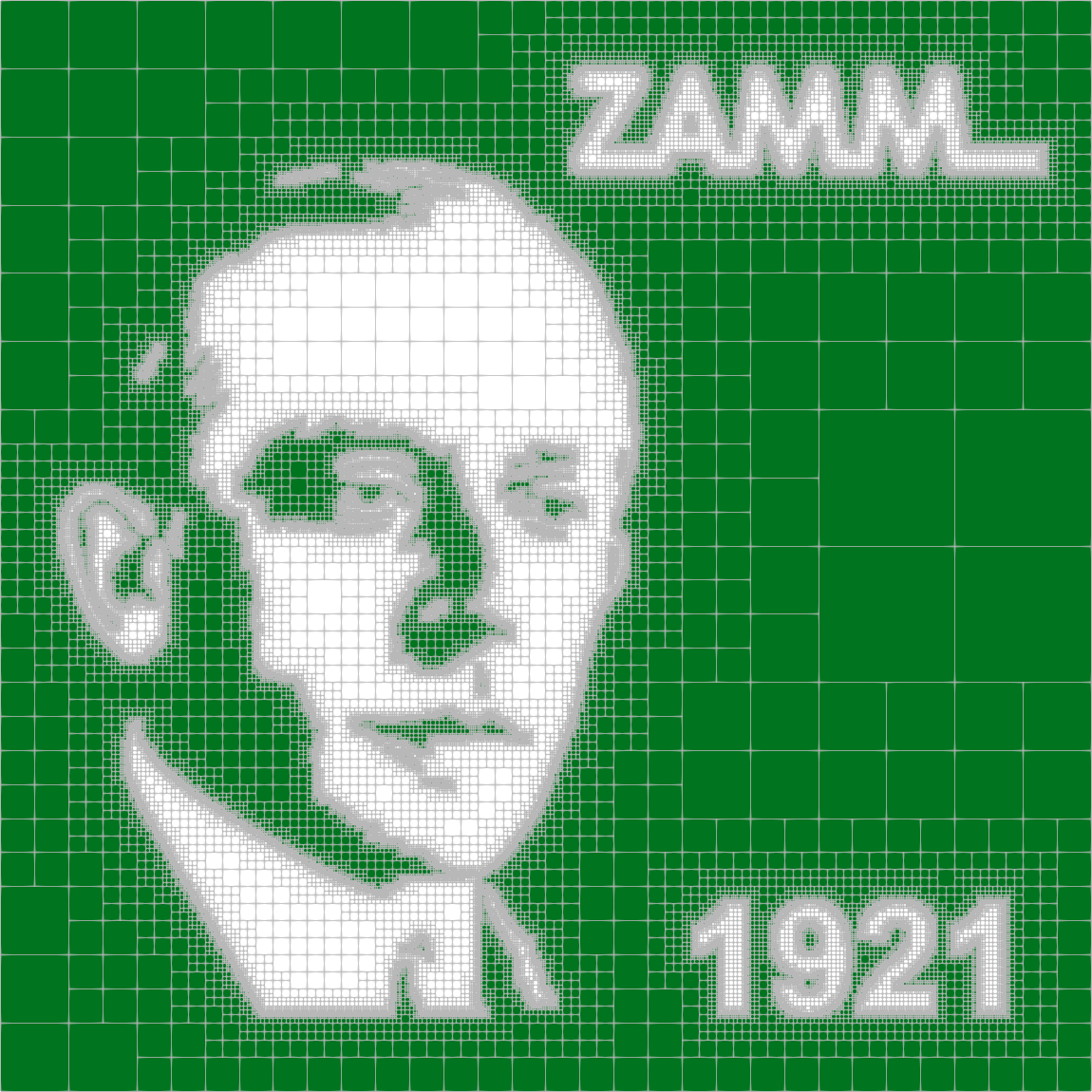}}
	\hspace*{0.01\linewidth}
	\subfloat[128$^2$px, adap4]
	{\includegraphics[height=5.1cm, angle=0]{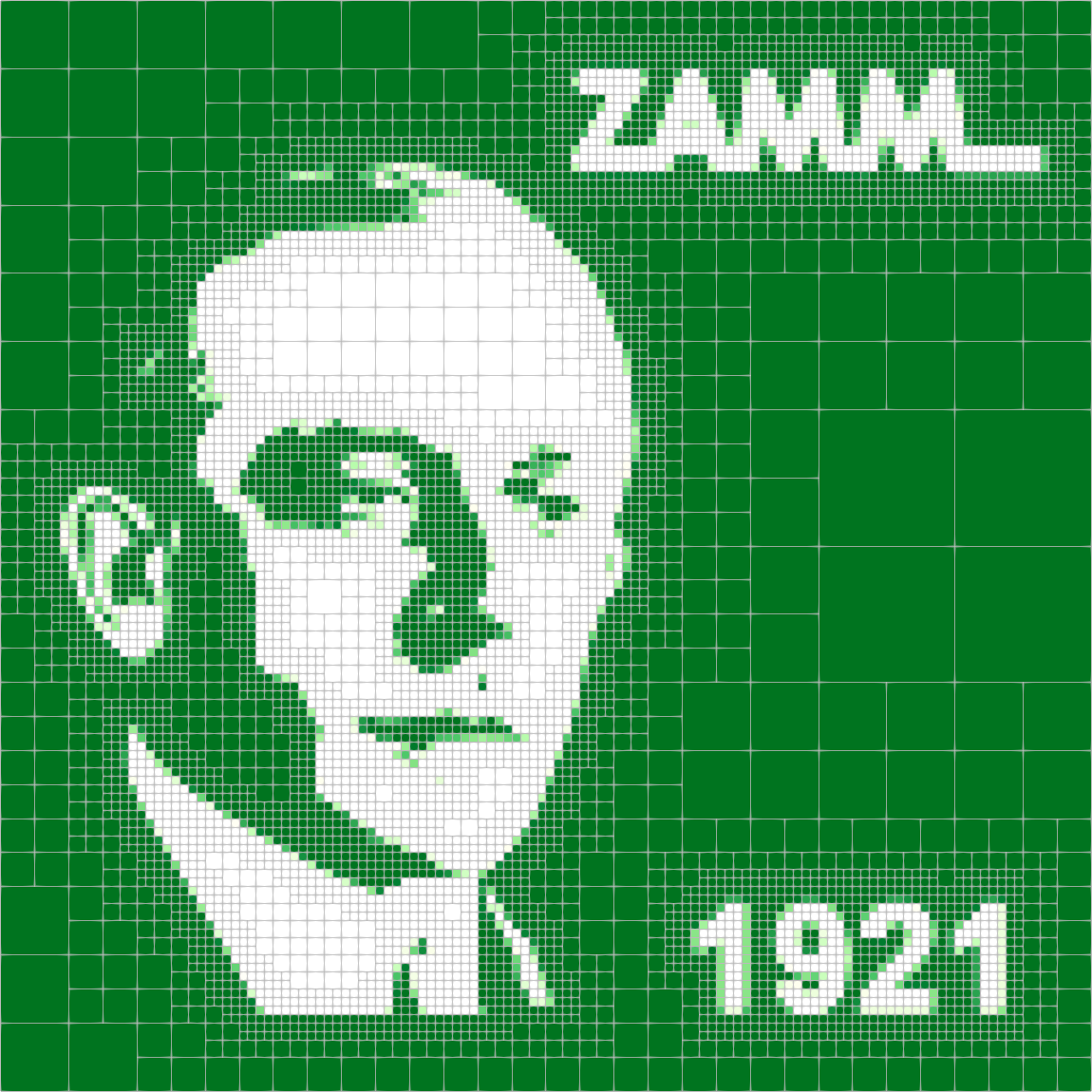}}
	\hspace*{0.01\linewidth}
	\subfloat[64$^2$px, adap3]
	{\includegraphics[height=5.1cm, angle=0]{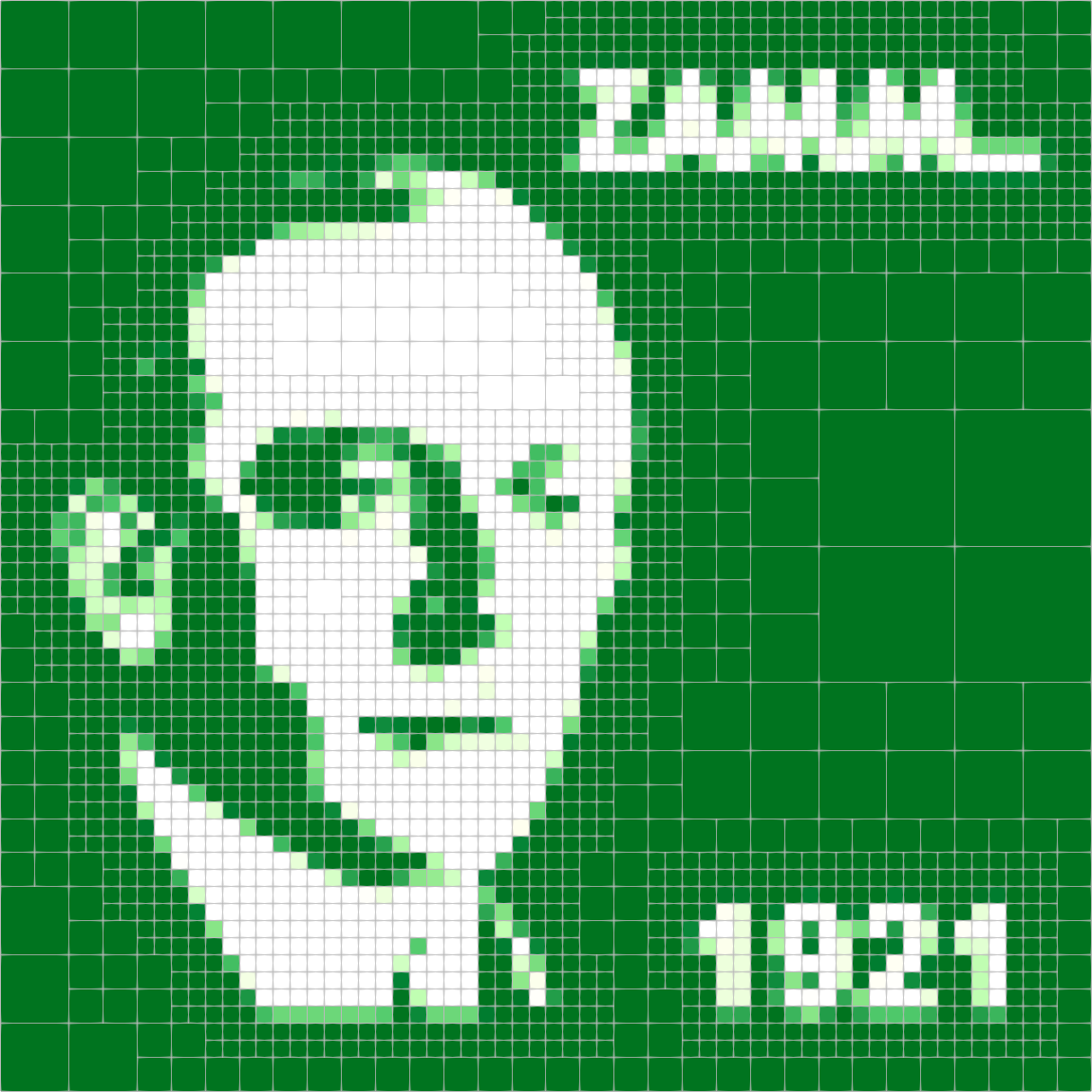}}
	\caption{\textbf{Adaptive mesh-coarsening} after uniform pixel coarsening along with new interphases in (b) and (c). The first entry in the subcaptions refers to the resolution, the second to the number of steps carried out in adaptive mesh coarsening.
 	\label{fig:Adaptive-meshcoarsening-vMises}}
\end{Figure} 

The diagrams in Figs. \ref{fig:RelErr-AdapComb-Interphase} and \ref{fig:RelErr-AdapComb-PhasePreserving} underpin that uniform resolution coarsening along with the same uniform discretization implies a discrete error jump, which is considerably larger than the error accumulating in several adaptive mesh coarsening steps at a fixed resolution. More precisely, adaptive mesh coarsening at discretization $h_{\square}$ is in general the better choice than choosing resolution and discretization $h_{\square}/2$, since the corresponding reduction of unknowns is at least as strong as the ndof-reduction of the pixel-coarsening step, the error increase throughout much smaller. 

\begin{Table}[htbp]
\begin{minipage}{16.5cm}  
\footnotesize
\centering
\renewcommand{\arraystretch}{1.2} 
\begin{tabular}{c r c c c c c c}
\hline
&  & \multicolumn{6}{c}{adaptive mesh coarsening} \\
\multicolumn{2}{r}{step no.} & 0  & 1  & 2  & 3  & 4 & 5   \\
\hline
&  & \multicolumn{6}{c}{additional interphases} \\
\hline
& ndof & $524\,288$ & $177\,298$ & $106\,272$ & $91\,398$ & $88\,668$ & $88\,264$   \\
& factor & $1.0000$ & $0.3382$ & $0.2027$ & $0.1743$ & $0.1691$ & $0.1684$ \\
\hline

%			& time (s) & $0\,000.0$ & $000.0$ & $000.0$ & $00.0$ & $00.0$ & $00.0$ \\
%			& speedup  & $1.0$ & $x.xx$ & $yy.yy$ & $zz.zz$ & $bbb.bb$ & $cc.cc$ \\
%			\hline
& $ e^{\epsilon}_{\text{mic}}$ & $16.5271$ & $16.7479$ & $16.9109$ & $17.1809$ & $17.5462$ & $18.0478$ \\
& $\bar e^{\epsilon \, h}_{\text{mic}}(\text{quad})$ & $6.2388$ & $6.8115$ & $7.2112$ & $7.8506$ & $8.6029$ & $9.4534$ \\
%& factor & $1.0000$ & $1.0918$ & $1.1559$ & $1.2583$ & $1.3789$ & $1.5152$ \\
& $\bar e^{\epsilon \, h}_{\text{mic}}(\text{ave})$ & $11.9762$ & $12.2845$ & $12.5108$ & $12.8904$ & $13.3627$ & $13.9266$ \\
%			\hline
& $e^{\epsilon \, h}_{\text{mic}}$ & $8.8569$ & $9.2594$ & $9.5493$ & $10.0171$ & $10.6284$ & $11.4299$ \\
%& factor & $1.0000$ & $1.0454$ & $1.0782$ & $1.1310$ & $1.2000$ & $1.2905$ \\
%			\hline
& $\Theta (\text{quad})$ & $0.7044$ & $0.7356$ & $0.7552$ & $0.7837$ & $0.8094$ & $0.8271$ \\
& $\Theta (\text{ave})$ & $1.3522$ & $1.3267$ & $1.3101$ & $1.2868$ & $1.2573$ & $1.2184$ \\
& $e^{\epsilon \, \square}_{\text{mic}}$ & \multicolumn{6}{c}{$14.4814$} \\[0mm]
\hline
&  & \multicolumn{6}{c}{phase-preserving} \\
\hline
& ndof & $524\,288$ & $172\,422$ & $100\,642$ & $85\,644$ & $82\,868$ & $82\,452$   \\
& factor & $1.0000$ & $0.3289$ & $0.1920$ & $0.1634$ & $0.1581$ & $0.1573$ \\
\hline
& $ e^{\epsilon}_{\text{mic}}$ & $23.1379$ & $23.3216$ & $23.4440$ & $23.6409$ & $23.9109$ & $24.2939$ \\
& $\bar e^{\epsilon \, h}_{\text{mic}}$ & $10.7376$ & $11.1378$ & $11.3956$ & $11.8141$ & $12.3389$ & $12.9726$ \\
% & factor & $1.0000$ & $1.0373$ & $1.0613$ & $1.1003$ & $1.1491$ & $1.2081$ \\
%			\hline
& $e^{\epsilon \, h}_{\text{mic}}$ & $12.6256$ & $12.9611$ & $13.1778$ & $13.5233$ & $13.9912$ & $14.6310$ \\
% & factor & $1.0000$ & $1.0266$ & $1.0437$ & $1.0711$ & $1.1082$ & $1.1588$ \\
%			\hline
& $\Theta$ & $0.8505$ & $0.8593$ & $0.8648$ & $0.8736$ & $0.8819$ & $0.8867$ \\
& $e^{\epsilon \, \square}_{\text{mic}}$ & \multicolumn{6}{c}{$21.4106$} \\[0mm]
\hline
& & \multicolumn{6}{c}{Reference solution for 2048$^2$ px, ndof$=8\,445\,856$.} \\
\hline
\end{tabular} 
\end{minipage}
\caption{\textbf{Adaptive mesh-coarsening}: for different adaptive mesh coarsening steps (1 to 5, 0 refers to the discretization with 512$^2$ px) the number of degrees of freedom (ndof) along with the ndof-factor compared to the uniform mesh, and the estimated errors for different methods are displayed. All error data in $10^{-4}$\,(MPa).
%{\color{red} @Andreas: Bitte erweitern wie in Tab.2.} {\color{blue}erledigt.}
	\label{tab:AdapCoarsening-Ndof-Error-speedup}} 
\end{Table}

\begin{Figure}[htbp]
\centering
\subfloat[$e^{\epsilon}_{mic}/||\bm u||_{A(\mathcal{B}_{\epsilon})}$]
{\resizebox{0.40\columnwidth}{!}{\includegraphics{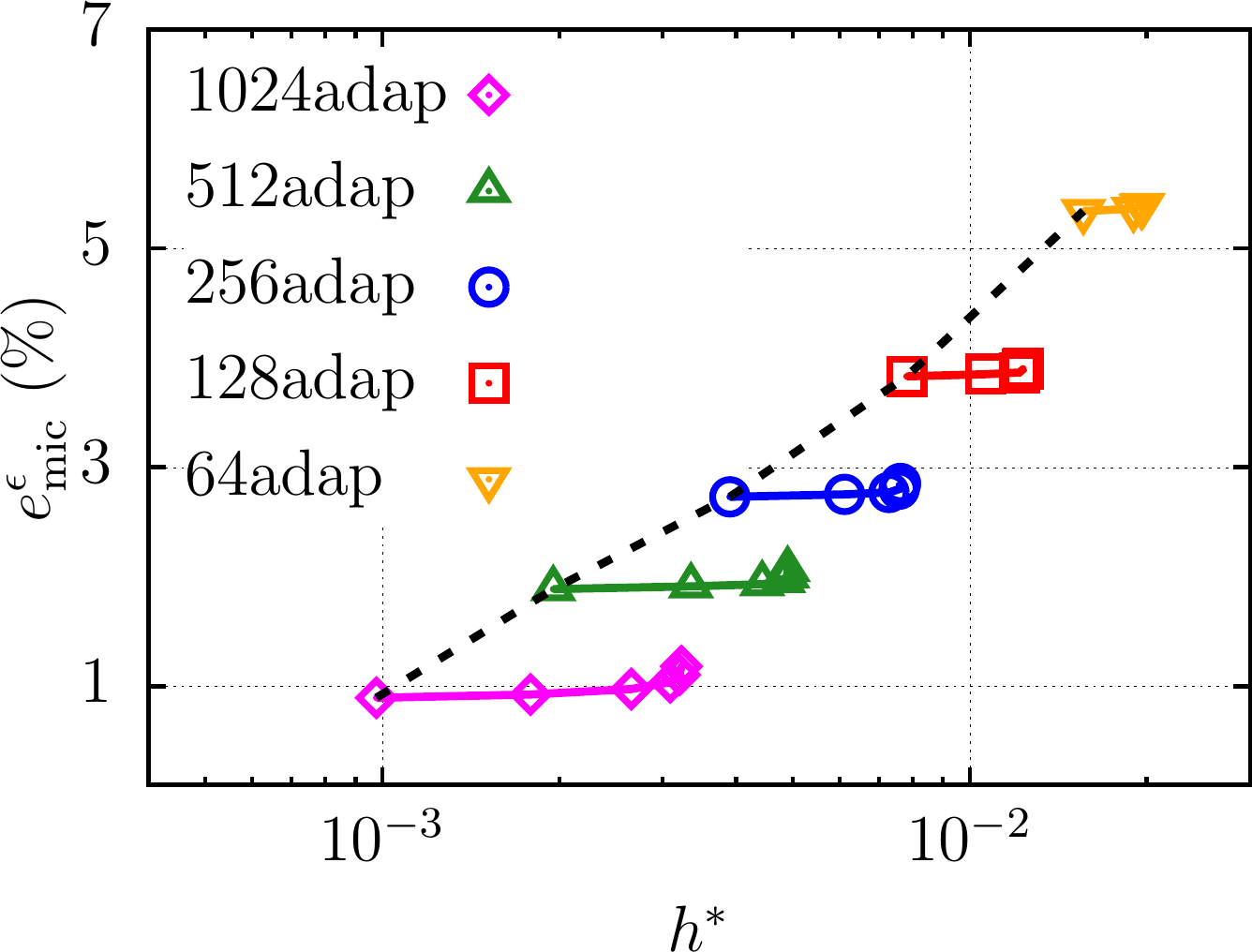}}}
\hspace*{0.02\linewidth}
	\subfloat[$e^{\epsilon\, h}_{mic}/||\bm u||_{A(\mathcal{B}_{\epsilon})}$]
{\resizebox{0.40\columnwidth}{!}{\includegraphics{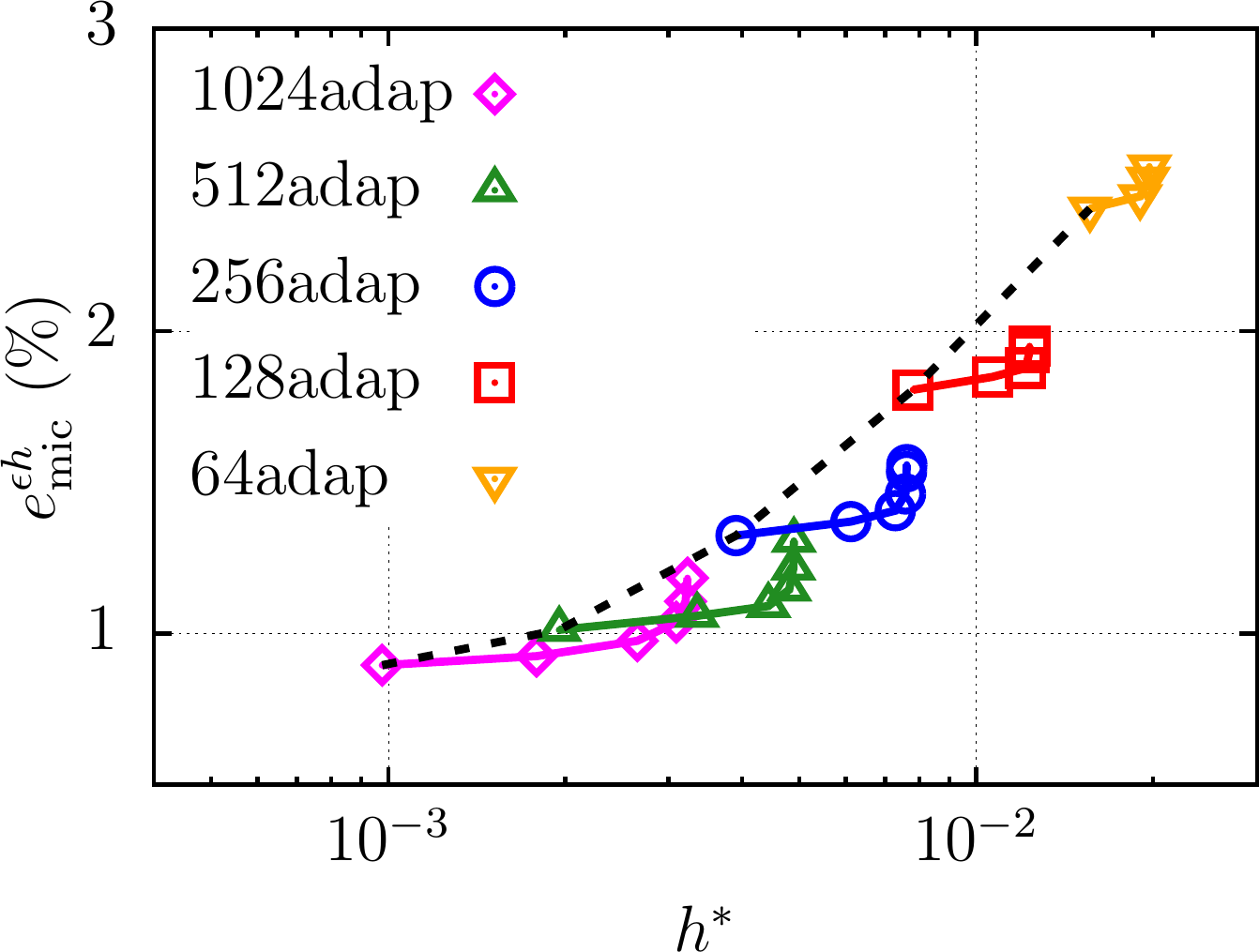}}}
\caption{\textbf{Adaptive mesh coarsening, additional interphases:} relative micro errors: (a) total error and (b) discretization error versus characteristic element length $h^{*}$. The dashed line marks the path of uniform pixel/mesh coarsening.}
\label{fig:RelErr-AdapComb-Interphase}
\end{Figure}

\begin{Figure}[htbp]
	\centering
    \subfloat[$e^{\epsilon}_{mic}/||\bm u||_{A(\mathcal{B}_{\epsilon})}$]
	{\resizebox{0.40\columnwidth}{!}{\includegraphics{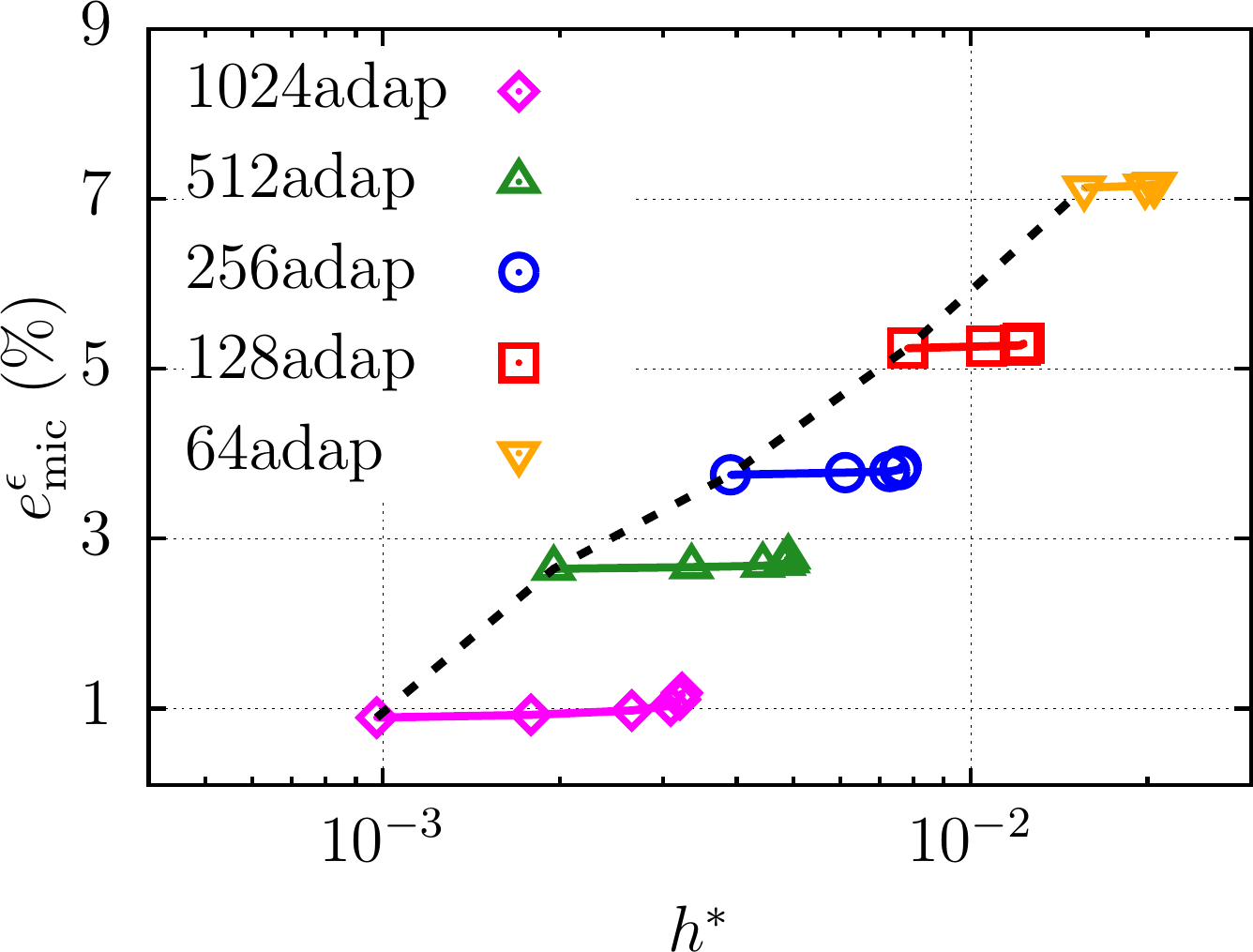}}}
	\hspace*{0.02\linewidth}
	\subfloat[$e^{\epsilon\, h}_{mic}/||\bm u||_{A(\mathcal{B}_{\epsilon})}$]
	{\resizebox{0.40\columnwidth}{!}{\includegraphics{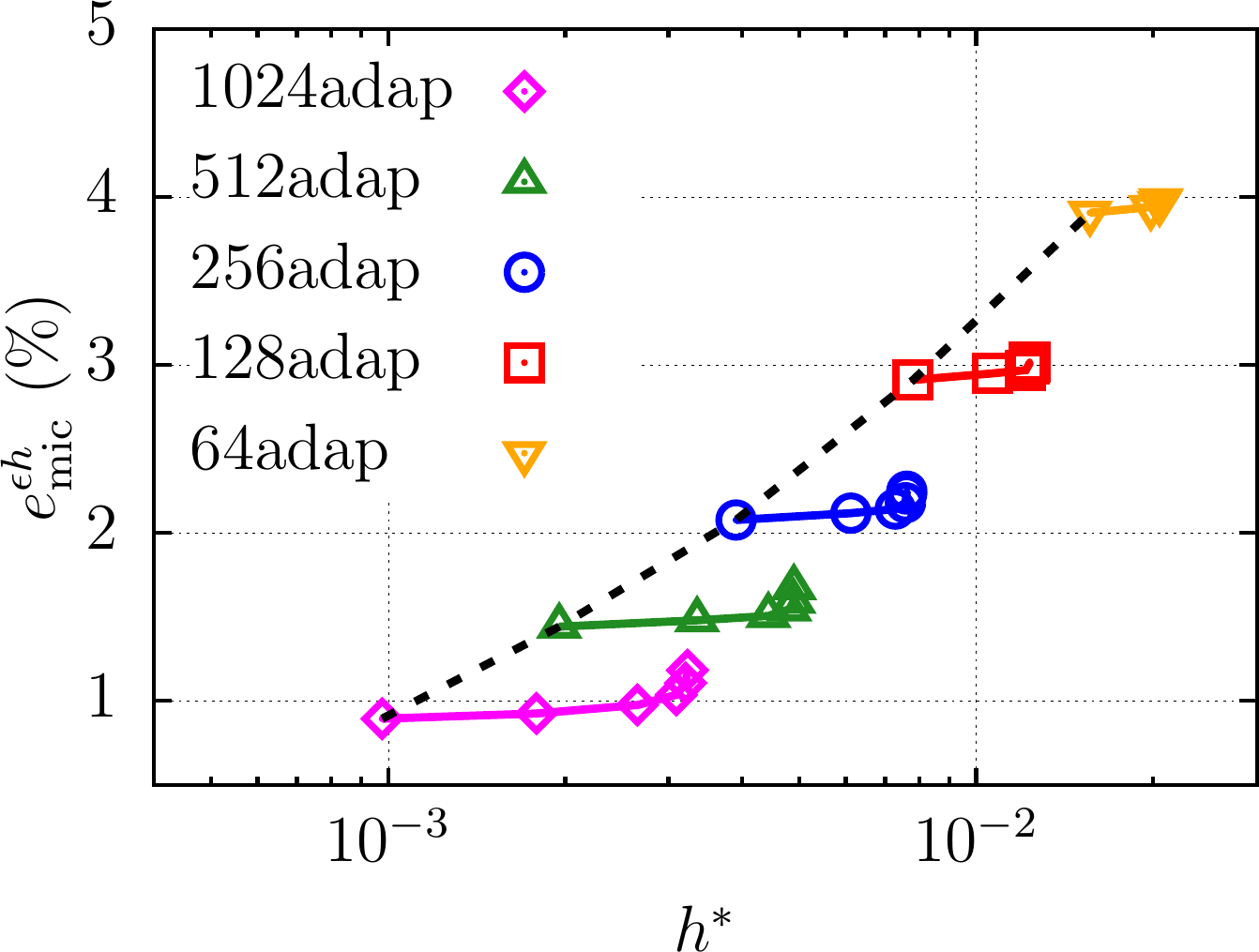}}}
	\caption{\textbf{Adaptive mesh coarsening, phase-preserving:} relative micro errors: (a) total error and (b) discretization error versus characteristic element length $h^{*}$. The dashed line marks the path of uniform pixel/mesh coarsening.}
	\label{fig:RelErr-AdapComb-PhasePreserving}
\end{Figure}

Notice that in Figs. \ref{fig:RelErr-AdapComb-PhasePreserving} and \ref{fig:RelErr-AdapComb-Interphase} a relative error is chosen defined as the absolute micro error over the total energy of the microdomain. In this representation the error is put into the perspective of the energy which is more descriptive than e.g. the percental error increase which was used in \cite{Fischer.2020}.
%, which is used in Tab.~\ref{tab:AdapCoarsening-Ndof-Error-speedup}.
 
\subsection{Comparison of strain}

\begin{Figure}[htbp]
	\centering
%1
	\subfloat[bulb][ndof=$2\,097\,152$ \\ \phantom{(a)} uni]
	{\includegraphics[width=0.32\linewidth, angle=0]{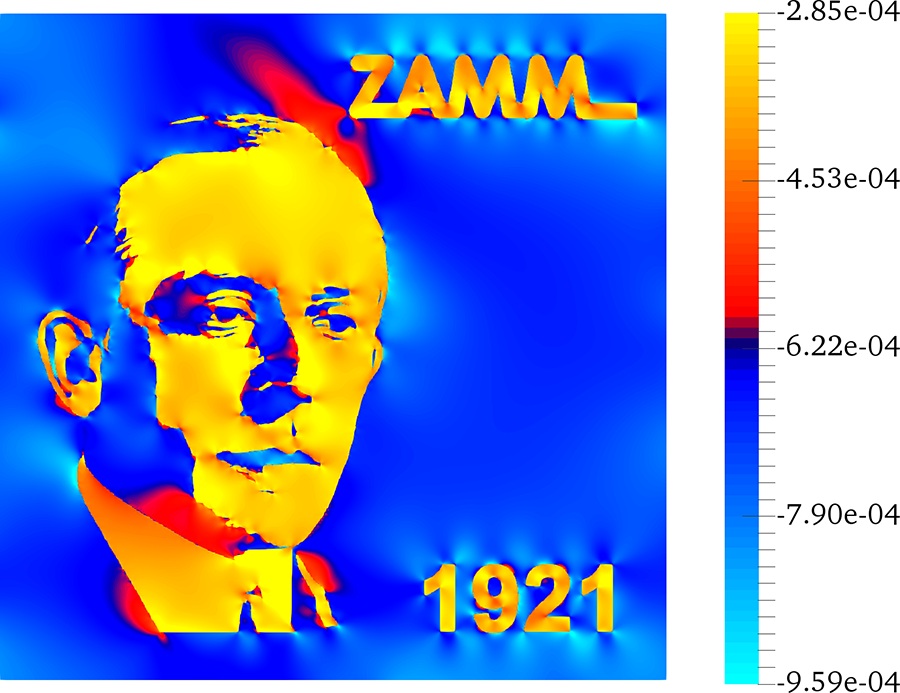}}
	\hfill
%3	
	\subfloat[bulb][ndof=$131\,072$ \\ \phantom{(a)} uni, preserv2]
	{\includegraphics[width=0.32\linewidth, angle=0]{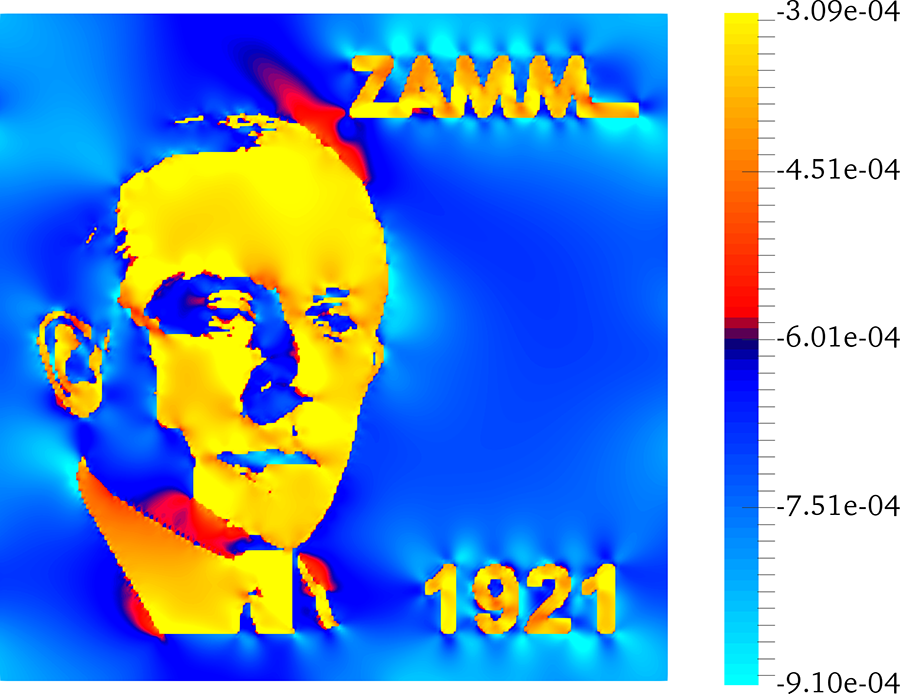}}
	\hfill
%5	
	\subfloat[bulb][ndof=$131\,072$ \\ \phantom{(a)} uni, interph2]
	{\includegraphics[width=0.32\linewidth, angle=0]{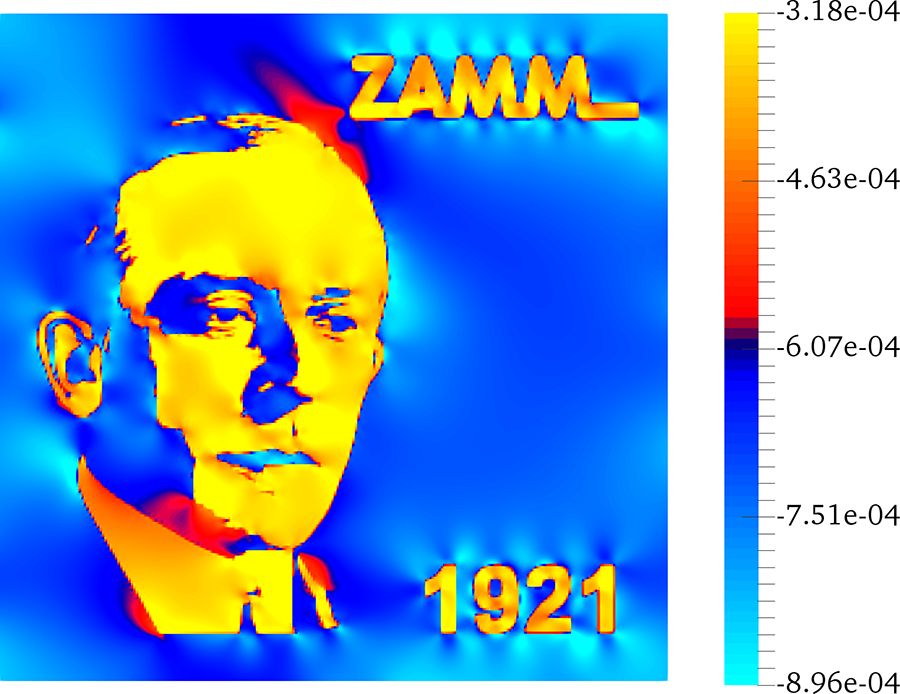}}
	\newline
%2	
	\subfloat[bulb][ndof=$209\,250$ \\ \phantom{(a)} adap3]
	{\includegraphics[width=0.32\linewidth, angle=0]{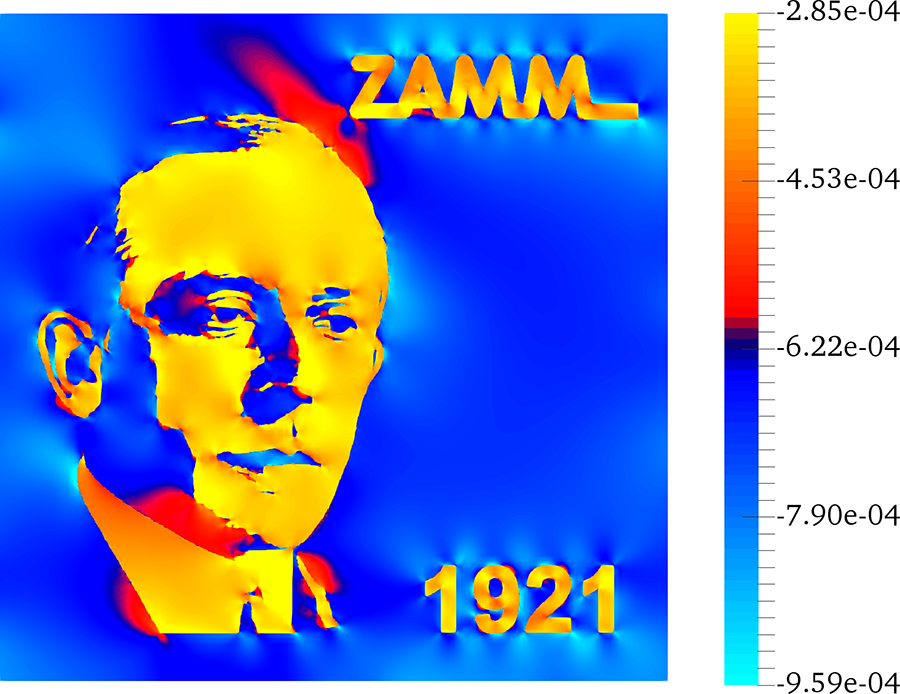}}
	\hfill
%4	
	\subfloat[bulb][ndof=$36\,886$ \\ \phantom{()} preserv2 \& adap2]
	{\includegraphics[width=0.32\linewidth, angle=0]{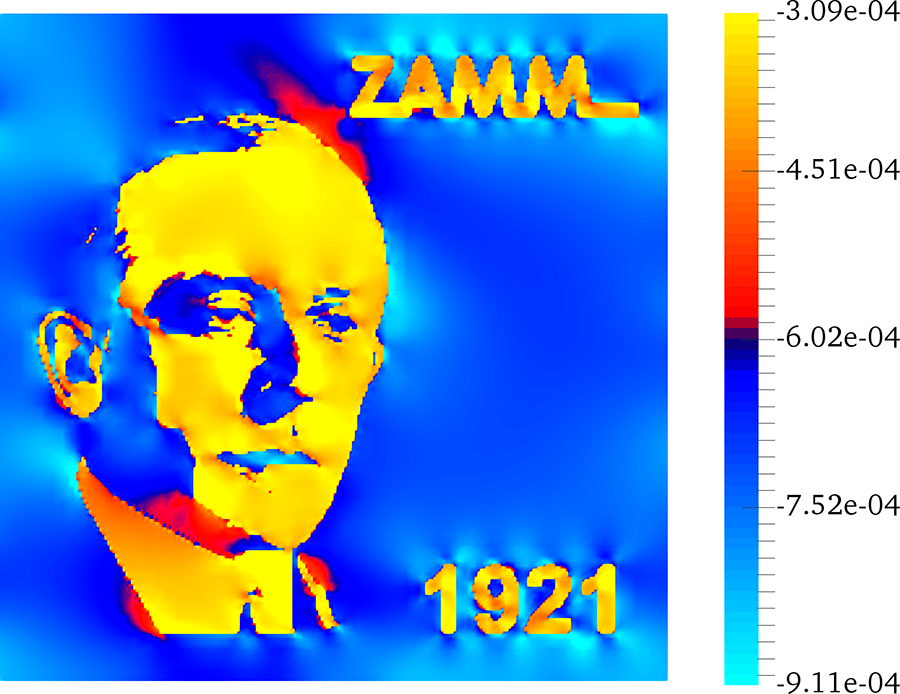}}
	\hfill
%6	
	\subfloat[bulb][ndof=$40\,322$ \\ \phantom{()} interph2 \& adap2]
	{\includegraphics[width=0.32\linewidth, angle=0]{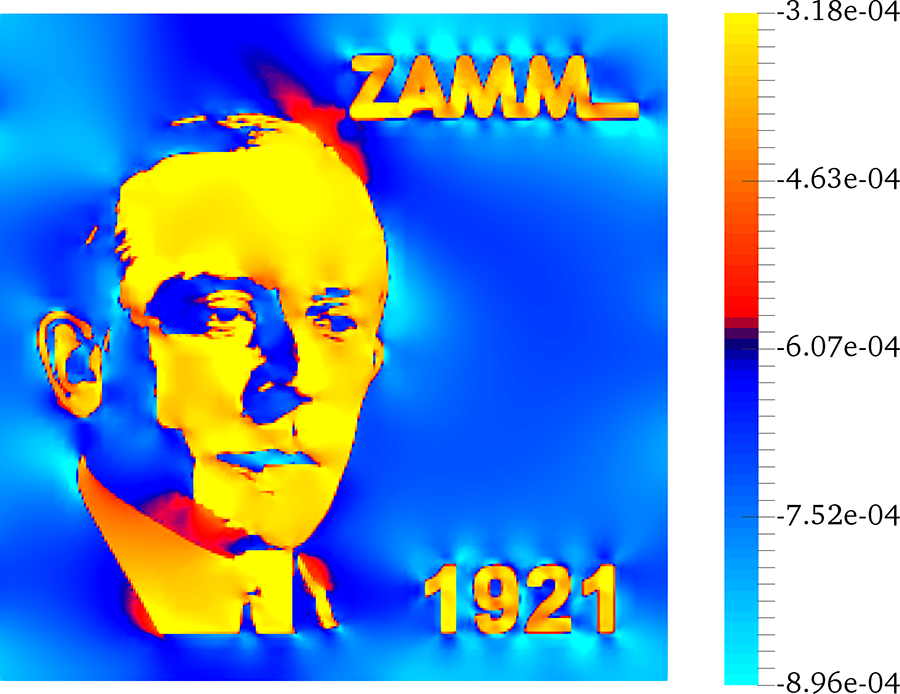}}
	
	\caption{\textbf{Shear strain distributions:} on microdomain at $x=2\,752.1$ mm, $y=552.1$ mm for (a) uniform mesh and its (d) three times adaptively coarsened counterpart, (b) phase-preserving, twicely uniform resolution coarsening and its (e) twicely adaptively coarsened counterpart, (c) twicely uniform resolution coarsening with interphases and its (f) twicely adaptively coarsened counterpart.}
	\label{fig:von-Mises_strain_plots}	
\end{Figure}

Figure \ref{fig:von-Mises_strain_plots} shows the distribution of shear strain at an exemplary microdomain for the initial uniform micro mesh and several coarsened meshes having an reduction of unknowns of at least 90\%. The best agreement --both qualitative and quantitative-- with the reference solution in Fig. \ref{fig:von-Mises_strain_plots} (a) is observed for the three times adaptively coarsened mesh in Fig. \ref{fig:von-Mises_strain_plots} (d) having the same pixel resolution as the reference. Compared to the other coarsened meshes this one exhibits the largest ndof due to the finely resolved interfaces, where coarsening is confined to the phase interior.

Nevertheless the solutions with uniformly coarsened meshes (phase-preserving in Fig. \ref{fig:von-Mises_strain_plots} (b) and with interphases in Fig. \ref{fig:von-Mises_strain_plots} (c)) also show a very good agreement with the initial uniform discretization. There are only minor deviations in the maximum and minimum values of the shear strain, the qualitative strain distribution is quite accurate. 

A further adaptive mesh coarsening on top of the already uniformly coarsened resolution/meshes (Fig. \ref{fig:von-Mises_strain_plots} (e) and (f)) renders an additional ndof reduction for very minor deviations in the resulting strain distributions. 

\subsection{Results and discussion}

The results can be summarized as follows:
\begin{enumerate}
	\item Pixel-coarsening along with new, artificial interphases, variant A and standard in bitmap image representation, introduces physical artefacts. They manifest in reduced interfacial stiffness mismatch and therefore in reduced maxima of stress and of their interfacial jumps. Error estimation along with the sound phase-distinction in stress computation introduces vanishing errors at interfaces, where they are truly largest. It is of all things the falsifying stress computation based on stress averaging, which reduces this effect. Hence the 'crime' in resolution coarsening can only be alleviated by an additional 'crime' in stress-computation. Although errors in the approximation of true micro errors are of the same magnitude as for variant B, errors in the effective elasticities are considerably larger. 
	\item Phase-preserving coarsening, variant B and the very standard in image segmentation, primarily relies on the knowledge of number and properties of individual phases. It maintains the phase contrast and therefore interfacial stress maxima and jumps, at least approximately, if the physically sound phase-distinction in stress computation for interface nodes is applied. 
	Another plus compared to variant A of coarsening is its accuracy in estimating the discretization error. The phase fraction ratio is not preserved, deviations typically increase the more, the coarser the resolution. 
	\item Adaptive mesh-coarsening --keeping fine resolution at interfaces, making mesh-coarsening the phases-- turns out to be a most effective means to achieve computational savings at very moderate error increases. The best strategy for achieving a target accuracy at minimal costs is to carry out resolution coarsening, which defines the accuracy level, and a consecutive adaptive mesh coarsening, which preserves this accuracy level.  
	\item The almost constant error in effective properties for micro resolutions 1024$^2$ px down to 64$^2$ px for the phase-preserving case indicates the trustworthiness in image acquisition of this microstructure. It is however restricted to macroscale results. 
	\item Microscale quantities are clearly much more sensitive to microscale resolution than macroscale, effective properties. 
%	{\color{blue} Der Satz ist etwas schwer zu verstehen so, finden Sie nicht?} Nö.
    A microstructure, coarse in resolution and discretization, which is likely rejected on the microscale for its large errors can show accurate macroscale properties as e.g. effective stiffness or maximal deflection of the macrostructure.    
	\item The augmented micro error analysis enables a detailed accuracy-efficiency balance. It quantifies the potential of resolution coarsening and adaptive mesh coarsening for a considerable efficiency gain at target accuracy. For the present microstructure in its initial, uniform 1024$^2$ px resolution, a reduction of unknowns to 0.6\% (from 2.1 mio. to 12k) can be achieved by a 128$^2$ px resolution along with adaptive mesh coarsening with the total micro error hardly above 5\%.
\end{enumerate}

With its newly introduced error analysis on the microscale the present work goes beyond the few existing investigations on the impact of image resolution on effective properties. They refer to diffusivity and conductivity of cathode catalyst layers \cite{Hutzenlaub.2013}, to the impact on global stress-strain curves, local damage initiation of foamed concrete \cite{Nguyen.2015} and on flow in porous media \cite{Shah.2016} and some applications in digital rock physics \cite{Berg.2017}.

% Lisa's statement in the MechMat-Article
%Overall, RVEs defined as such real microstructures usually require a complex discretization and thus high computational costs. This issue can be circumvented using statistically similar RVEs (SSRVEs) instead of conventional RVEs, cf. Balzani et al. (2014b), see also Balzani et al. (2014a) for two-dimensional results. These SSRVEs consist of an artificial inclusion morphology which is less complex than the real one but which still ensures similar morphological properties in a statistical sense. The method for the construction of SSRVEs is based on the that the material’s microstructure may be represented by a periodic microstructure composed of SSRVEs periodically put together. A further important assumption is that the macroscopic response is only dependent on the microstructure morphology provided that the behavior of the individual microscopic constituents is known.

\section{Conclusion and Outlook}
\label{sec:summary}

In this work we have introduced a rationale to assess microstructure images of finite resolution for use in computational solid mechanics. The analysis of multiphase solids distinguishes between a finite resolution modeling error and a discretization error. It is consistently embedded in the unified framework of errors in two-scale finite element methods for numerical homogenization (FE-HMM and FE$^2$).  

For the considered example the best choice is a combination of phase-preserving resolution coarsening with a consecutive adaptive mesh-coarsening; it results in considerable computational savings at controlled errors. We consider that combination as advantageous for other examples likewise, but the range of resolution and discretization at a target accuracy is case-dependent, its identification therefore left to an explicit error analysis. 

It is worth to note that for the first steps of uniform, phase-preserving pixel coarsening, the modeling error is approximately of the same magnitude as the discretization error. This suggests that the estimate of the discretization error could approximate the modeling error and thereby (an upper bound of) the total micro error. Since the modeling error is not directly accessible to error estimation and for that reason the total micro error neither, this aspect deserves further investigation. It provides an interesting link to the finite cell method (FCM) where the quadrature error (which corresponds to the modeling error in the present work) shall be bound to the magnitude of the discretization error. Only very recently error estimation for FCM has been introduced \cite{DiStolfo.2019}.
 
Worth for additional investigations is the influence of pixel size with respect to a characteristic microstructural length scale and the impact of phase contrast on errors. The gradual deviation from the true RVE characteristics by means of continuous resolution coarsening could be measured by statistical similarity measures thus making a link to the computed resolution errors.

The analysis of modeling and discretization errors in image data is similarly applicable to 3D problems following from reconstructed voxel-microstructures. Moreover, the concept opens the door to other applications, in the mechanics of solids to nonlinear models and methods such as hyperelasticity or inelastic constitutive laws in nonlinear computational homogenization. 

We consider a two-phase von-Mises-Prandtl-microstructure\footnote{Richard von Mises (1883--1953) and Ludwig Prandtl (1875--1953) founded the Gesellschaft f\"ur Angewandte Mathematik und Mechanik (GAMM) in 1922.}  subject to PBC. The rectangular unit cell of edge lengths $\epsilon_x=2$~mm, $\epsilon_y=1$~mm exhibits an initial uniform pixel resolution of 2048$\times$1024. The elastic properties coincide with the previous example. The matrix phase is purely elastic, the inclusion phase follows von-Mises elasto-plasticity \cite{Mises.1928}, a model not only frequently used for metal plasticity but also a source of inspiration for further original work till today \cite{Dyck.2020}. The initial yield stress is $y_0=0.1$~MPa, the constant modulus of isotropic hardening is $H=5$~MPa. For a uniaxial stretch in $x$-direction of $\varepsilon_{xx}= 0.088$ using 22 uniform load steps, Fig. \ref{fig:vonMisesPrandtlGAMM-eqvplstr} displays for two different resolutions along with their adaptively coarsened discretizations the results in terms of equivalent plastic strain as the quantity of interest. Image resolution has a strong effect on maximal values, adaptive coarsening maintaining the resolution however only a very minor one. This suggests that elasto-plastic analyses are more demanding in resolution for results of high fidelity than merely elastic analyses are. Results in depth are discussed elsewhere. 

\bigskip

{\bf Acknowledgements.} Bernhard Eidel acknowledges support by the Deutsche Forschungsgemeinschaft (DFG) within the Heisenberg program (grant no. EI 453/2-1). Simulations were performed with computing resources granted by RWTH Aachen University under project ID BUND0005. 

%\vfill
%\newpage

\begin{Figure}[htbp]
	\begin{minipage}{16.5cm}  
		\subfloat[512$\times$256 resolution, uniform mesh, ndof=1\,191\,192\,: equiv. pl. strain $\alpha$.]
		{\includegraphics[width=7.2cm,angle=0]{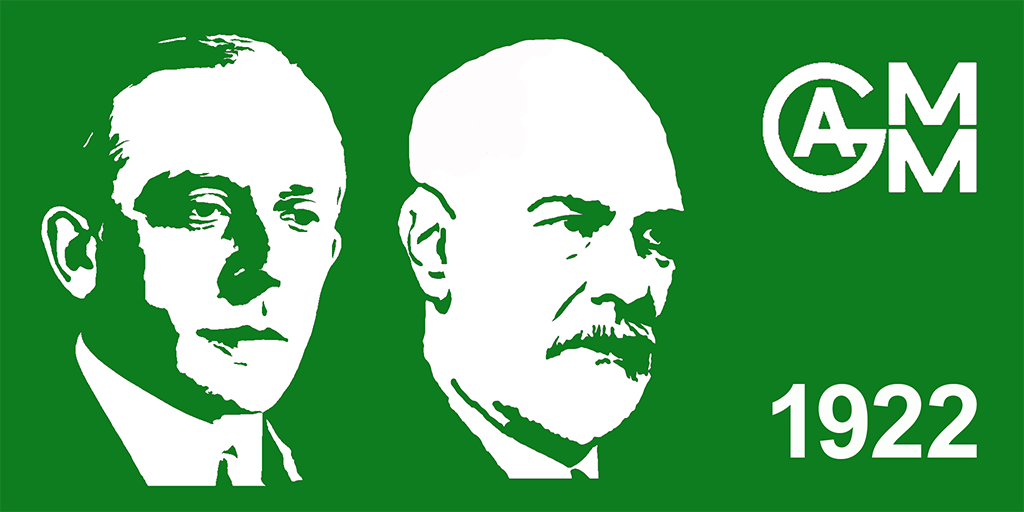}      
		 \includegraphics[width=8.5cm,angle=0]{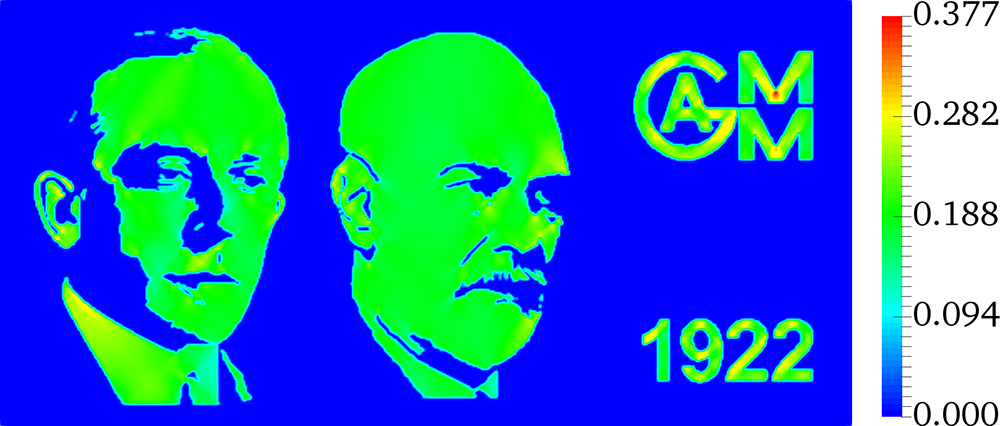}}
	    \\
	    \subfloat[512$\times$256 resolution, 5 adaptive coarsening steps, ndof=365\,400\,: equiv. pl. strain $\alpha$.]
		{\includegraphics[width=7.2cm,angle=0]{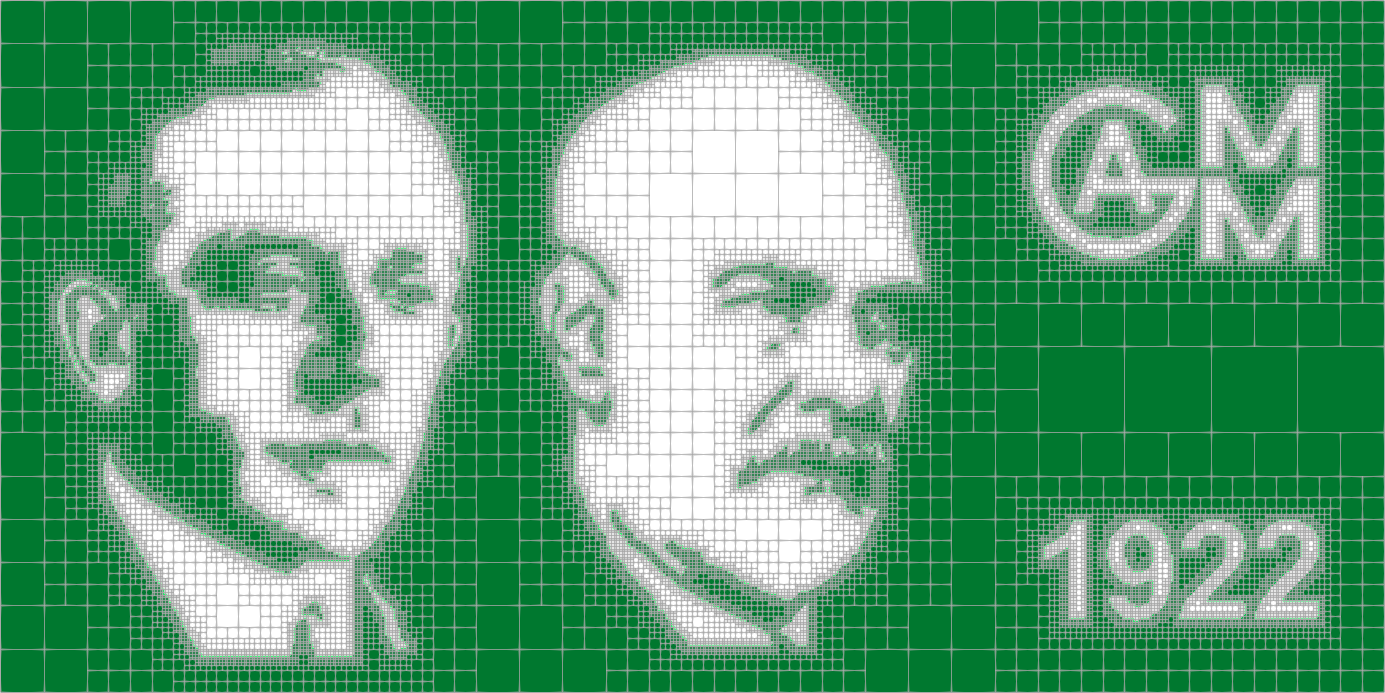} 
		 \includegraphics[width=8.5cm,angle=0]{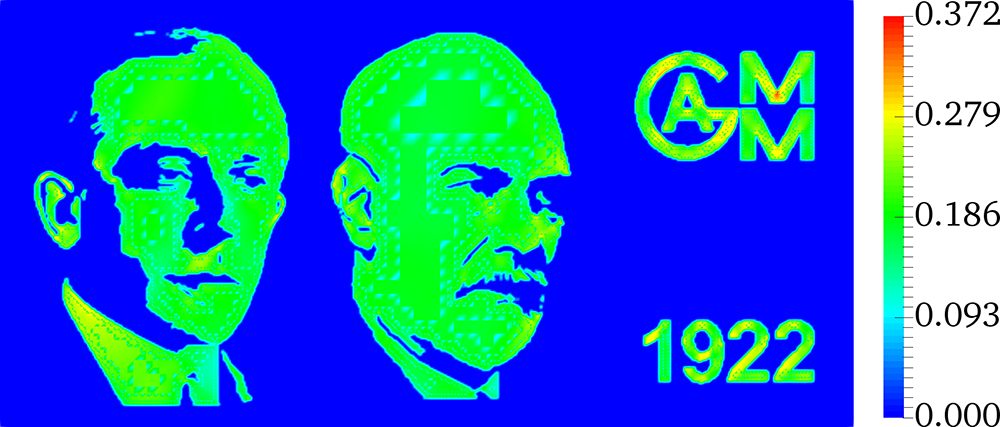}}
	    \\
	    \subfloat[128$\times$64 resolution, uniform mesh, ndof=76\,632\,: equiv. pl. strain $\alpha$.]
		{\includegraphics[width=7.2cm,angle=0]{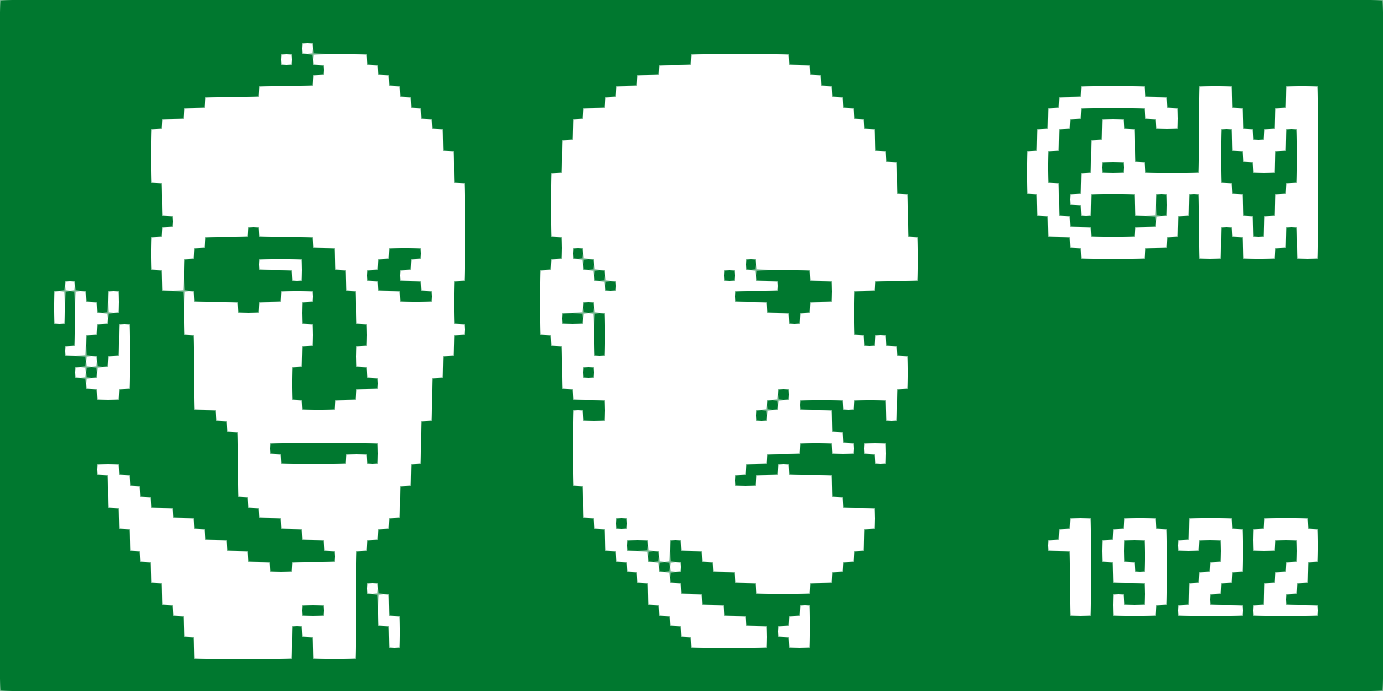}   
		 \includegraphics[width=8.5cm,angle=0]{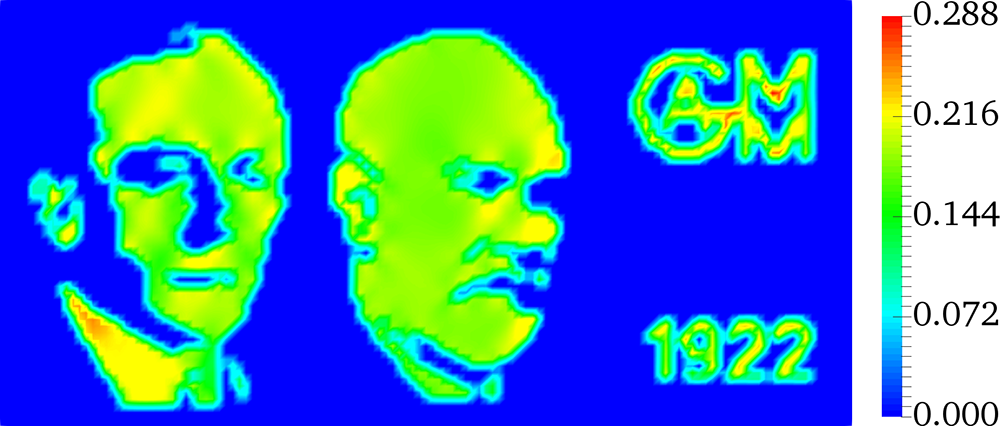}}
	    \\
	    \subfloat[128$\times$64 resolution, 3 adaptive coarsening steps, ndof=49\,206\,: equiv. pl. strain $\alpha$.]
		{\includegraphics[width=7.2cm,angle=0]{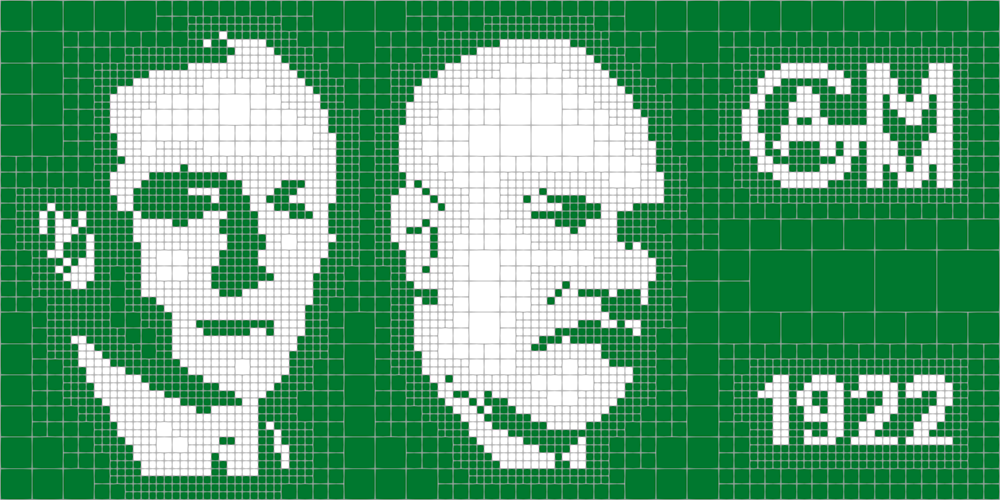}  
		 \includegraphics[width=8.5cm,angle=0]{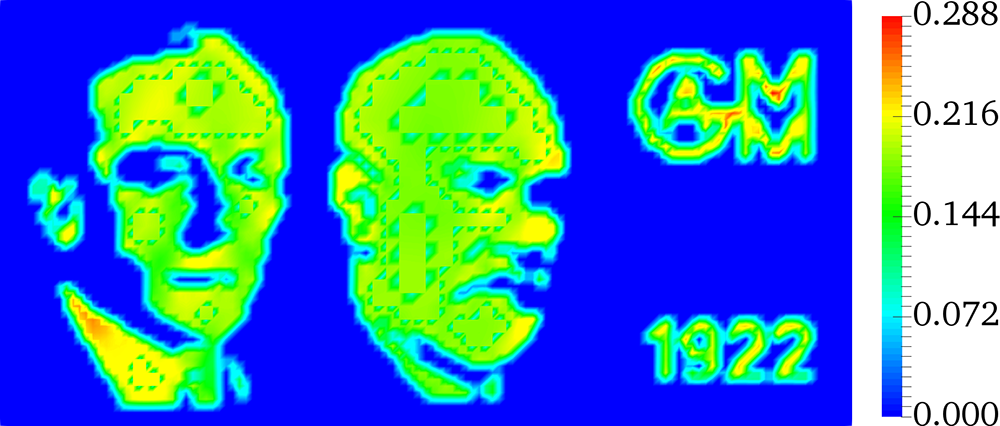}}
	\end{minipage}
	\caption{\textbf{GAMM-founders microstructure}: (a,c) different uniform, two-phase resolutions along with their (b,d) adaptively coarsened meshes result in (a--d, right column) different distributions of equivalent plastic strain.
    \label{fig:vonMisesPrandtlGAMM-eqvplstr}}
\end{Figure}
  
{\bf Declaration of Interest.} None.

%=========================================================================================

%\vfill
%\newpage

\begin{appendix}

\addcontentsline{toc}{section}{Appendix}
\renewcommand{\thesubsection}{\Alph{section}.\arabic{subsection}}
\renewcommand{\theequation}{\Alph{section}.\arabic{equation}}
\renewcommand{\thefigure}{\Alph{section}.\arabic{figure}}
\renewcommand{\thetable}{\Alph{section}.\arabic{table}}
\newcommand {\ssectapp}{
                        \setcounter{equation}{0}
                        \setcounter{figure}{0}
                        \setcounter{table}{0}
		                \subsection
                        }

\setcounter{equation}{0}

%=========================================================================================
%\vfill
%\newpage
%\input{appendix_AsymptoticExpansion}
%\input{appendix_NumTables}
%\input{appendix_Miscellaneous}

\end{appendix}

\bibliographystyle{abbrv} 
\bibliography{octree}
%\bibliographystyle{plaindin}
%\bibliographystyle{plaindin_shortname2}
%\bibliographystyle{elsarticle-num-names}
% ------- bib-datei --------------
%\bibliography{js_master_2009}
%\begin{appendix}
%\include{appendixA}
%\end{appendix}
\end{document}